\pdfoutput=1
\pdfoutput=1
\documentclass[final,1p,times]{elsarticle}

\usepackage{epsfig}
\usepackage{float}

\usepackage{amssymb}

\usepackage{amsmath, amssymb}

\usepackage{lineno}

\newcommand{\pt}{\,\partial_t\,}
\newcommand{\px}{\,\partial_x\,}
\newcommand{\py}{\,\partial_y\,}
\newcommand{\pz}{\,\partial_z\,}

\newcommand{\F}{\,\mathcal F\,}
\newcommand{\C}{\,\mathcal C\,}
\newcommand{\G}{\,\mathcal G\,}
\renewcommand{\P}{\,\mathcal P\,}
\renewcommand{\H}{\,\mathcal H\,}
\renewcommand{\S}{\,\mathcal S\,}
\newcommand{\dt}{\Delta t }
\newcommand{\dx}{\Delta x }
\newcommand{\dy}{\Delta y }
\newcommand{\dz}{\Delta z}

\journal{}

\begin{document}

\begin{frontmatter}



\title{Modelling and numerical approximation of a 2.5D set of equations for mesoscale atmospheric processes}


\author[dm]{Dante Kalise}
\address[dm]{Dipartimento di Matematica, Universit\`a di Roma ``La Sapienza'', P. Aldo Moro 2, 00185 Roma,
Italy.}
\ead{kalise@mat.uniroma1.it}
\author[storm]{Ivar Lie}
\address[storm]{Research and Development Department, StormGeo AS, Universitetsgata 8, 0164 Oslo,  Norway}
\ead{ivli@online.no}
\begin{abstract}
The set of 3D inviscid primitive equations for the atmosphere is dimensionally reduced by a Discontinuous Galerkin discretization in one horizontal direction. The resulting model is a 2D system of balance laws where with a source term depending on the layering procedure and the choice of coupling fluxes, which is established in terms of upwind considerations. The ``2.5D'' system is discretized via a WENO-TVD scheme based in a flux limiter centered approach. We study four tests cases related to atmospheric phenomena to analyze the physical validity of the model.
\end{abstract}

\begin{keyword}
primitive equations \sep layering \sep discontinuous Galerkin \sep upwind flux \sep WENO-TVD schemes \sep test cases for dynamical cores



\end{keyword}

\end{frontmatter}



\section{Introduction}

The so-called primitive equations used to model atmospheric and ocean dynamics
are not well-posed for any reasonable set of boundary conditions.
This important, and relatively old result, was obtained by Oliger and
Sundstr\"{o}m in \cite{Oliger}. The problem of formulating a well-posed set of
 ``primitive'' equations is of obvious importance both from the mathematical
 and numerical side, and has recently attracted substantial research activity.
One example of such research was undertaken in \cite{temam}, where the authors have
formulated a type of dimensionally-reduced equations, called 2.5D equations, where
well-posedness was analyzed. However, this paper uses linearized primitive equations and the layering procedure was
performed via a Continuous Galerkin approach using orthogonal piecewise linear functions. This imposes some limitations regarding the linear character of the system and the number of elements to be considered in the expansion.
What we want to communicate here is an alternative derivation of a 2.5D set of full {\em nonlinear} equations and its approximation by a high-order finite volume scheme. We consider the analysis of well-posedness for these equations to be a far too
difficult task. Instead we perform numerical experiments with well-known test
cases normally used to test dynamical cores of atmospheric models \cite{giraldo,skamarock}.
Good results from such experiments are a strong indication that our formulation
of a 2.5D set of equations is well-posed.
The dimensional reduction procedure that we propose is based on a Discontinuous Galerkin approach (DG) \cite{cockburn}; this class of schemes have proved to be successful in the resolution 2D atmospheric models \cite{giraldo}. We use the same basic principles to divide the domain into 2D layers where the equations are locally approximated and coupled by suitable fluxes. This approach leads to a system of 2D balance laws, which is solved via a WENO-TVD scheme \cite{titarevweno}, using a flux limiter centered approach \cite{torobillett} for space discretization and a Runge-Kutta TVD scheme for time discretization.
Since the primitive equations are not well-posed, one may wonder what approach
is taken in the many atmospheric and oceanic models that exist.
The answer is that a pragmatic, and hence not strictly rigorous, approach is taken.
In practice this means that one imposes boundary conditions that are not boundary
conditions in the ordinary sense. An example is the use of relaxation zones near
the boundary, where the solution in the inner of the domain is relaxed towards a
certain boundary value. While this is stable numerically (in most cases), it can
hardly be said to be satisfactory in the physical sense.
And not in the mathematical sense either.

The remaining part of the paper is organized as follows: In section 2 we present the
full 3D governing equations in conservative form for a dry atmosphere, and its
''semi-discretization'' in a spatial dimension so the equations are formulated
in 2.5D as explained above.
The finite volume-based numerical scheme is presented in detail in section 3, and
in section 4 we present a number of numerical experiments using both well-known
tests for atmospheric model dynamical cores, as well as new tests.  The results
are discussed for each test, and we summarize the results along with directions
for further work in section 5.

\section{A dimensionally-reduced system for atmospheric dynamics}
This section is concerned with the derivation of a dimensionally-reduced set of equations describing the atmosphere. We first present the original three-dimensional system of primitive equations which is then layered via a DG ansatz in order to obtain a 2.5D model.

\subsection{A three-dimensional system of balance laws}

Our starting point is the set of equations modelling the evolution in time of a dry  atmosphere.
We impose conservation of mass, momentum and energy, and consider effects of rotation
and gravity, together with neglecting friction; this leads to a set of 3D inviscid primitive equations cast in conservative form
\begin{equation}\label{euler3d}
\pt Q + \px \F + \py \G +\pz \H = \S,
\end{equation}
where
\begin{equation}
Q=\left[\begin{array}{c}\rho\\ \rho u\\ \rho v\\ \rho w\\ \rho\theta\end{array}\right],
\F=\left[\begin{array}{c}\rho u\\ \rho u^2 +\P\\ \rho uv\\ \rho uw\\ \rho u \theta\end{array}\right],
\G=\left[\begin{array}{c}\rho v\\ \rho uv\\ \rho v^2+\P\\ \rho vw\\ \rho v \theta\end{array}\right],
\end{equation}
\begin{equation}
\H=\left[\begin{array}{c}\rho w\\ \rho wv\\ \rho wv\\ \rho w^2+\P\\ \rho w \theta\end{array}\right],
\S=\left[\begin{array}{c}0\\ f \rho v\\ -f \rho u\\ -\rho g\\ 0\end{array}\right].
\end{equation}
The dynamics are expressed in terms of the density of the fluid $\rho$, the 3D velocity field $(u,v,w)$, the pressure $\P$ and the potential temperature $\theta$, which relates to the thermodynamic temperature $\mathcal{T}$ via\begin{equation}
\theta\equiv \mathcal{T}\left(\frac{\P}{\P_0}\right)^{-R_d/c_p}.
\end{equation}
The  system is closed by the equation of state for an ideal gas,
\begin{equation}
\P=C_0(\rho\theta)^{\gamma},\qquad C_0=\frac{R_d^\gamma}{\P_0^{R_d/c_v}}.
\end{equation}
The model parameters are: the Coriolis parameter $f=1\times10^{-4}[s^{-1}]$,
the gravitational acceleration $g=9.81 [ms^{-2}]$,
the atmospheric pressure at sea level $\P_0=10^5 [Pa]$,
the gas constant for dry air $R_d=287 [JK^{-1}kg^{-1}]$,
the specific heat of dry air at constant pressure and volume $C_p=1004 [JK^{-1}kg^{-1}]$,
the specific heat of dry air at constant volume $C_v=717 [JK^{-1}kg^{-1}]$
and its ratio $\gamma=C_p/C_v=1.4$. Moreover, by defining the Exner pressure
\begin{equation}
\pi\equiv\left(\frac{\P}{\P_0}\right)^{R/c_p},
\end{equation}
the expression for the total energy of the system is given by,
\begin{equation}
   e=c_v\theta\pi+\frac12(u^2+v^2+w^2)+gz.
\end{equation}
The system (\ref{euler3d}) governs the time evolution of the conserved variables in $Q$ by the interaction between the physical fluxes $\F,\G,\H$ and the source term $S$. A more detailed representation of the atmospheric dynamics could include an equation for
 conservation of moisture that will affect the thermodynamics and hence the dynamics
(expressed by  the fluxes). Physical parameterizations in a weather model, like surface parameterizations and
condensation processes, together with coordinate transformations, like the extensively used terrain-following coordinates will generate terms appearing in the source $\S$.
In the next subsection we present a dimensionally-reduced version of the original 3D set of equations where we characterize the behavior
of the variables in one spatial dimension, in a way such that the action of its associated flux is transferred to the source term.
Although this can be seen as a step in the numerical discretization of the original set of 3D equations, our study will be focused on the physical significance of the dimensionally-reduced (or 2.5D) set of equations as a system
of balance laws on its own.

\subsection{Towards a 2.5D model}

In \cite{temam}, the authors obtained a 2.5D model from the linearized primitive equations via a Continuous Galerkin expansion of the variables in the $y$-direction. It has to be emphasized that the derivation made considerable use of the linearity of the system, and that the basis functions were a set of orthogonal piecewise linear functions; removing any of those aspects leads to a model whose complexity is not clearly related those of the underlying physics. On the other hand, there is an obvious limitation related to the trade-off between orthogonality and locality of the basis functions; the aforementioned article considered two elements with a basis of three orthogonal piecewise functions, which is just a slight modification of the classical hat functions. An increment in the number of elements, together with the orthogonality condition, will produce a set of basis function with increasing support, significantly different from the analogous set of hat functions. The approach that we present pays a cost by removing the continuity requirement on the solution (in the classical Galerkin framework), and switches to a DG formulation that allows, in a straightforward and tractable manner, to deal with nonlinearity and locality at the same time.
We start by defining our set of equations over
$\,\Omega =[0,\,L^x]\times[0,\,L^y]\times[0,\,L^z]\times[0,\,T_f]$.
Initial and boundary conditions  will be specified later.
We consider a uniform partition in the y-direction into $N_y$ elements, i.e.
\begin{equation}
Y=[0,\,L^y]=\bigcup_{i=1}^{N_y}Y_i,\qquad Y_i=[y_{i-1/2},\,y_{i+1/2}]
    \qquad y_{i+1/2}-y_{i-1/2}=\dy,\;\forall i.
\end{equation}
Next, in every element we consider a local classical Galerkin ansatz by multiplying eq.
(\ref{euler3d}) by a test function $V$ and integrating with respect to $y$, leading to
\begin{equation}
\int_{y_{i-1/2}}^{y_{i+1/2}} V^T\{\pt Q + \px \F + \py \G +\pz \H - \S\}\,dy=0.
\end{equation}
Mapping every element into the canonical element $(-1,1)$ by means of
\begin{equation}
y=\frac{1-\xi}{2}\,y_{i-1/2}+\frac{1+\xi}{2}\,y_{i+1/2},
\end{equation}
and, as usual, performing integration by parts, leads to
\begin{equation}
\frac{\dy}{2}\int_{-1}^1V^T\{\pt Q + \px \F +\pz \H - \S\}\,d\xi+V^T\G|_{-1}^1=
          \int_{-1}^1\,\partial_{\xi}\,V^T\G\,d\xi.
\end{equation}
We consider a local set of basis functions and expand $Q$ restricted to $Y_i$:
\begin{equation}
Q_i(x,\xi,z,t)=\sum_{j=0}^{N_{i}}q_{ij}(x,z,t)\phi_j(\xi).
\end{equation}
The choice of the set of basis functions is an open question and will depend on
the application and the complexity of the computations.
In this particular case, in order to be able
to carry on in a tractable manner the calculations for a nonlinear model,
it is important that we chose a set of orthogonal functions, and therefore selecting
a set of orthogonal polynomials (as the Legendre polynomials for instance) will allow us
to preserve a reasonable model but also to include high-order approximations.
Throughout this article however, we will restrict ourselves to the most basic case
in this framework, i.e., to consider Legendre polynomials up to degree 0,
which is nothing but to consider a piecewise constant approximation of the variables
inside every element. In \cite{klt}, the piecewise linear extension of this procedure
was performed in a similar framework, and poses no additional difficulties except
for the increase in the computational cost associated with the amount of local expansion coefficients.
When we consider 2 elements in the $y$-direction, $Y_1=[0,\,0.5\,L^y]$ and
$Y_2=[0.5\,L^y,\,L^y]$ and the first Legendre polynomial $\phi_0\equiv 1$,
the 2.5D system of equations reads,
\begin{align}
\pt Q_1+\px\F(Q_1)+\pz \H(Q_1)&=\dy\{\G(y=0.5\,L^y)-\G(y=0)\}+\S(Q_1),\\
\pt Q_2+\px\F(Q_2)+\pz \H(Q_2)&=\dy\{\G(y=L^y)-\G(y=0.5\,L^y)\}+\S(Q_2).
\end{align}
At this point we must close the system by giving computable expressions
for the evaluation of $\G$ at the boundaries and the intermediate value.
We first address the problem of the intermediate value, and then we expand this
procedure to the boundary fluxes.
Following the same ideas of the finite volume framework, we seek to express
the intermediate state as a function of adjacent cells, i.e.
\begin{equation}
\G(y=0.5\,L^y)\approx \,G(Q_1,Q_2).
\end{equation}
As it is extensively described in the literature, the two main choices for such procedure
are centered (averaged) values, or states obtained via upwind considerations.
The typical centered fluxes can be interpreted as particular averaging operators
in space and time, which does not fit in a proper manner in our scheme, as
time integration will be considered only for the resulting 2.5D model;
therefore we opt for upwind fluxes to determine the coupling at the interface.
We recall that the flux $\G$ in the Euler equations has the homogeneity property
\begin{equation}
\G(Q)=\C(Q)Q,
\end{equation}
with
\begin{equation}
\C(Q)=\left[\begin{array}{ccccc}
0 & 0 & 1 & 0 & 0\\
-uv & v & u & 0 & 0\\
-v^2& 0 & 2v & 0 & \frac{c_s^2}{\gamma\theta}\\
-wv & 0 & w & v & 0\\
-\theta v & 0 & \theta & 0 & v\end{array}\right] ,
   \quad c_s^2=\partial_{\rho}\P=C_0(\rho\theta)^{\gamma-1}\gamma\theta,
\end{equation}
where $c_s$ the speed of sound in the fluid. This matrix admits a representation of the form
\begin{equation}
\C(Q)=R\,\Lambda\,R^{-1},
\end{equation}
with
\begin{equation}
R=\left[\begin{array}{ccccc}
1 & 1 & 0 & 0 & 1\\
u & 0 & 1 & 0 & u\\
v-\frac{c_s}{\sqrt{\gamma}}& v & 0 & 0 & v+\frac{c_s}{\sqrt{\gamma}}\\
w & 0 & 0 & 1 &w\\
\theta & 0 & 0 & 0 & \theta\end{array}\right],
\Lambda=\left[\begin{array}{ccccc} v-\frac{c_s}{\sqrt{\gamma}}&0&0&0&0\\
0&v&0&0&0\\
0&0&v&0&0\\
0&0&0&v&0\\
0&0&0&0&v+\frac{c_s}{\sqrt{\gamma}}\end{array}\right].
\end{equation}
We decompose $\Lambda$ into positive and negative parts, by defining
\begin{equation}
v^+=\max(v,0)\,,\quad v^-=\min(v,0)\,,
\end{equation}
and we obtain
\begin{equation}
\Lambda^+=\left[\begin{array}{ccccc} 0&0&0&0&0\\
0&v^+&0&0&0\\
0&0&v^+&0&0\\
0&0&0&v^+&0\\
0&0&0&0&v^++\frac{c_s}{\sqrt{\gamma}}\end{array}\right],
\Lambda^-=\left[\begin{array}{ccccc} v^--\frac{c_s}{\sqrt{\gamma}}&0&0&0&0\\
0&v^-&0&0&0\\
0&0&v^-&0&0\\
0&0&0&v^-&0\\
0&0&0&0&0\end{array}\right].
\end{equation}
Note that in the last expression we have assumed a low Mach number for the fluid
($\,|v|<<c_s\,$), as expected in atmospheric dynamics. Thus, the intermediate state is defined as
\begin{equation}
\G(y=0.5\,L^y)\approx \,G(Q_1,Q_2)=\G^+(Q_1)+\G^-(Q_2),
\end{equation}
where
\begin{eqnarray}
\G^+(Q_1)&=&R\Lambda^+\,R^{-1}Q_1=\frac12\left(v+\frac{c_s}{\sqrt{\gamma}}\right)\left[\begin{array}{c}
\rho\\
\rho u\\
\rho \left(v+\frac{c_s}{\sqrt{\gamma}}\right)\\
\rho w\\
\rho \theta\end{array}\right]\,,\\
\G^-(Q_2)&=&R\Lambda^-\,R^{-1}Q_2=\frac12\left(v-\frac{c_s}{\sqrt{\gamma}}\right)\left[\begin{array}{c}
\rho\\
\rho u\\
\rho \left(v-\frac{c_s}{\sqrt{\gamma}}\right)\\
\rho w\\
\rho \theta\end{array}\right].
\end{eqnarray}
The same argument can be applied at the boundaries; the tests performed in section
\ref{tests} considered solid wall boundary conditions, which were implemented
by defining ghost cells values $Q_0$ and $Q_3$ where,
\begin{eqnarray}
\rho_0&=&\rho_1\,,\,u_0=u_1\,,\,v_0=-v_1\,,\,w_0=w_1\,,\,\theta_0=\theta_1,\\
\rho_3&=&\rho_2\,,\,u_3=u_2\,,\,v_3=-v_2\,,\,w_3=w_2\,,\,\theta_3=\theta_2.
\end{eqnarray}
This final definition allows us to write a direct expression for the 2.5D model:
\begin{eqnarray}\label{sbl3}
\pt Q_1+\px\F(Q_1)+\pz \H(Q_1)&=&\S_1(Q_1,Q_2),\label{sbl31}\\
\pt Q_2+\px\F(Q_2)+\pz \H(Q_2)&=&\S_2(Q_1,Q_2),\label{sbl32}
\end{eqnarray}
\begin{eqnarray}
\S_1(Q_1,Q_2)&=&\dy\{\G^+(Q_1)+\G^-(Q_2)-\G^+(Q_0)-\G^-(Q_1)\}+\S(Q_1),\\
\S_2(Q_1,Q_2)&=&\dy\{\G^+(Q_2)+\G^-(Q_3)-\G^+(Q_1)-\G^-(Q_2)\}+\S(Q_2).
\end{eqnarray}
As we previously pointed out, the resulting dimensionally-reduced model maintains
the same structure for the physical fluxes in the $x$ and $z$-directions,
but has transformed the flux $\G$ into a ``reactive'' source term that generates
a coupling between the layers. The next section deals with the numerical approximation
of this coupled system of balance laws.

\section{The numerical scheme}

In this section we present a finite volume scheme for system of balance laws derived
in the previous section. For the sake of simplicity, we develop the scheme
for a single set of equations of the form,
\begin{equation}
\pt Q+\px \F(Q)+ \pz \H(Q)=\S(Q),
\end{equation}
being eqns. (\ref{sbl31})-(\ref{sbl32}) a particular case with the augmented state $Q=[Q_1,\, Q_2]^T$.
We first indicate that our strategy will be based in a splitting scheme, as it suggested
in \cite{TOROCHEN} given the flux choice that we will make. Thus, we will first establish a numerical scheme for the system of conservation laws,
\begin{equation}\label{scl3}
\pt Q+\px \F(Q)+ \pz \H(Q)=0,
\end{equation}
to be combined with a procedure for the resolution of the source term dynamics
\begin{equation}\label{st3}
\pt Q=\S(Q).
\end{equation}
In order to approximate eq. (\ref{scl3}), we begin by meshing the spatial domain
$\Omega_{x,z}$ into uniform control volumes
$\Omega_{i,j}=[x_{i-1/2},\,x_{i+1/2}]\times [z_{j-1/2},\,z_{j+1/2}]$ of size $\dx\,\dz$.
Inside every control volume we average with respect to $x$ and $z$ leading to
the semi-discrete scheme
\begin{equation}\label{sclh3}
\frac{d Q_{i,j}(t)}{dt}=-\frac{1}{\dx}(F_{i+1/2,j}-F_{i-1/2,j})-
   \frac{1}{\dz}(H_{i,j+1/2}-H_{i,j-1/2})\equiv L_{i,j}(Q),
\end{equation}
where,
\begin{equation}
Q_{i,j}=\frac{1}{\dx}\frac{1}{\dz}\int_{x_{i-1/2}}^{x_{i+1/2}}\,
        \int_{z_{j-1/2}}^{z_{j+1/2}}Q(x,z,t)\,dz\,dx,
\end{equation}
\begin{align}
F_{i+1/2,j}&=\frac{1}{\dz}\int_{z_{j-1/2}}^{z_{j+1/2}}\F(Q(x_{i+1/2},z,t))\,dz\,,\label{fluxes31}\\
H_{i,j+1/2}&=\frac{1}{\dx}\int_{x_{i-1/2}}^{x_{i+1/2}}\H(Q(x,z_{j+1/2},t))\,dx.\label{fluxes32}
\end{align}
We approximate the expressions in eqns. (\ref{fluxes31})-(\ref{fluxes32}) by suitable Gaussian quadrature formulas,
\begin{align}
F_{i+1/2,j}\approx\frac{1}{2}\sum_{z_{Gp}}w_{z_{Gp}}\F(Q(x_{i+1/2},z_{Gp},t))\,dz\,,\label{fluxesh31}\\
H_{i,j+1/2}\approx\frac{1}{2}\sum_{x_{Gp}}w_{x_{Gp}}\H(Q(x_{Gp},z_{j+1/2},t))\,dz,\label{fluxesh32}
\end{align}
where $x_{Gp}$ and $z_{Gp}$ are prescribed Gauss points with corresponding weights $w_{x_{Gp}}$ and $w_{z_{Gp}}$ respectively.
The computation of eqns. (\ref{fluxesh31})-(\ref{fluxesh32}) is performed via a high-resolution approach that makes use of a WENO reconstruction procedure; after this step is completed, a polynomial of prescribed order is obtained at every cell, and therefore, at every cell interface, accurate flux calculations can be performed by taking extrapolated boundary values.
We briefly describe the WENO reconstruction procedure that is used in this article; we opted for the technique described in \cite{balsara1,balsara2} in its third order (quadratic reconstruction) version. This technique makes extensive use of the structure of the reconstruction procedure in one dimension, adding some additional mixed terms (``cross terms'') that are efficiently computed by reduced stencils. It is an optimal and easy way to implement the algorithm for achieving high-order reconstructions in 2 and 3 dimensions; it also defines an unique polynomial in every cell, which is particularly useful when space dependent source terms such as viscosity are considered.
At a given time $t$ (the subscript indicating time is omitted throughout this derivation), given the set of averaged values $\{Q_{i,j}\}$ for the whole domain, at every cell, the reconstruction procedure seeks a quadratic expansion upon a linear combination of Legendre polynomials rescaled in local coordinates $(x,z)=[-1/2,\,1/2]\times [-1/2,\,1/2]$ expressed in the form
\begin{equation}\label{rec}
Q(x,z)=Q_0+Q_xP_1(x)+Q_{xx}P_2(x)+Q_zP_1(z)+Q_{zz}P_2(z)+Q_{xz}P_1(x)P_1(z),
\end{equation}
\begin{equation}
P_1(x)=x \qquad P_2(x)=x^2-\frac{1}{12}.
\end{equation}
Except for the last term in (\ref{rec}), every coefficient can be computed by performing a dimension-by-dimension reconstruction, which we now illustrate. We assign the subscript ''0'' to the cell where we are computing the coefficients, other values indicating location and direction with respect to $Q_0$ (note that the notation is coherent with the fact that the first coefficient in the expansion $Q_0$, holds $Q_0=Q_{ij}$, i.e., the centered value). Next, for this particular problem we define three stencils
\begin{equation}
S^1=\{Q_{-2},Q_{-1},Q_{0}\}\,,\quad S^2=\{Q_{-1},Q_{0},Q_{1}\}\,,\quad S^3=\{Q_{0},Q_{1},Q_{2}\}\,,
\end{equation}
and in every stencil we compute a polynomial of the form
\begin{equation}
Q^{(i)}(x)=Q_0^{(i)}+Q_x^{(i)}P_1(x)+Q_{xx}^{(i)}P_2(x)\qquad i=1,2,3.
\end{equation}
The coefficients are given by
\begin{eqnarray}
S^{1}&:&Q^{(1)}_x=-2Q_{-1}+Q_{-2}/2+3Q_0/2,\quad Q^{(1)}_{xx}=(Q_{-2}-2Q_{-1}+Q_0)/2\,,\\
S^{2}&:&Q^{(2)}_x=(Q_1-Q_{-1})/2,\quad Q^{(2)}_{xx}=(Q_{-1}-2Q_{0}+Q_1)/2\,,\\
S^{3}&:&Q^{(3)}_x=-3Q_{0}/2+2Q_{1}-Q_{2}/2,\quad Q^{(3)}_{xx}=(Q_{0}-2Q_{-1}+Q_2)/2.
\end{eqnarray}
For every polynomial we calculate a smoothness indicator defined as
\begin{equation}
IS^{(i)}=\left(Q^{(i)}_x\right)^2+\frac{13}{3}\left(Q^{(i)}_{xx}\right)^2\,,
\end{equation}
leading to the following WENO weights:
\begin{equation}
\omega^{(i)}=\frac{\alpha^{(i)}}{\sum_{i=1}^3\alpha^{(i)}}\,,\quad \alpha^{(i)}=\frac{\lambda^{(i)}}{(\epsilon+IS^{(i)})^r}\,^,
\end{equation}
where $\epsilon$ is a parameter introduced in order to avoid division by zero; usually $\epsilon=10^{-12}$. The scheme is rather insensitive to the parameter $r$, which we set $r=5$. The parameter $\lambda$ is usually computed in an optimal way to increase the accuracy of the reconstruction at certain points; we opt for a centered approach instead, thus $\lambda^{(1)}=\lambda^{(3)}=1$, while $\lambda^{(2)}=100$. The 1D reconstructed polynomial is given by
\begin{equation}
Q(x)=\omega^{(1)}Q^{(1)}(x)+\omega^{(2)}Q^{(2)}(x)+\omega^{(3)}Q^{(3)}(x).
\end{equation}
Next, a 1D reconstruction in the $z$ direction is performed in a totally analogous way. Finally, we address the computation of the mixed term $Q_{xz}$, which is calculated in a 2D fashion. Keeping the same convention regarding location subscripts as in 1D, \cite{balsara1} considers 4 formulas for the cross term upon taking all the moments around the cell. The expressions for the cross term are:
\begin{eqnarray}
Q_{xz}^{(1)}&=&Q_{1,1}-Q_{0,0}-Q_x-Q_z-Q_{xx}-Q_{zz},\\
Q_{xz}^{(2)}&=&-Q_{1,-1}+Q_{0,0}+Q_x-Q_z+Q_{xx}+Q_{zz},\\
Q_{xz}^{(3)}&=&-Q_{-1,1}+Q_{0,0}-Q_x+Q_z+Q_{xx}+Q_{zz},\\
Q_{xz}^{(4)}&=&Q_{-1,-1}-Q_{0,0}+Q_x+Q_z-Q_{xx}-Q_{zz},
\end{eqnarray}
and the corresponding smoothness indicators are given by
\begin{equation}
IS^{(i)}=4\left(Q^{(i)}_{xx}\right)^2+4\left(Q^{(i)}_{zz}\right)^2+\left(Q^{(i)}_{xz}\right)^2.
\end{equation}
Note that in the first part of the reconstruction, when the weights were computed, a larger suboptimal weight was assigned to the central stencil, which is a way to ensure stability and robustness of the algorithm by sacrificing additional order in the approximation (for more details, see \cite{dumbserkaser}). However, for this term, the numerators assigned to the corresponding $\alpha$'s remains the same for every expression. The computation of this term concludes the reconstruction procedure, and now we have at our disposal one polynomial per cell that can be used to calculate values at the boundaries or inside the cell. The next step in our numerical scheme consists of the calculation of the numerical fluxes (\ref{fluxesh31})-(\ref{fluxesh32}), which will use extrapolated boundary values of the reconstructed polynomials. Rather than the use of the classical WENO scheme (as in \cite{liu} or \cite{shu1}), which performs this calculation via a first order flux, we opt for the WENO-TVD approach described in \cite{titarevweno}. We make use of the 2D extension of the flux-limiter-centered scheme (FLIC) approach presented in \cite{torobillett,torolibro}, which is a second-order, centered and non-oscillatory flux. In our case, it consists of a flux-limited version of a generalized Lax Wendroff flux, using  as a low-order flux the GFORCE (generalized first order centered) flux \cite{toroforce}, which can be interpreted as a convex combination of Lax-Friedrichs and Lax-Wendroff-type of fluxes:
\begin{equation}
F_{i+1/2,j}^{FLIC}=F_{i+1/2,j}^{GFORCE}+\psi_{i+1/2,j}\left(F_{i+1/2,j}^{LW}-F_{i+1/2,j}^{GFORCE}\right),
\end{equation}
where
\begin{align}
F_{i+1/2,j}^{GFORCE}&=F_{i+1/2,j}^{GFORCE}\left(Q_{i+1/2,j}^L,Q_{i+1/2,j}^R\right)=\omega F_{i+1/2,j}^{LW}+(1-\omega)F_{i+1/2,j}^{LF}\,,\\
F_{i+1/2,j}^{LF}&=\frac12\left(\F\left(Q_{i+1/2,j}^L\right)+\F\left(Q_{i+1/2,j}^R\right)-\frac12\frac{\dx}{\dt}\left(Q_{i+1/2,j}^R-Q_{i+1/2,j}^L\right)\right)\,,\\
F_{i+1/2,j}^{LW}&=\F\left( Q_{i+1/2,j}^* \right)\,,\\
Q_{i+1/2,j}^*&=\frac12\left(Q_{i+1/2,j}^L+Q_{i+1/2,j}^R\right)-\frac{\dt}{\dx}\left(\F\left(Q_{i+1/2,j}^R\right)-\F\left(Q_{i+1/2,j}^L\right)\right).
\end{align}
The parameter $\omega$ varies between 0 and 1, and is chosen in a compatible manner with the CFL number in order to ensure monotonicity; throughout the numerical experiments presented in this article, we will restrict ourselves to the usual FORCE flux, i.e., $\omega=0.5$.
We have omitted the formulas for the remaining cell boundaries, but they can be derived in a straightforward manner. Also note that even though the formulas are written along the boundary '$i+1/2,j$', the use of the Gaussian quadrature formula will replace the axes'$j$' by Gauss points and therefore this subscript must be understood in that sense.
It is important to notice that so far we are deriving expressions for the semi-discrete approximation of the system of conservation laws, however, the fluxes include the parameter $\dt$ which arises from the averaging operators that originate these fluxes. Thus, in the spatial discretization of the system, the time stepping enters just as a parameter. At the end of the derivation of the scheme, when we present the time discretization of eq. (\ref{sclh3}), $\dt$ will be considered as ``marching parameter'' in the sense that its inclusion in the formulas will generate an updated state in time.
The function $\psi_{i+1/2,j}=\psi_{i+1/2,j}(r^L_{i+1/2,j},r^R_{i+1/2,j})$ is a flux limiter; a slight variation of the usual limiters has to be considered in this context since we use a centered flux instead of an upwind approach (the reader can refer to \cite[Ch 13.]{torolibro} for more details); in our case we mainly use the VanLeer limiter, which on its centered version reads:
\begin{equation}
\psi(r)=\begin{cases}
0 &\mbox{if } r\leq 0,\\
\frac{2r}{1+r} &\mbox{if } 0\leq r \leq 1,\\
\phi_g+\frac{2r(1-\phi_g)}{1+r} &\mbox{if } r\geq 1,
\end{cases}
\qquad \phi_g=\frac{1-|c|}{1+|c|},
\end{equation}
where $c$ corresponds to the Courant number which depends on the problem.
To implement the limiter in the ''conservation paradigm'', the flow parameter $r$ is defined
via the total energy of the system,
\begin{equation}
e=c_v\theta\pi+\frac12(u^2+v^2+w^2)+gz.
\end{equation}
The left and right flow parameters are then defined as,
\begin{equation}
r^L_{i+1/2,j}=\frac{e^R_{i-1/2,j}-e^L_{i-1/2,j}}{e^R_{i+1/2,j}-e^L_{i+1/2,j}}\,,
     \qquad r^R_{i+1/2,j}=\frac{e^R_{i+3/2,j}-e^L_{i+3/2,j}}{e^R_{i+1/2,j}-e^L_{i+1/2,j}},
\end{equation}
and finally,
\begin{equation}
\psi_{i+1/2,j}=\min(\psi(r^L_{i+1/2,j}),\psi(r^R_{i+1/2,j})).
\end{equation}
The above described procedure starts with a set of averaged values and ends with a numerical approximation of the space operators involved in eq. (\ref{scl3}). The resulting scheme is still continuous in time, and we conclude this section by discretizing this operator in a manner that is consistent with the choices that we have made in the generation of the space discretization operator. At a given starting time $t^n$, we begin by considering the semi-discrete scheme
\begin{equation}
\frac{d Q_{i,j}(t)}{dt}= L_{i,j}(Q),
\end{equation}
bringing the system to a final state $t^{n+1}$ with a time stepping $\dt$. In order to preserve high-order and non-oscillatory properties in time, we consider the well-known family of explicit TVD Runge-Kutta schemes \cite{shutvd}, in particular its third order version
\begin{eqnarray}
Q_{i,j}^{n+\frac13}&=&Q_{i,j}^n+\dt \,L_{i,j}(Q_{i,j}^n),\\
Q_{i,j}^{n+\frac23}&=&\frac34Q^n_{i,j}+\frac14Q_{i,j}^{n+\frac13}+\frac14\dt\, L_{i,j}(Q_{i,j}^{n+\frac13}), \\
Q_{i,j}^{n+1}&=&\frac13 Q^n_{i,j}+\frac23 Q_{i,j}^{n+\frac23} +\frac23\dt\, L_{i,j}(Q_{i,j}^{n+\frac23}). \\
\end{eqnarray}
We conclude this section with the inclusion of the source term. The source term appearing in eq. (\ref{st3}) will not depend on space nor space derivatives, and therefore it can be averaged in space and solved in the same manner as the above presented time discretization, by replacing $L_{i,j}(Q_{i,j}^n)$ by $\S(Q_{i,j}^n)$. If we denote by the $\mathfrak{L}(\dt)$ the fully discrete operator that brings the system of conservation laws (\ref{scl3}) $\dt$ units ahead in time, and by $\mathfrak{S}(\dt)$ the fully discrete operator that updates the source term (\ref{st3}) in $\dt$ units, we preserve, at least, second order accuracy in time by implementing a Strang splitting \cite{strang} in the form
\begin{equation}
Q^{n+1}_{i,j}=\mathfrak{S}(\dt/2)\mathfrak{L}(\dt)\mathfrak{S}(\dt/2)Q_{i,j}^n.
\end{equation}

In this way we have derived a fully discrete, high-order and total variation diminishing
numerical scheme  for treating a system of balance laws arising from dimensional reduction
of the original set of 3D primitive equations. In the next section,
numerical experiments are performed with the aim of understanding the consequences of the dimensional
reduction and the coupling effects, and its ability to represent atmospheric phenomena
in a plausible way.
The test are motivated by standard study cases for non-hydrostatic model
development (see \cite{skamarock,robert,knothjcp}),
and therefore a qualitative comparison with published model outputs is possible.

\section{Numerical experiments}\label{tests}

The purpose of the numerical experiments is to use well-known test cases that are
used extensively to test dynamical cores of atmospheric models so that our
results can easily be compared with results in the literature.
There are quite a few such test cases published, and we have selected three.
In addition we have constructed a new test case that is used to show the capabilities
of the dimensionally-reduced model.

\subsection{Convective bubble in a neutral atmosphere}\label{section:test1}

\begin{figure}[H]
\caption{Potential temperature for the convective bubble test case. Colormap of the first layer at $t=0$, 120, 300 and 600[s]. $\dx=\dz=125[m]$, $160\times 80$ elements.}
\begin{minipage}[b]{0.5\linewidth}
\centering
\includegraphics[scale=0.45]{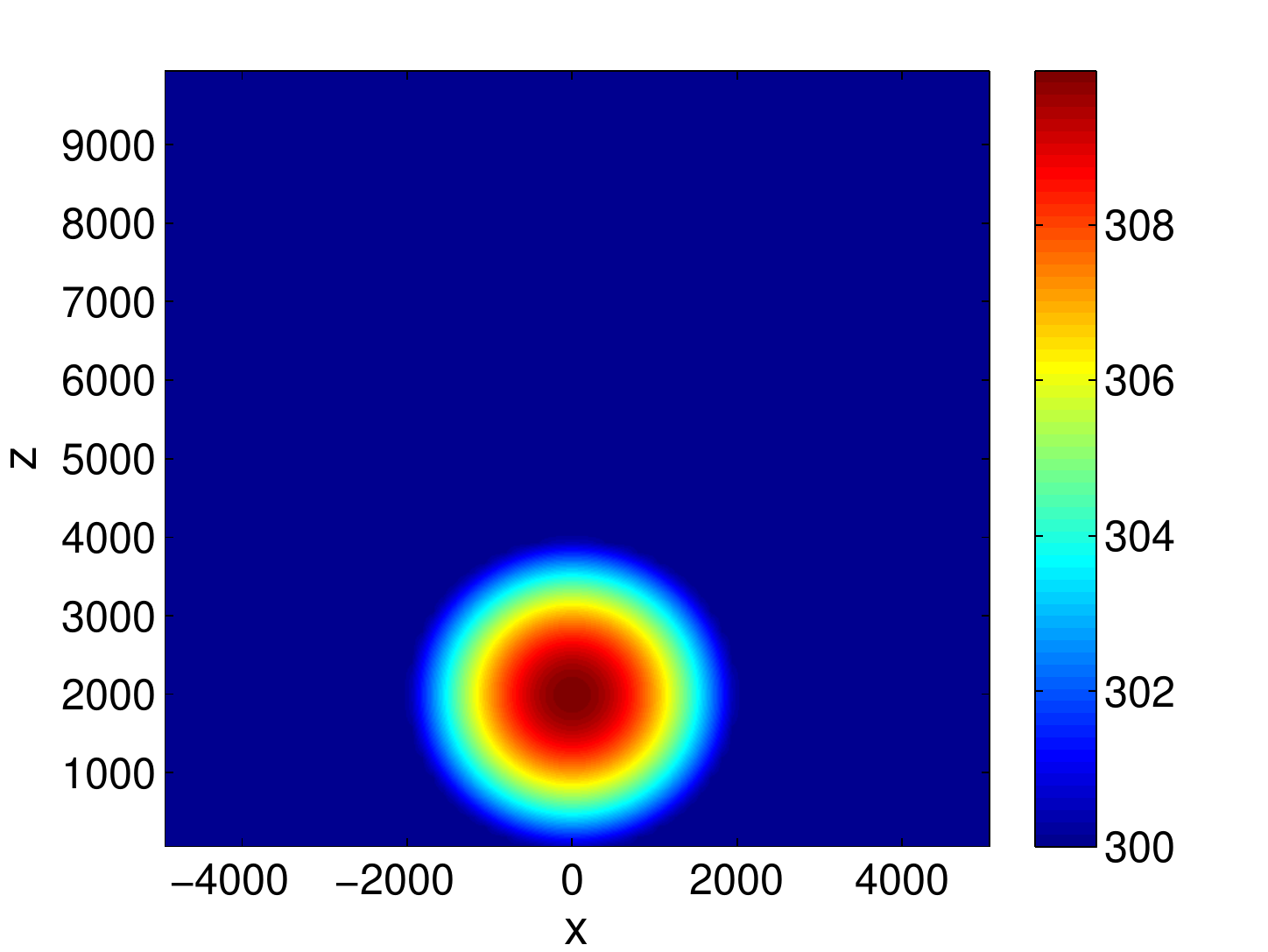}
\includegraphics[scale=0.45]{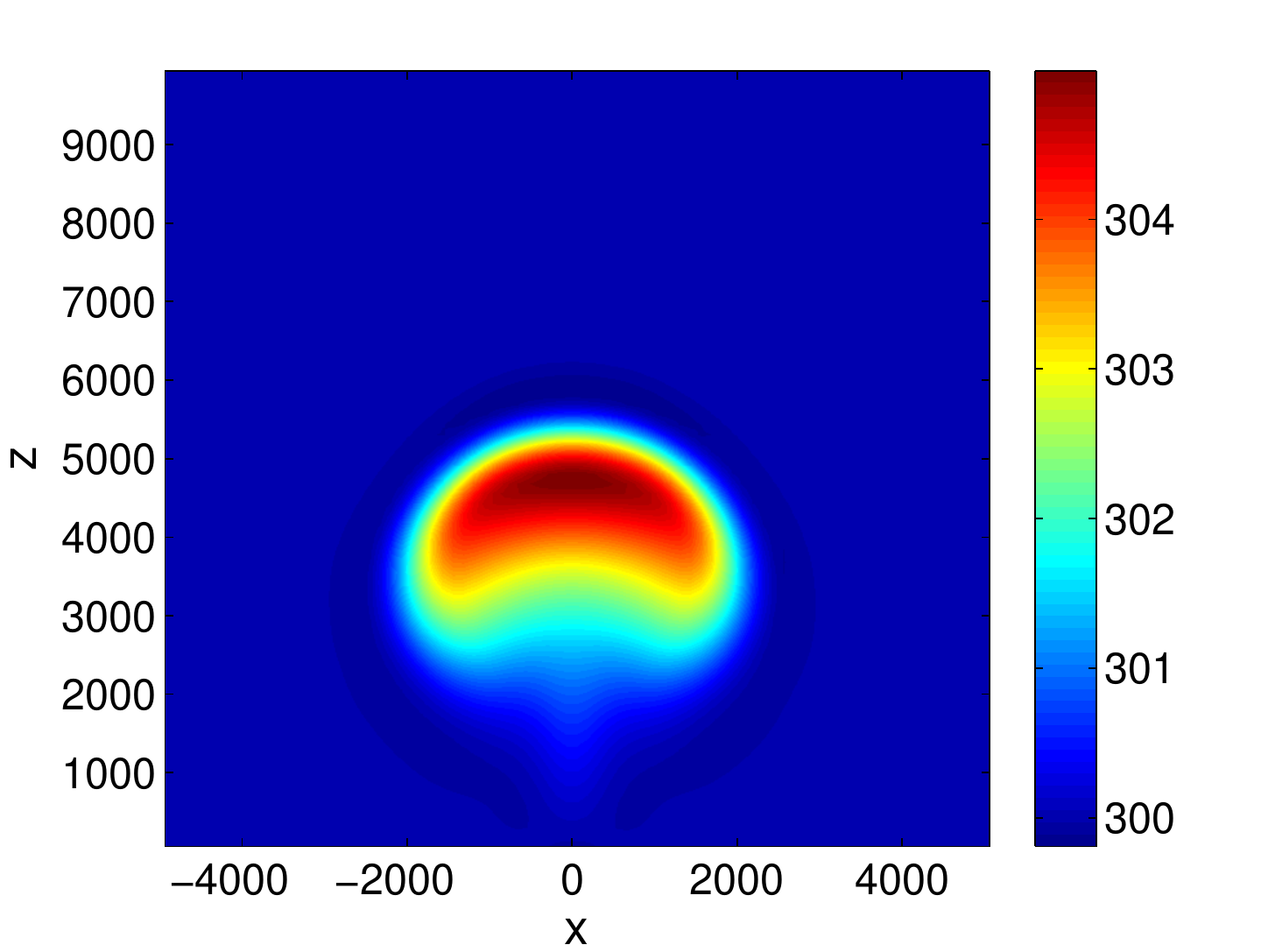}
\end{minipage}
\hspace{0.0cm}
\begin{minipage}[b]{0.5\linewidth}
\centering
\includegraphics[scale=0.45]{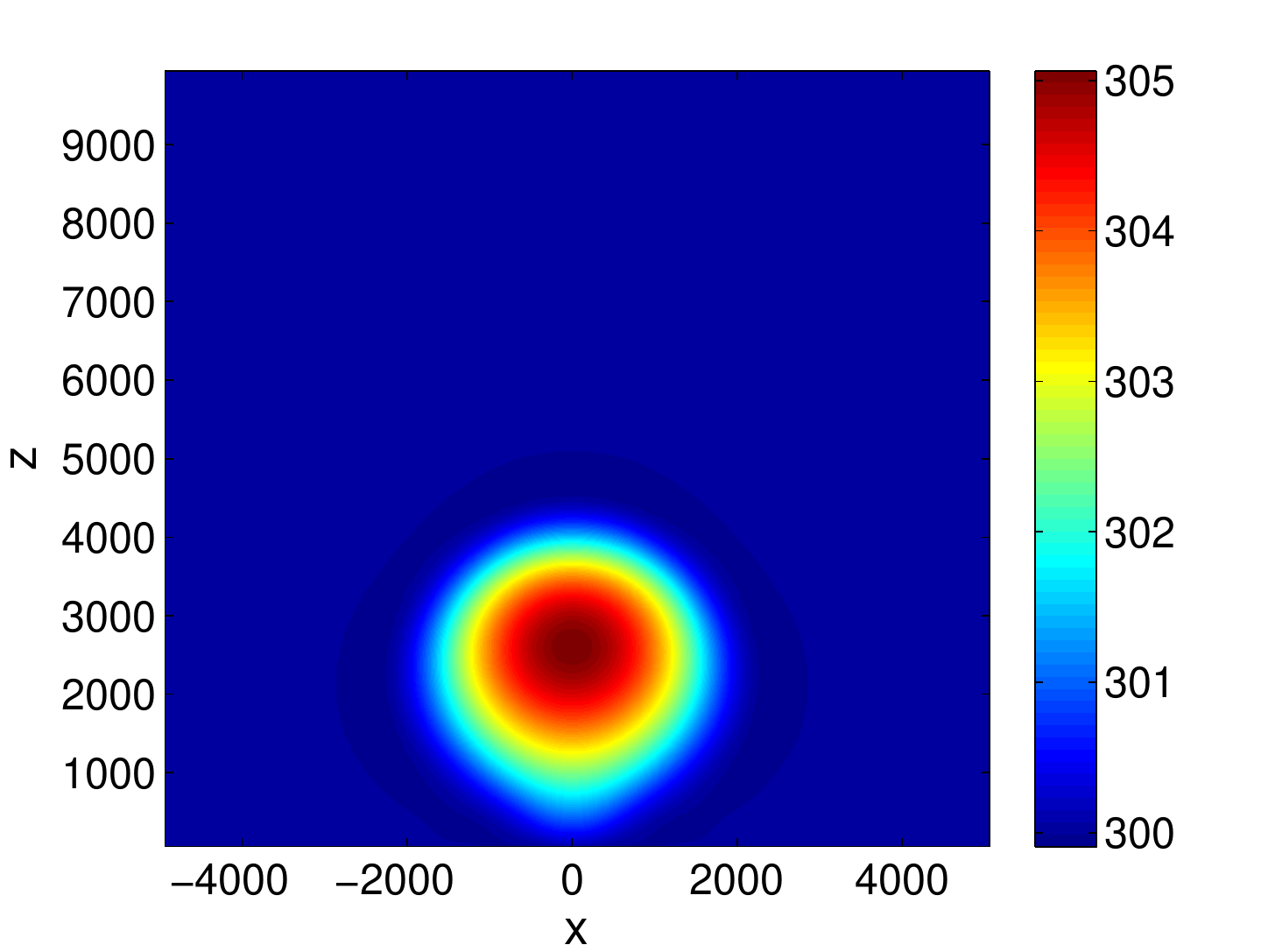}
\includegraphics[scale=0.45]{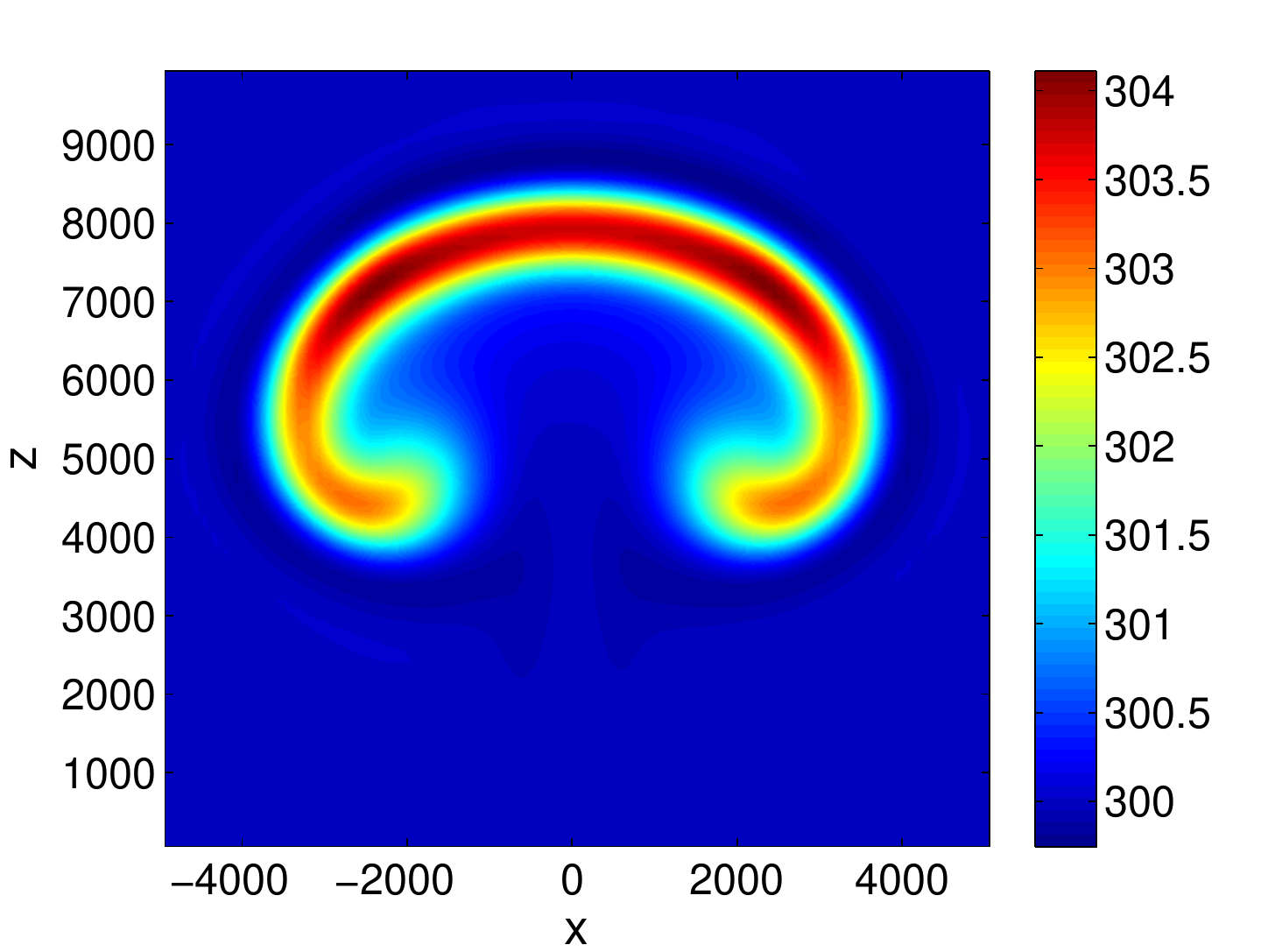}
\end{minipage}
\label{fig:figure1}
\end{figure}

\begin{figure}[H]
\caption{Potential temperature for the convective bubble test case. Colormap of the residual ($\theta_1-\theta_2$) at $t=120$ and 600[s]. $\dx=\dz=125[m]$, $160\times 80$ elements.}
\begin{minipage}[b]{0.5\linewidth}
\centering
\includegraphics[scale=0.45]{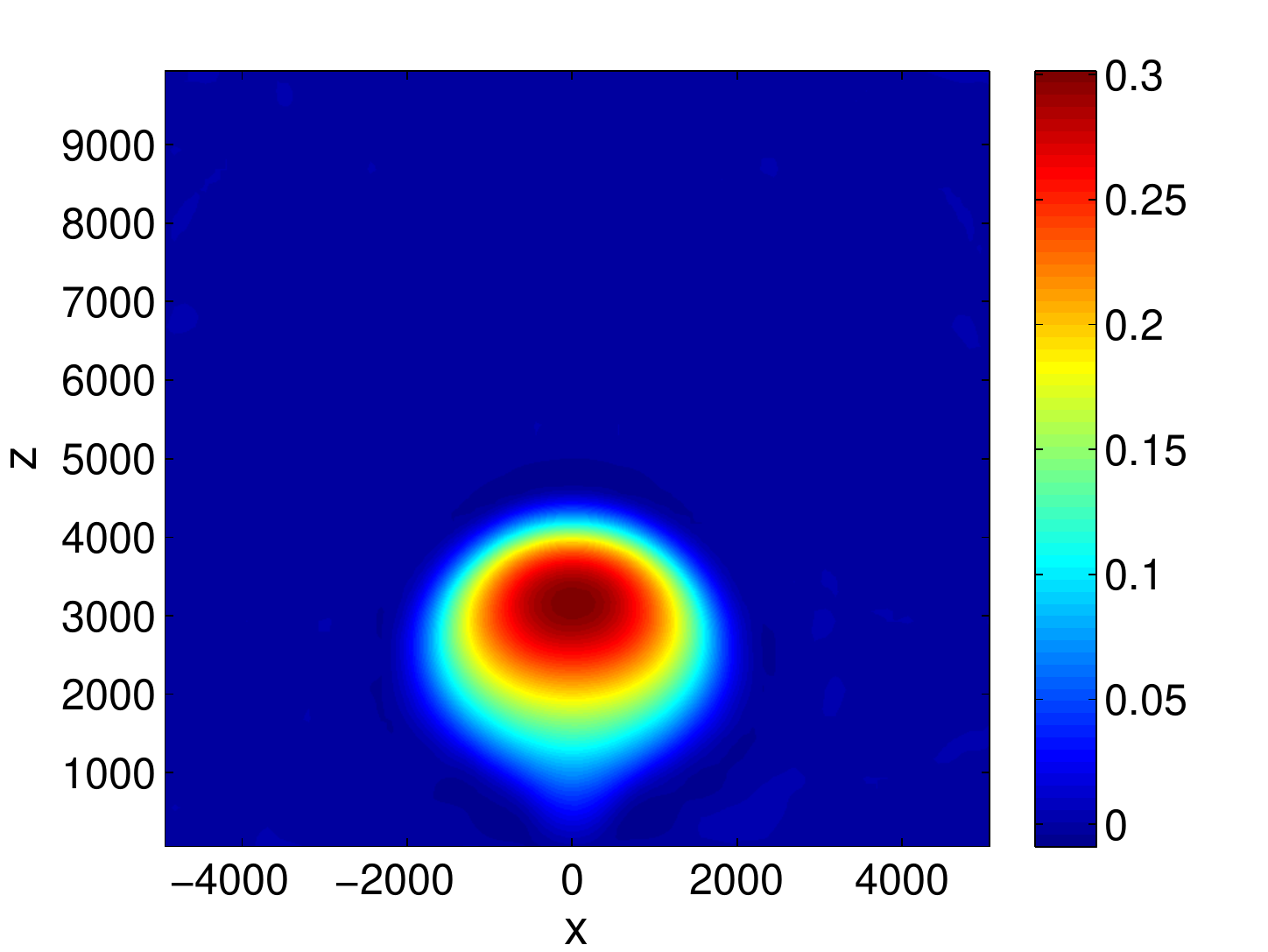}
\end{minipage}
\begin{minipage}[b]{0.5\linewidth}
\centering
\includegraphics[scale=0.45]{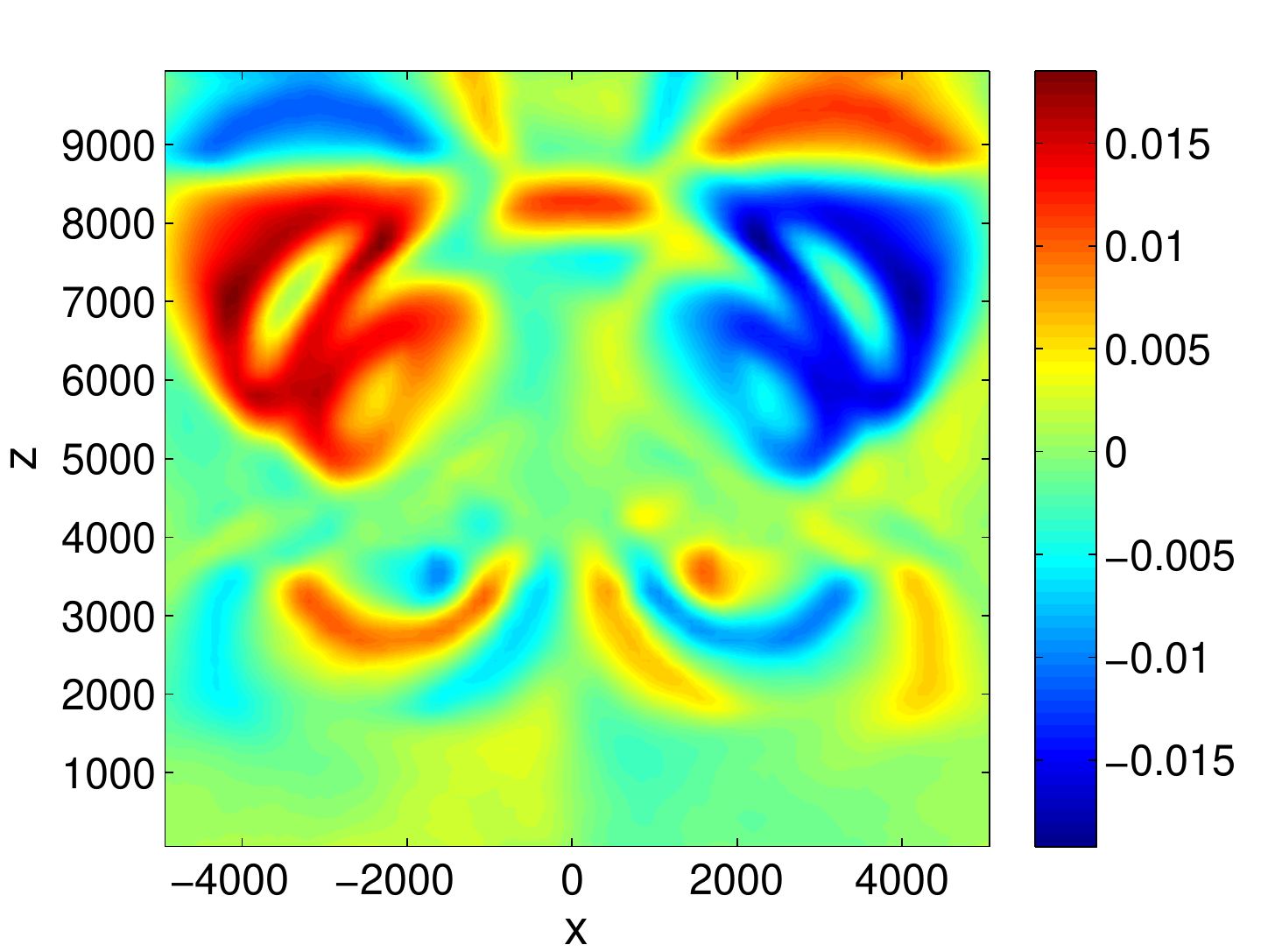}
\end{minipage}
\label{fig:figure2}
\end{figure}

\begin{figure}[H]
\caption{Velocity field colormaps at $t=120$, 300 and 600 [s] for the first layer. Left: horizontal velocity. Right: vertical velocity. $\dx=\dz=125[m]$, $160\times 80$ elements.}
\begin{minipage}[b]{0.5\linewidth}
\centering
\includegraphics[scale=0.45]{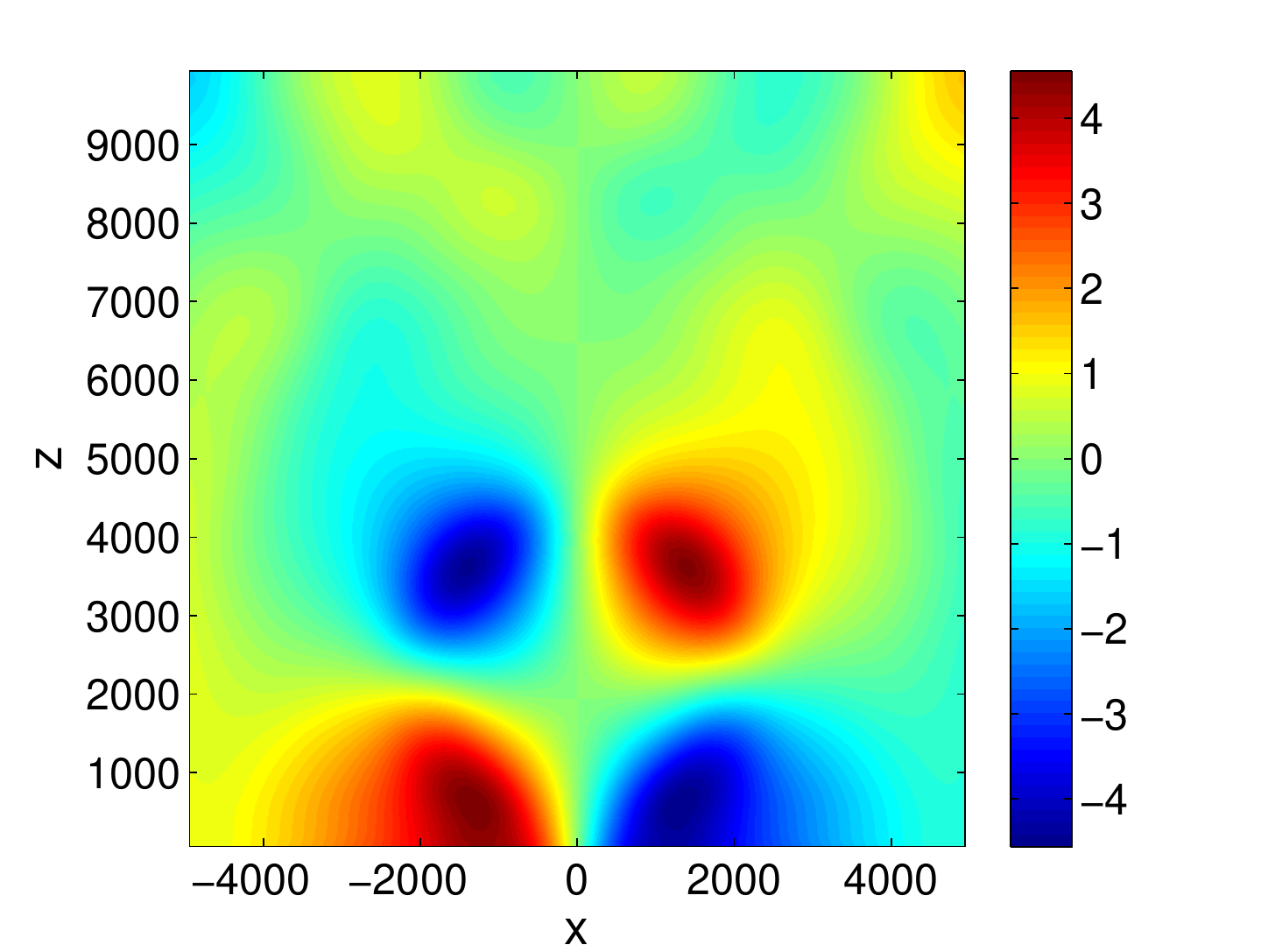}
\includegraphics[scale=0.45]{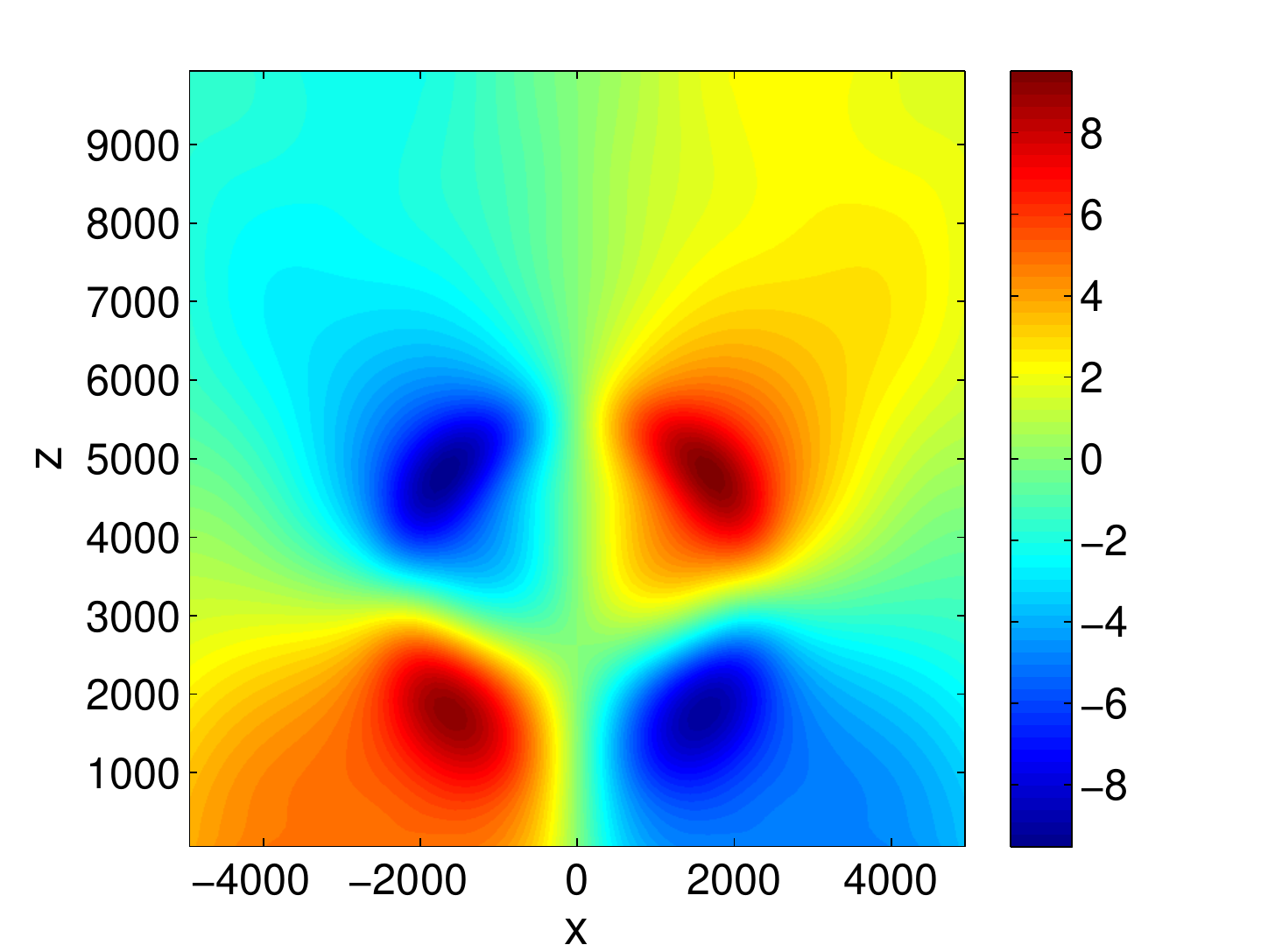}
\includegraphics[scale=0.45]{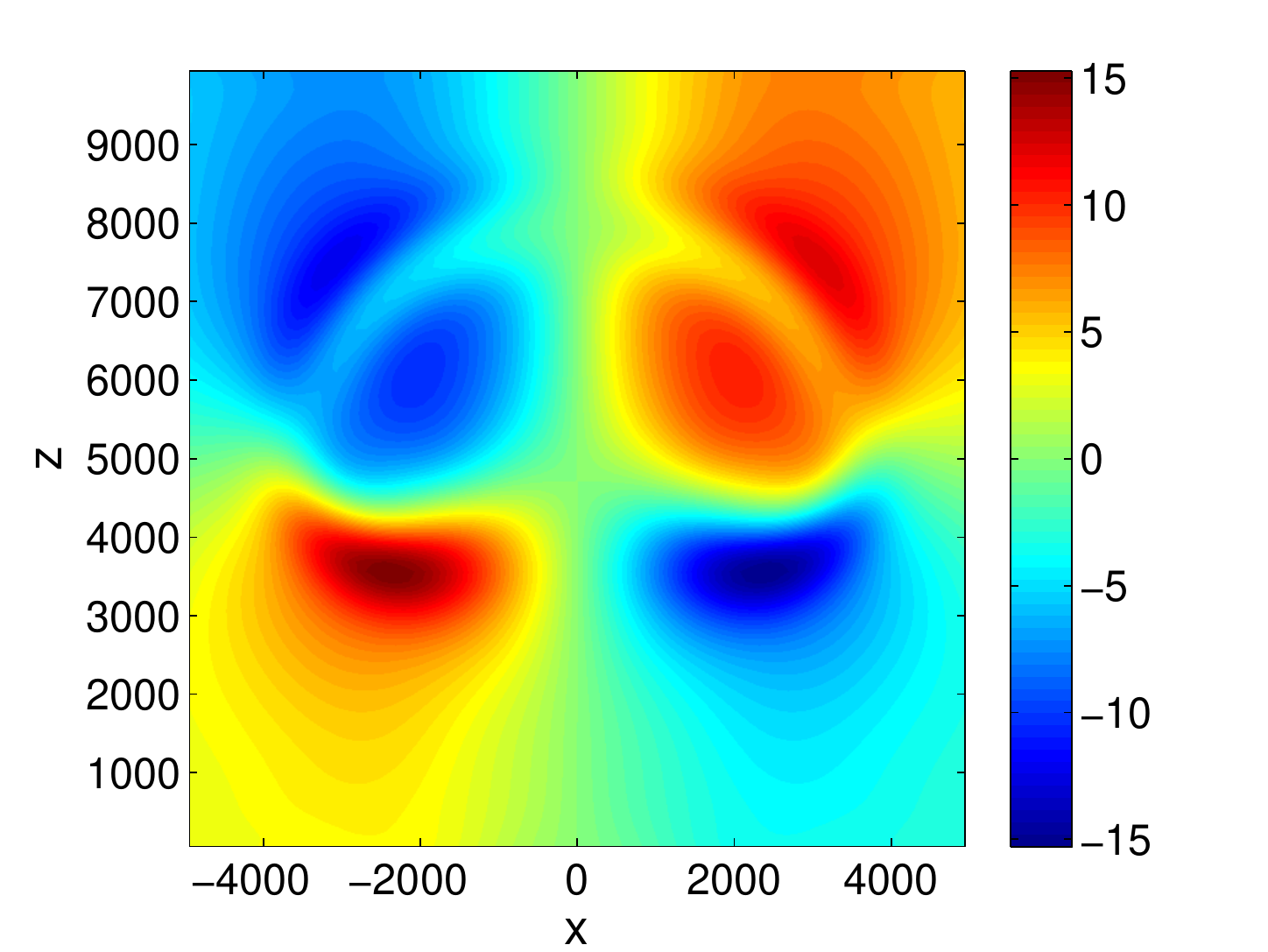}
\end{minipage}
\hspace{0.0cm}
\begin{minipage}[b]{0.5\linewidth}
\centering
\includegraphics[scale=0.45]{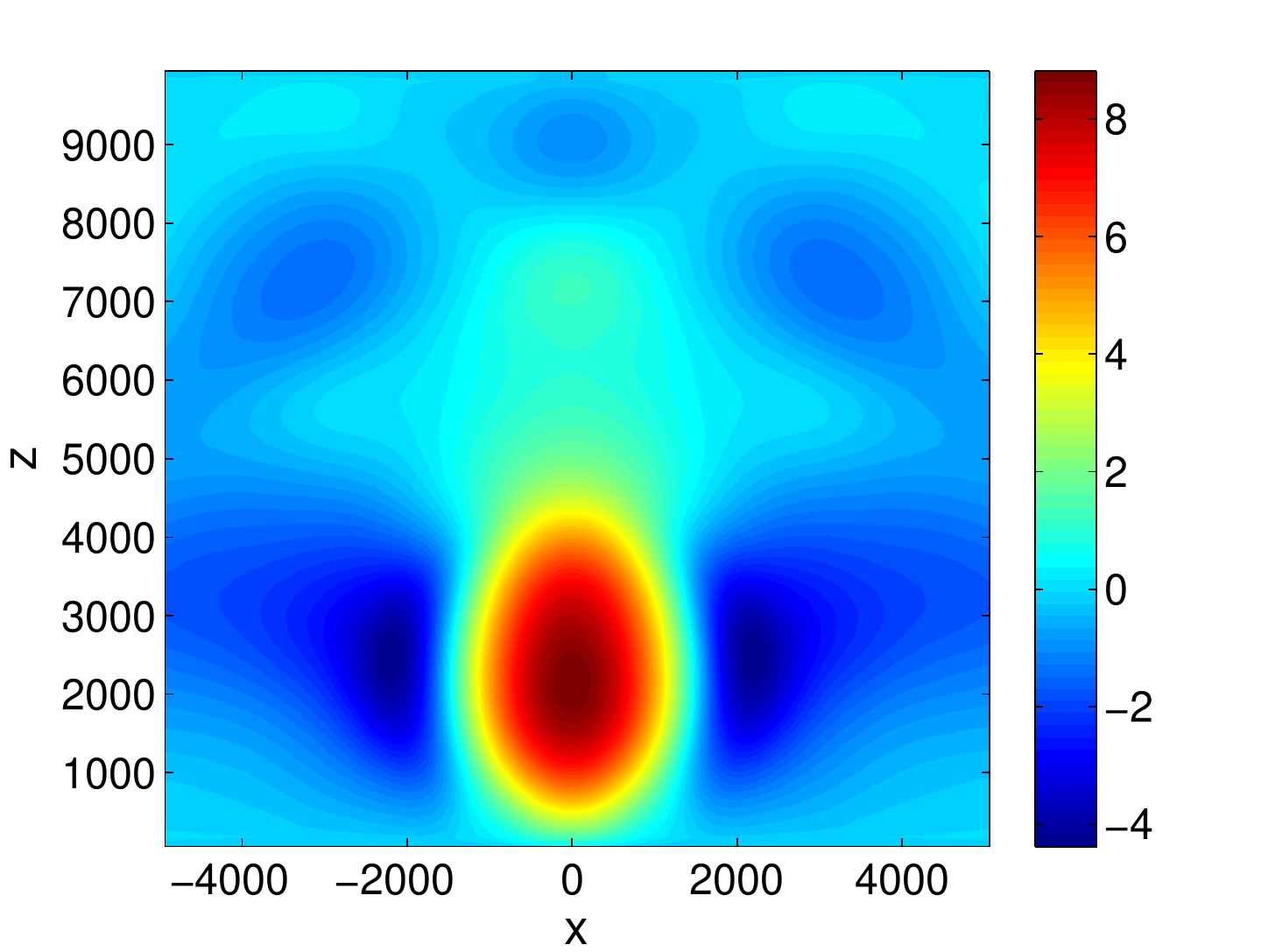}
\includegraphics[scale=0.45]{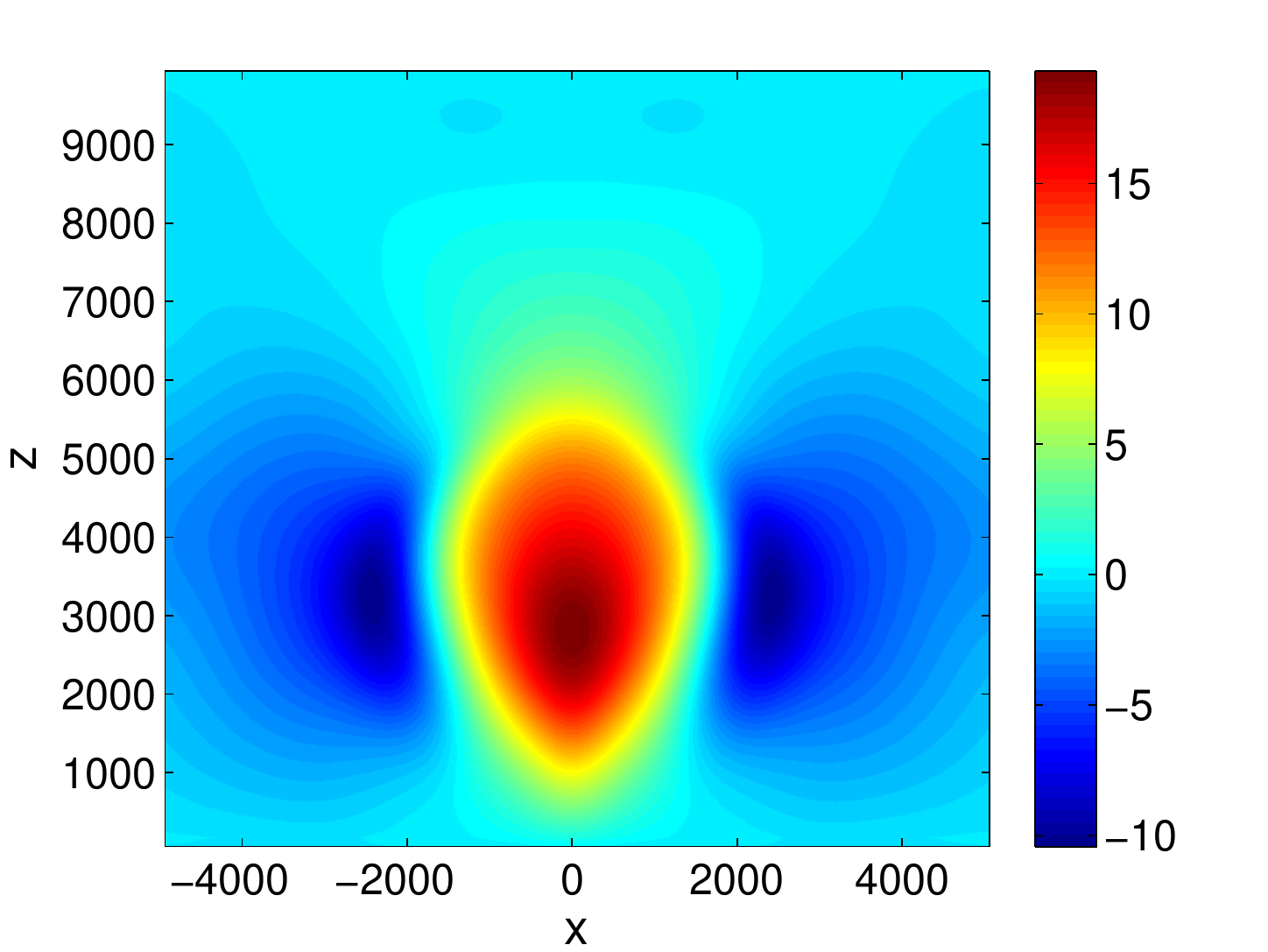}
\includegraphics[scale=0.45]{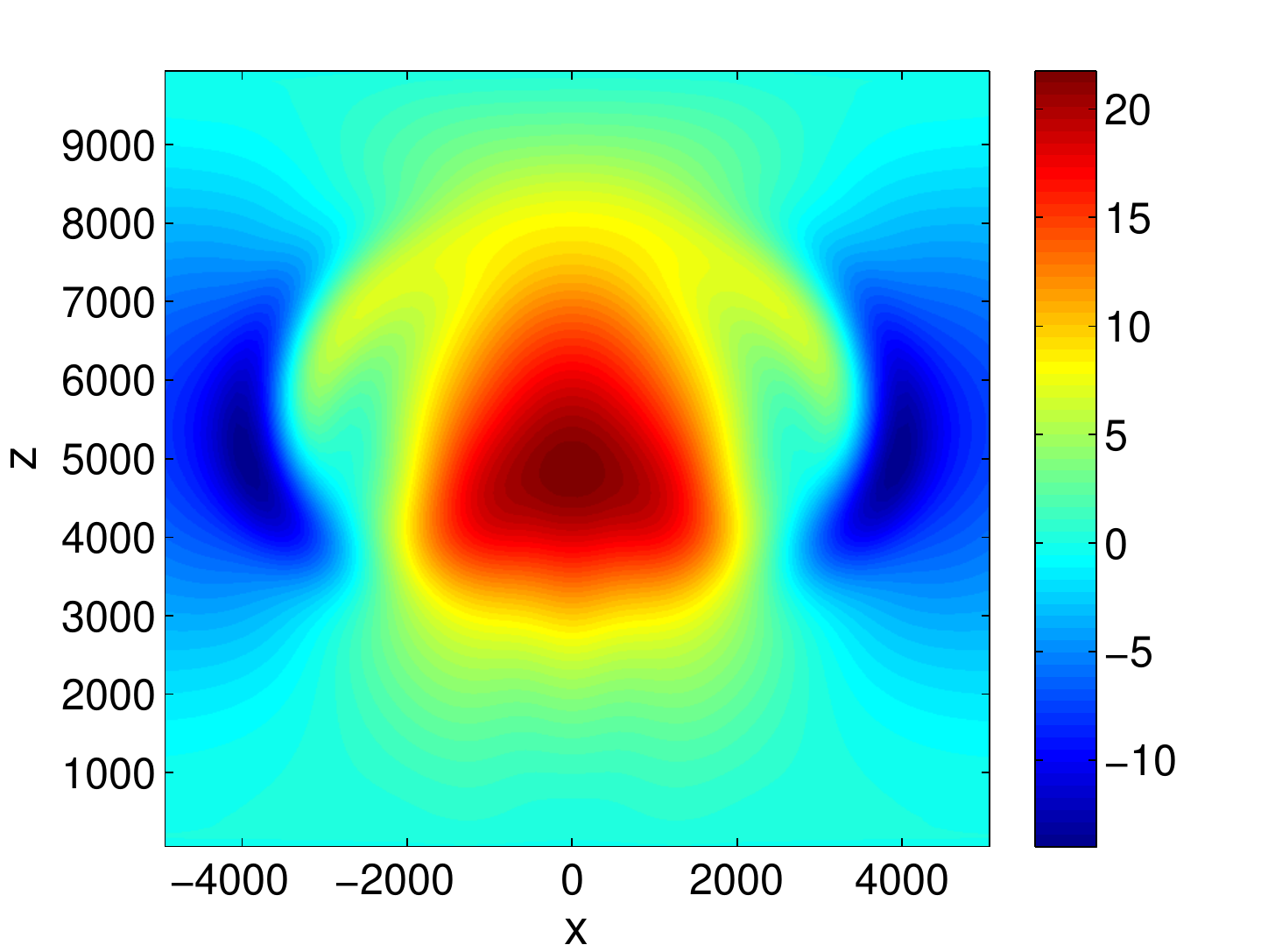}
\end{minipage}
\label{fig:figure3}
\end{figure}

\begin{figure}[H]
\caption{Velocity field Vector plots at $t=0$, 120, 300 and 600 [s] for the first layer. $\dx=\dz=125[m]$, $160\times 80$ elements.}
\begin{minipage}[b]{0.5\linewidth}
\centering
\includegraphics[width=\textwidth,height=0.25\textheight]{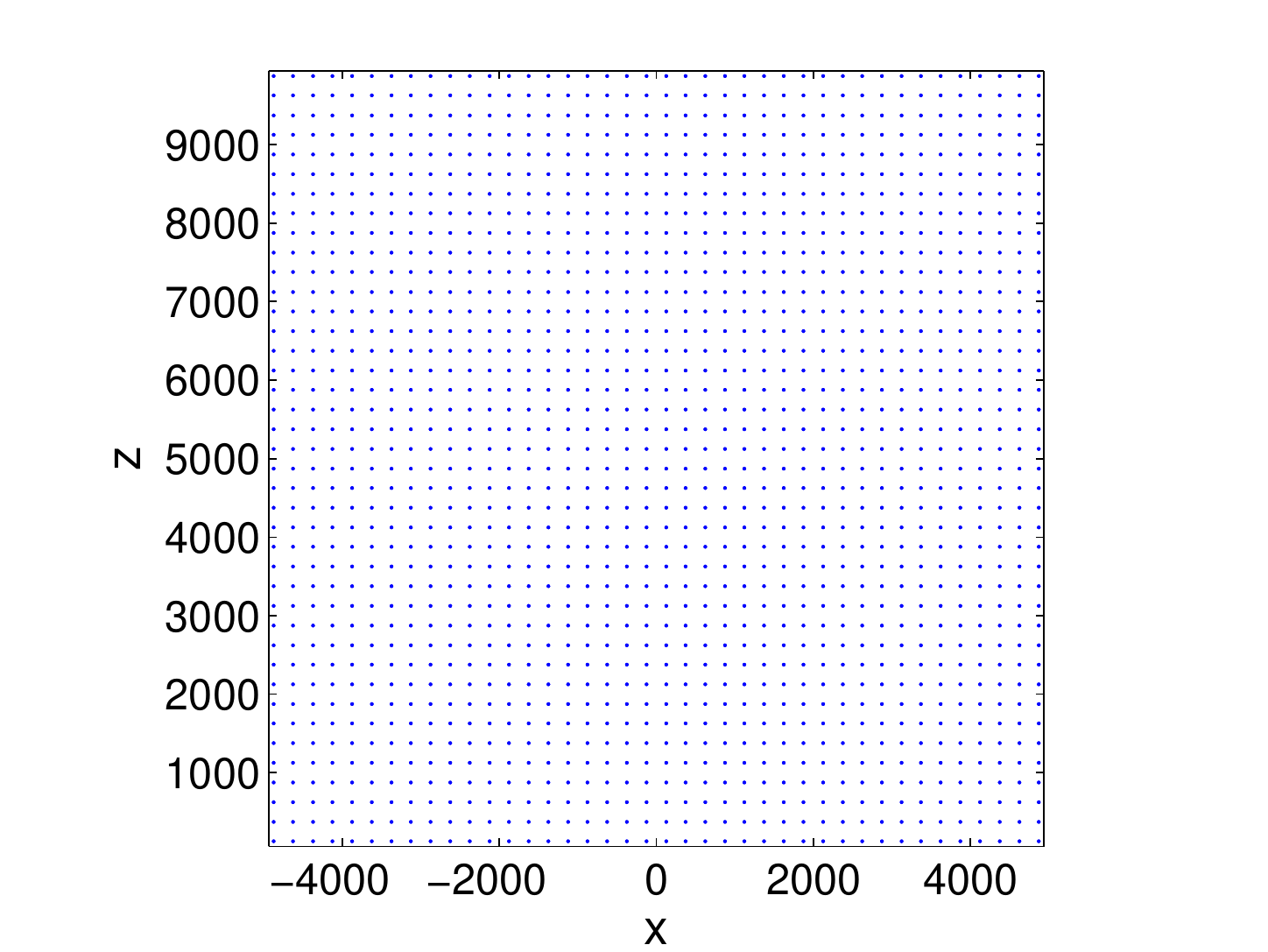}
\includegraphics[width=\textwidth,height=0.25\textheight]{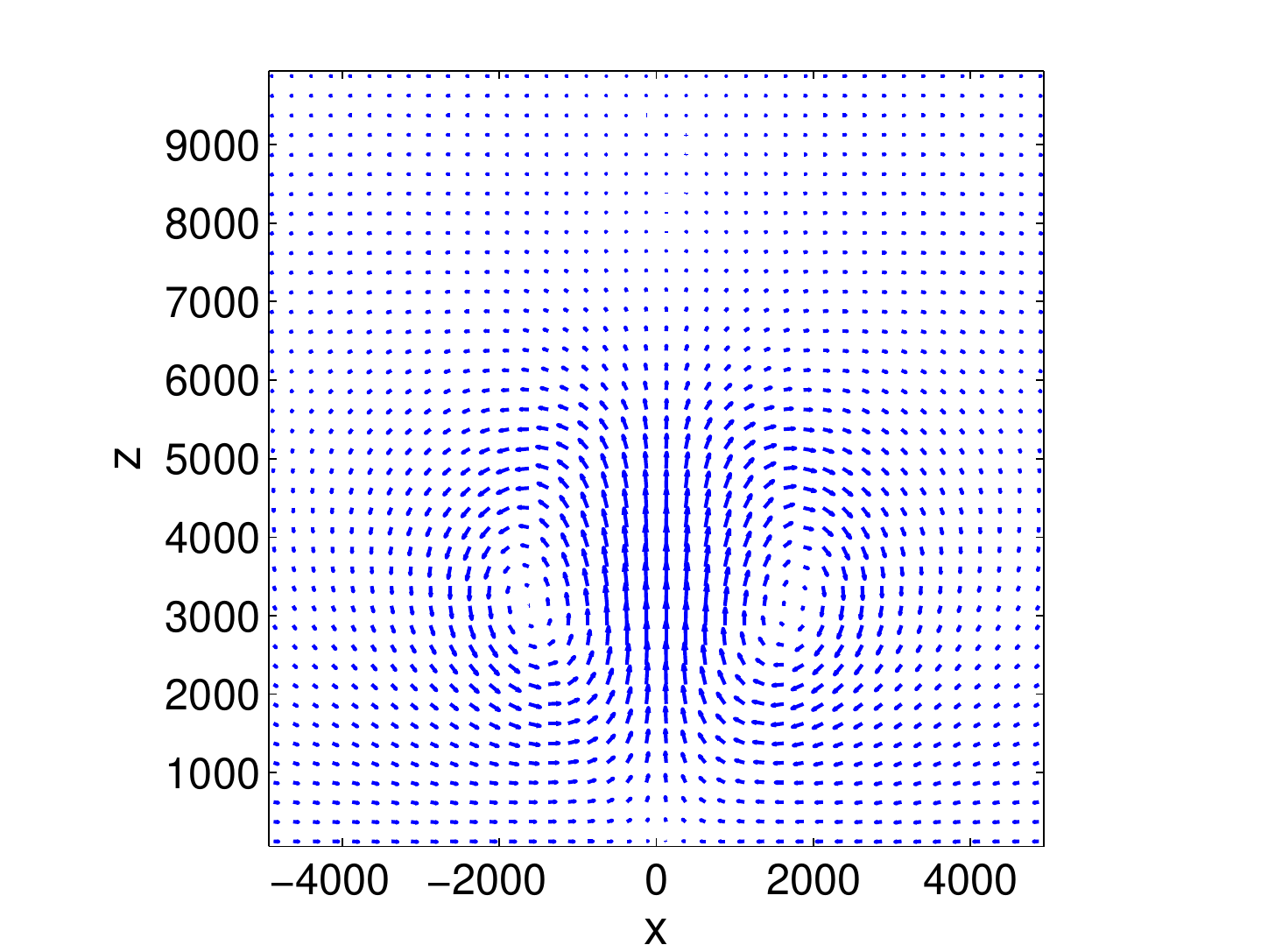}
\end{minipage}
\hspace{0.0cm}
\begin{minipage}[b]{0.5\linewidth}
\centering
\includegraphics[width=\textwidth,height=0.25\textheight]{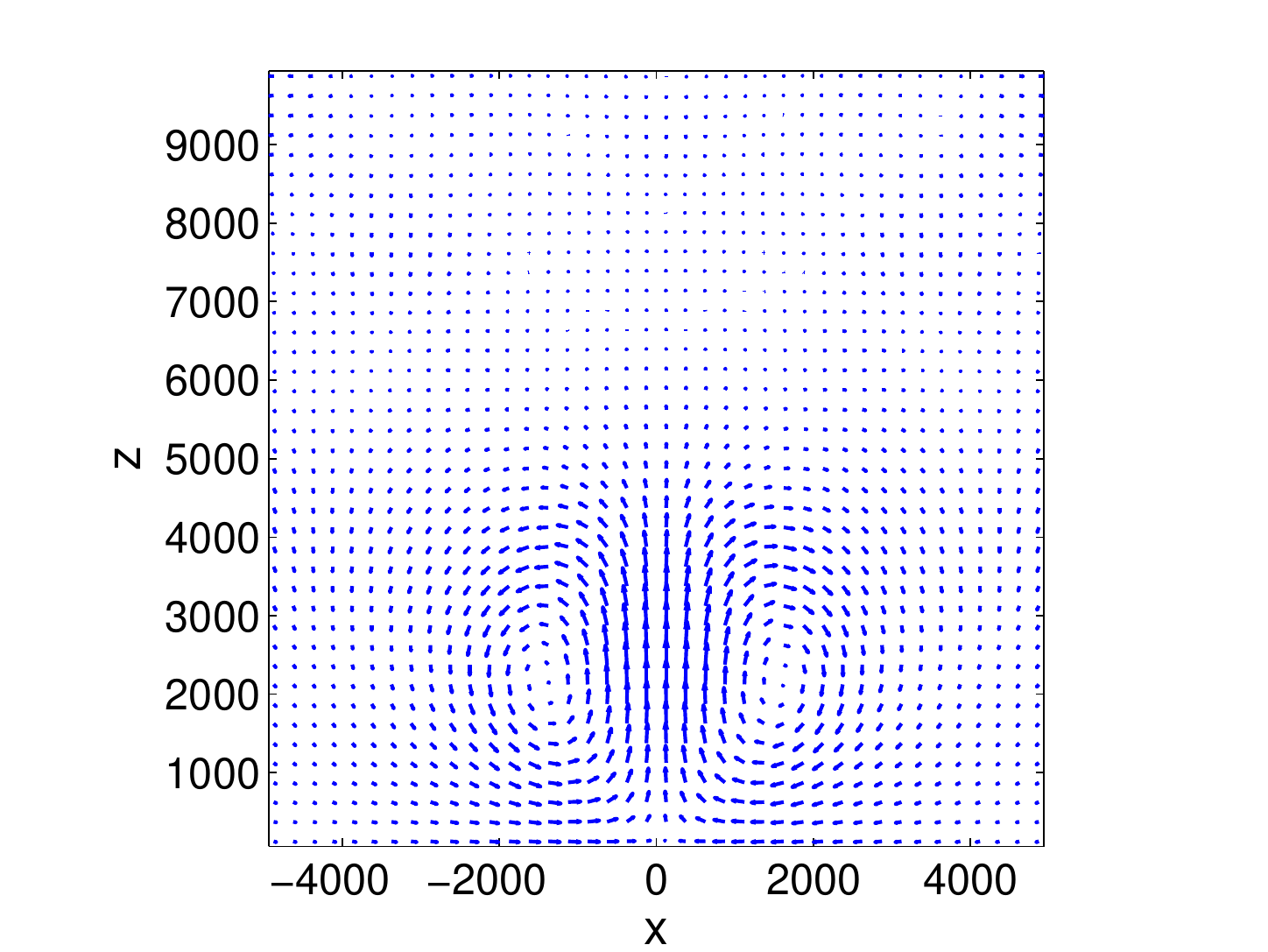}
\includegraphics[width=\textwidth,height=0.25\textheight]{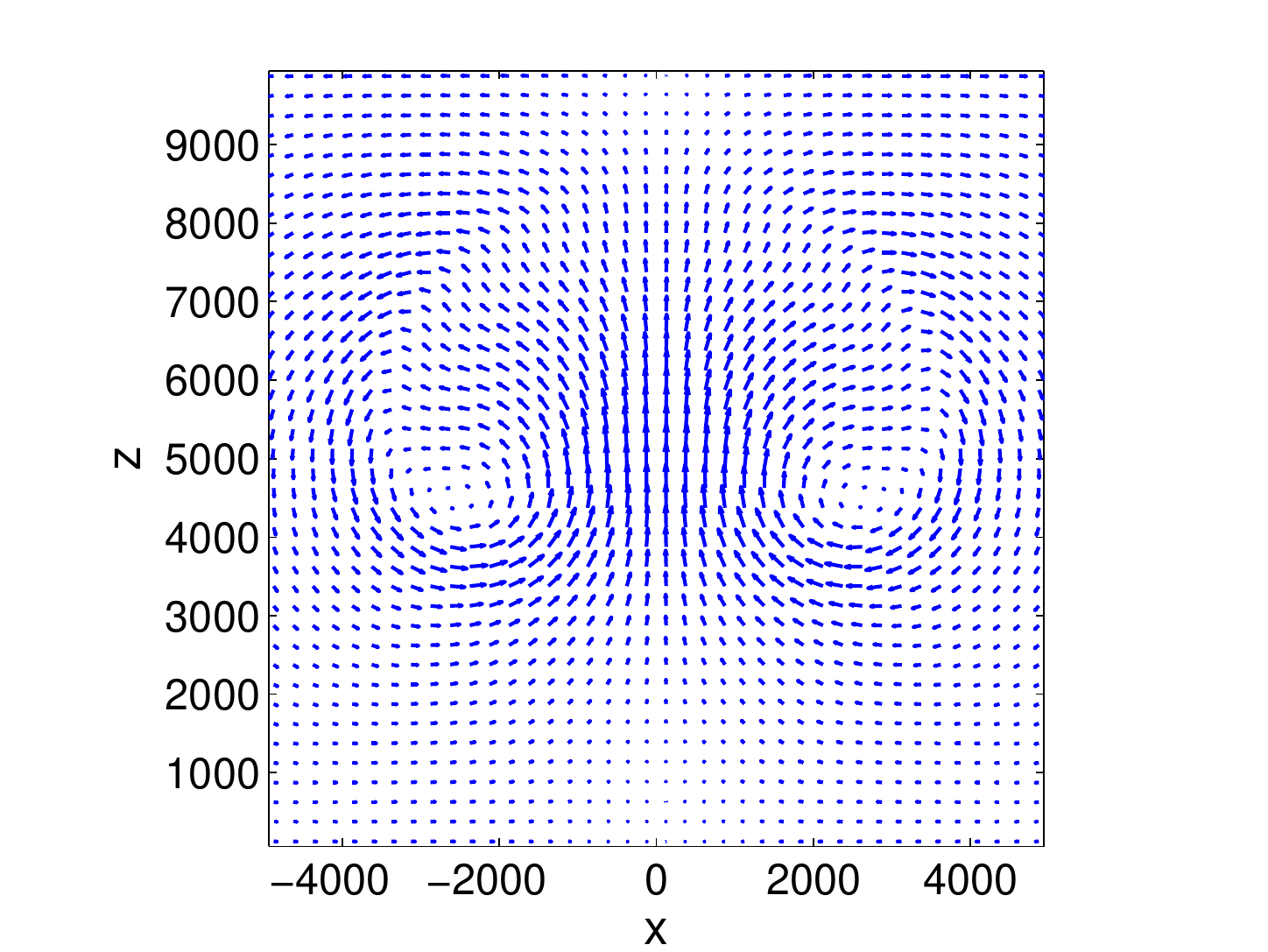}
\end{minipage}
\label{fig:figure4}
\end{figure}

\begin{figure}[H]
\caption{Total energy of the system for the convective bubble test case.}
\centering
\includegraphics[width=0.8\textwidth,height=0.17\textheight]{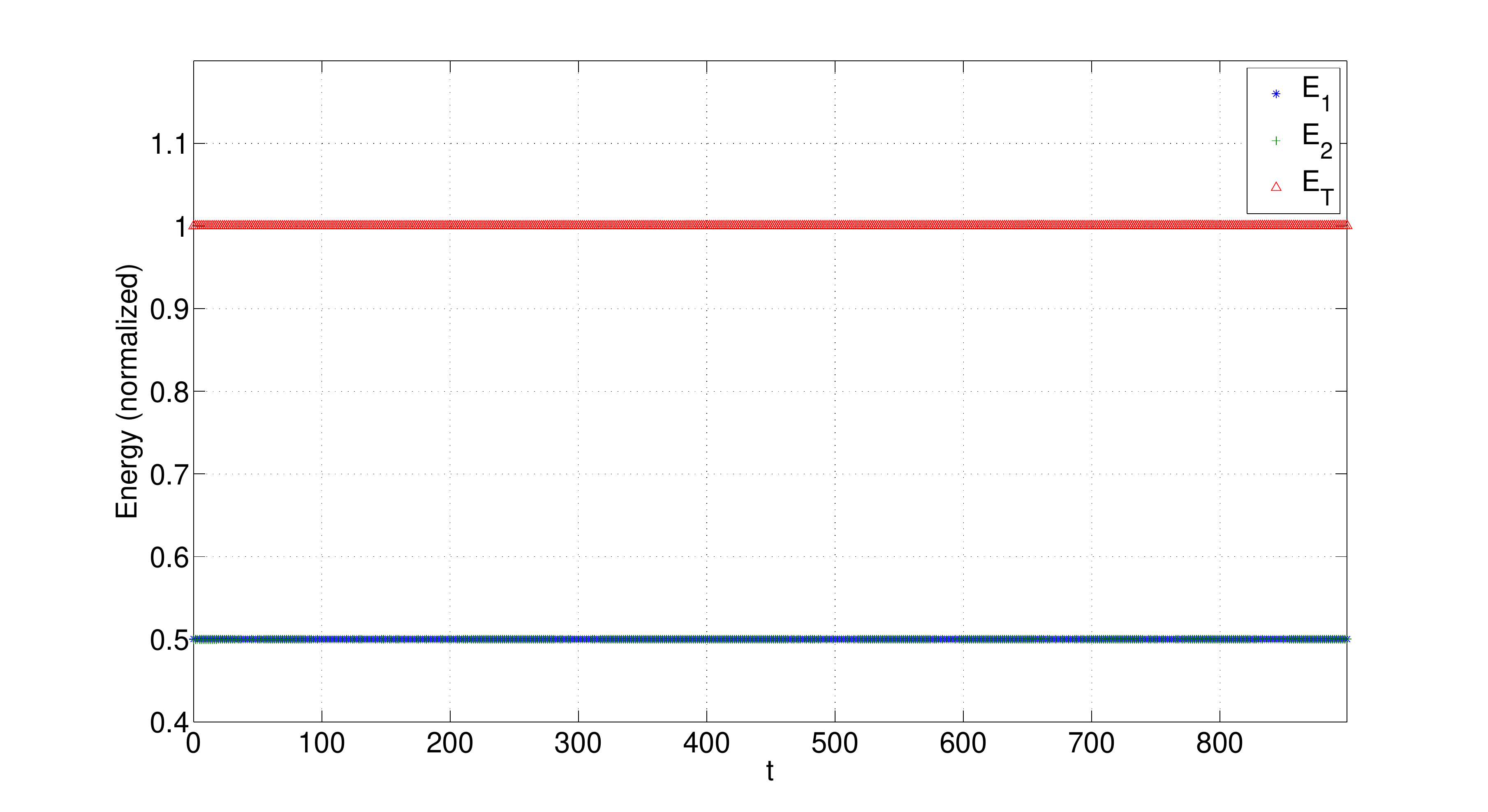}
\includegraphics[width=0.8\textwidth,height=0.17\textheight]{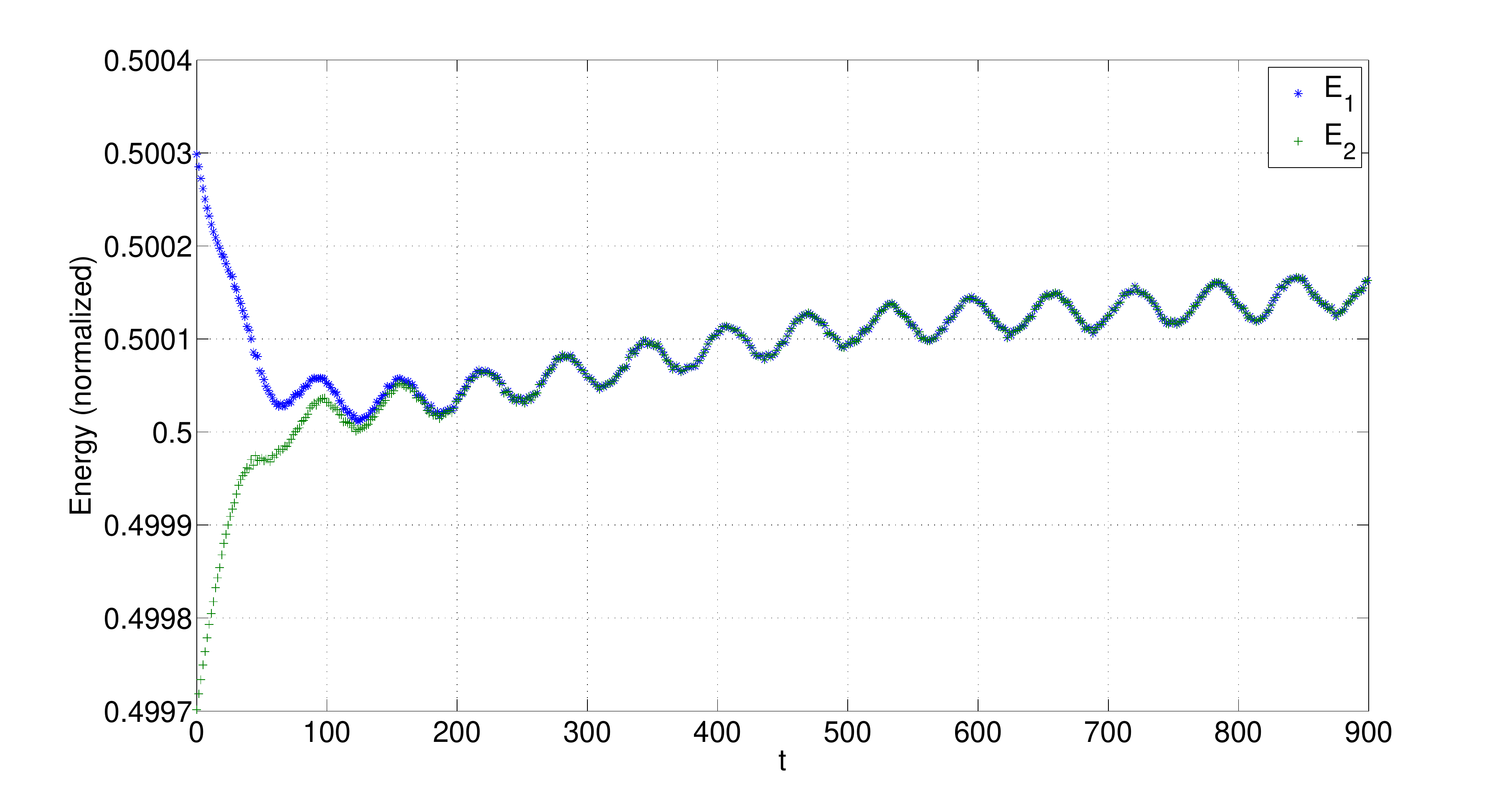}
\label{fig:figure5}
\end{figure}

This first test case, that has been previously studied in \cite{robert,smolarbubble,giraldo} (among others),
is performed by using the scales and settings as in \cite{knothjcp}.
The domain is $\Omega=[-10000,\,10000]\times[-10000,\,10000]\times[0,\,10000]\times[0,\,1000]$,
with $\dx=\dz=125[m]$ and $\dt$ chosen according to,
\begin{equation}
\dt=CFL\frac{\dx}{\displaystyle{\max_{i,j}}{\,\,(|vel-c_s|,|vel+c_s)|}},
\end{equation}
where the $CFL$ number is set to $0.4$, and $vel=\sqrt{u^2+w^2}$ is
the module of the $x-z$ velocity field (this setting of the time step will be kept for
the remaining tests). Both layers are initialized at rest, and we consider the
following potential temperature profiles:
\begin{equation}
\theta_1=\bar\theta_1+\theta'_1,\quad\theta_2=\bar\theta_2,
\end{equation}
with,
\begin{equation}
\bar\theta_1=\bar\theta_2=300[K],\quad
\theta'_1=\begin{cases}10\cos\left(\frac{\pi L}{2}\right) & L\leq 1,\\ 0 & i.o.c.\end{cases},\quad L=\frac{1}{2000}\sqrt{x^2+(z-2000)^2}.
\end{equation}
The reference states $\bar\theta_1$ and $\bar\theta_2$ are in hydrostatic balance at their respective layers, and therefore the density $\rho$ is initialized via the hydrostatic balance law,
\begin{equation}
c_p\theta\frac{d\pi}{dz}=-g,
\end{equation}
together with the equation of state for an ideal gas,
\begin{equation}
\rho=\frac{P_0}{R_d\theta}\pi^{\frac{c_v}{R_d}},
\end{equation}
both evaluated at the corresponding reference state $\theta=\bar\theta$.
We implement solid wall boundary conditions along the domain.
The underlying physics of this test dictate that, as the temperature perturbation is warmer than the background state, a buoyancy force will push the perturbation upwards. As it starts rising, it also starts experiencing horizontal expansion because of the same buoyancy effect, eventually generating a mushroom-shaped cloud. This has been observed in many tests based on 2D sets of Euler equations. We observe a similar behavior in the layered model, as shown in figure \ref{fig:figure1}. Figure \ref{fig:figure2} shows the evolution of the potential temperature residual between both layers, which preserves horizontal symmetry but decreases in magnitude; recall that the second layer is initialized without a perturbation and therefore the system tends to an equilibrium state between both layers. Velocity fields in figure \ref{fig:figure3} are consistent with results obtained in 2D tests performed in the references mentioned in the description of the test; vector plots in figure \ref{fig:figure4} illustrate the deformation process to which the temperature perturbation is subjected. Finally, \ref{fig:figure5} exhibits conservation of energy of the total system, and also confirms the equilibrium reached by the layers.

\subsection{Interaction between hot and cold bubbles}

This test is a variation of the test presented in \cite{robert},
which studied the interaction between positive and negative potential temperature
perturbations, modified to keep the same scale of our first test case; all the domain
parameters are set as in test \ref{section:test1}. Horizontal velocity $u$ is initialized in both layers
with a value of $20[ms^{-1}]$ in order to make the test more stringent.
A periodic boundary condition is set in the lateral $x$-direction while solid boundary walls
are kept at the bottom and the top of the domain. Vertical velocity in both layers is set to zero.
Thermodynamic variables are computed in exactly the same way as in the first test,
but we add a cold perturbation to the second layer,
\begin{equation}
\theta'_2=\begin{cases}-15\cos\left(\frac{\pi L}{2}\right) & L\leq 1,\\ 0 & i.o.c.\end{cases},\quad L=\frac{1}{2000}\sqrt{x^2+(z-8000)^2}.
\end{equation}
\begin{figure}[H]
\caption{Initial potential temperatures for the cold and hot bubbles test case ($\theta_1$ on the left, $\theta_2$ on the right).}
\begin{minipage}[b]{0.5\linewidth}
\centering
\includegraphics[scale=0.45]{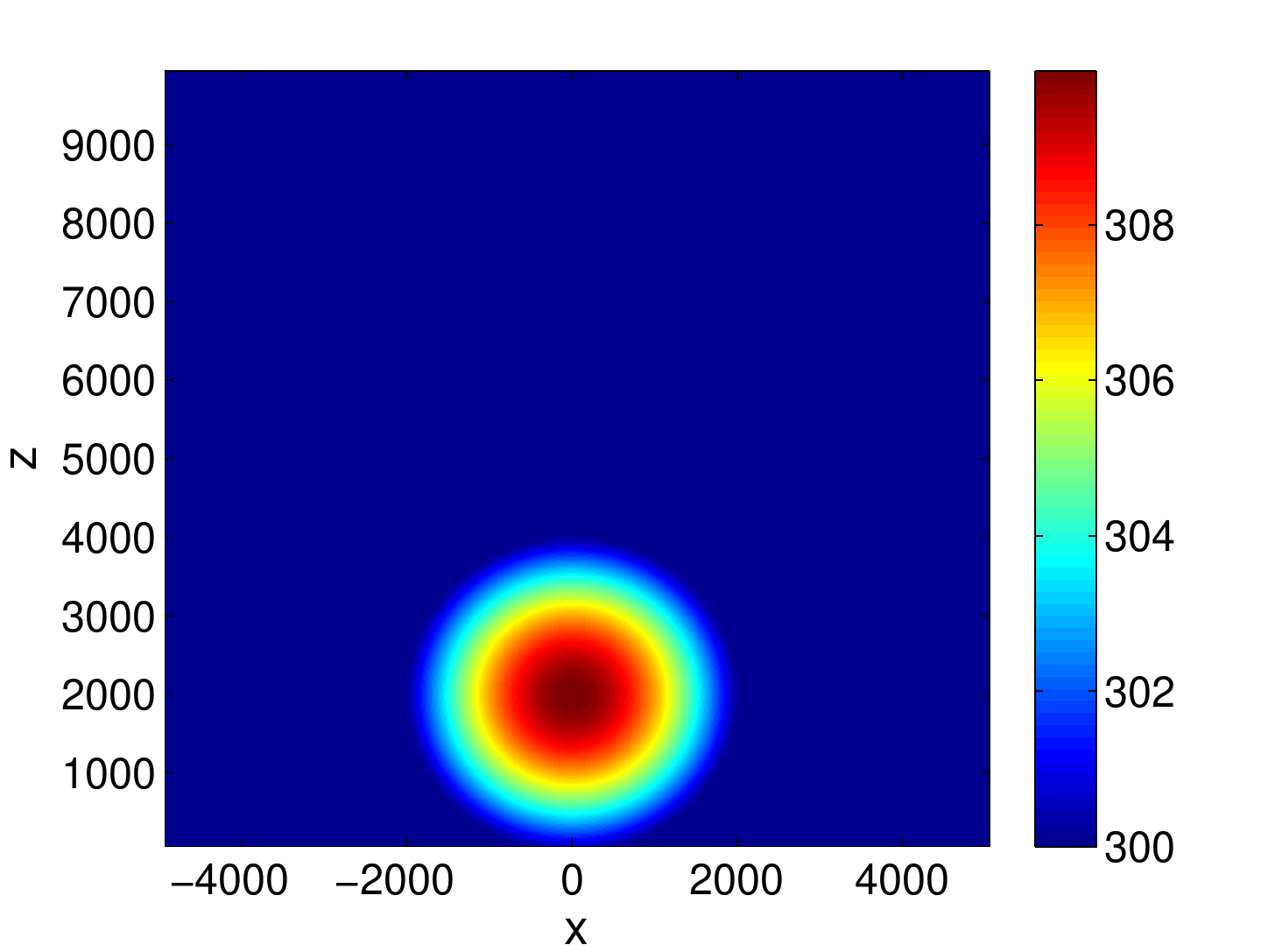}
\end{minipage}
\begin{minipage}[b]{0.5\linewidth}
\centering
\includegraphics[scale=0.45]{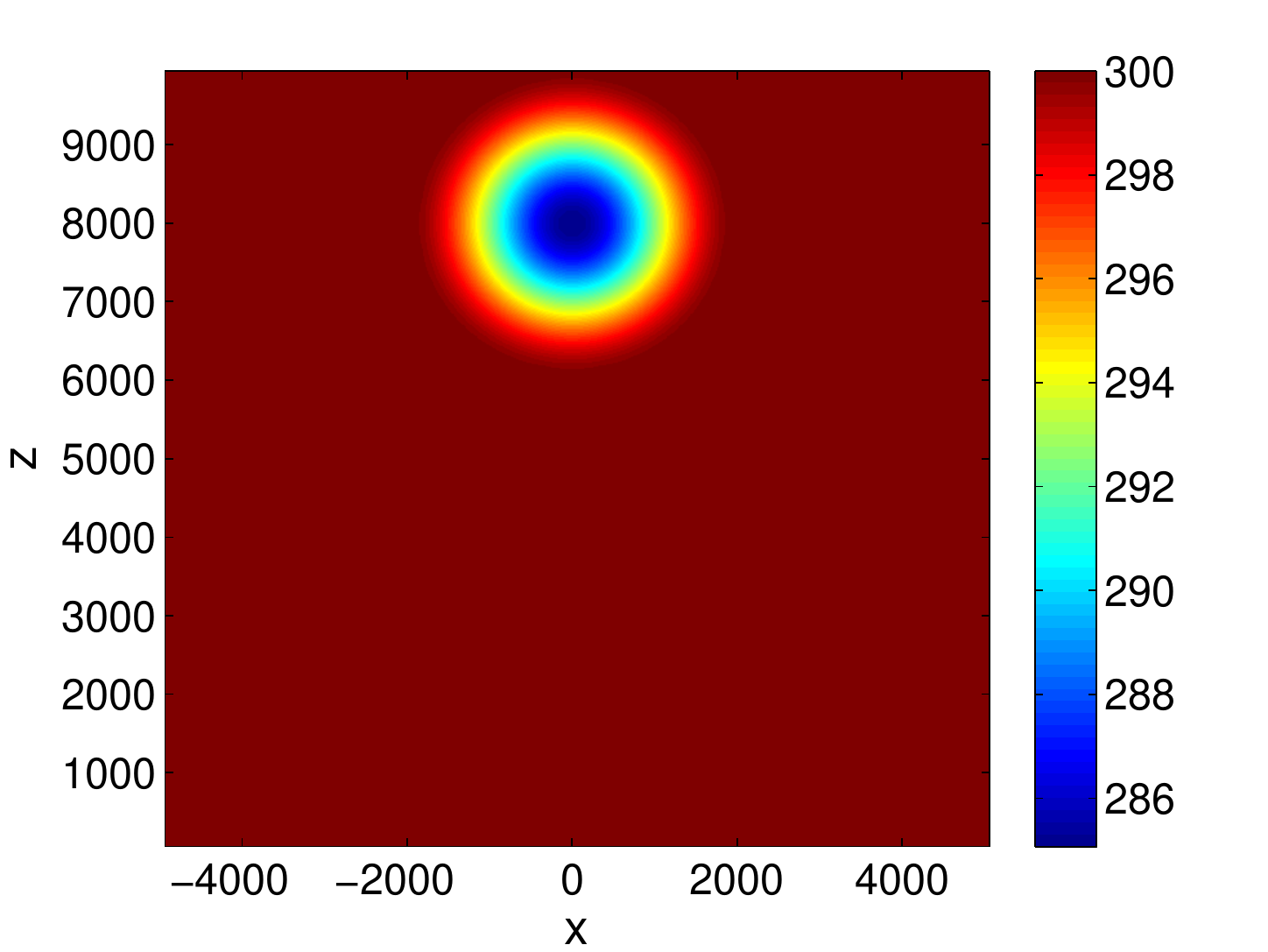}
\end{minipage}
\label{fig:figure6}
\end{figure}
As we have previously described in the first test case, in a neutral atmosphere, a warm temperature perturbation is expected to raise; by the same principles, a cold perturbation is expected to fall. Figure \ref{fig:figure6} shows the initial temperature for both layers. Figure \ref{fig:figure7} illustrates the experiment: placed in the same vertical axis, the warm bubble raises while the cold bubble falls, and at certain instant both bubbles collide starting an eddy interaction which is governed by the same buoyancy effects. In this experiment we also included horizontal velocity and we can observe that the perturbations are horizontally translated with a proper speed (both layers were initialized at the same velocity) and preserving symmetry all along the experiment. Figure \ref{fig:figure8} shows a decrease in the magnitude of the potential temperature residual and figures \ref{fig:figure9} and \ref{fig:figure10} exhibit the associated velocity field, where we can observe the formation of eddies. Figure \ref{fig:figure11} illustrates the evolution of potential temperature extreme values in every layer, and it can be seen that for both maximum and minimum, the layers have a tendency to reach equilibrium values. Total energy is preserved as shown in figure \ref{fig:figure12} in a similar manner as in the first test.

\begin{figure}[H]
\caption{Potential temperature for the hot and cold bubbles test case. Colormap of the first layer at $t=120$, 300, 600 and 1000[s]. $\dx=\dz=125[m]$, $160\times 80$ elements.}
\begin{minipage}[b]{0.5\linewidth}
\centering
\includegraphics[scale=0.45]{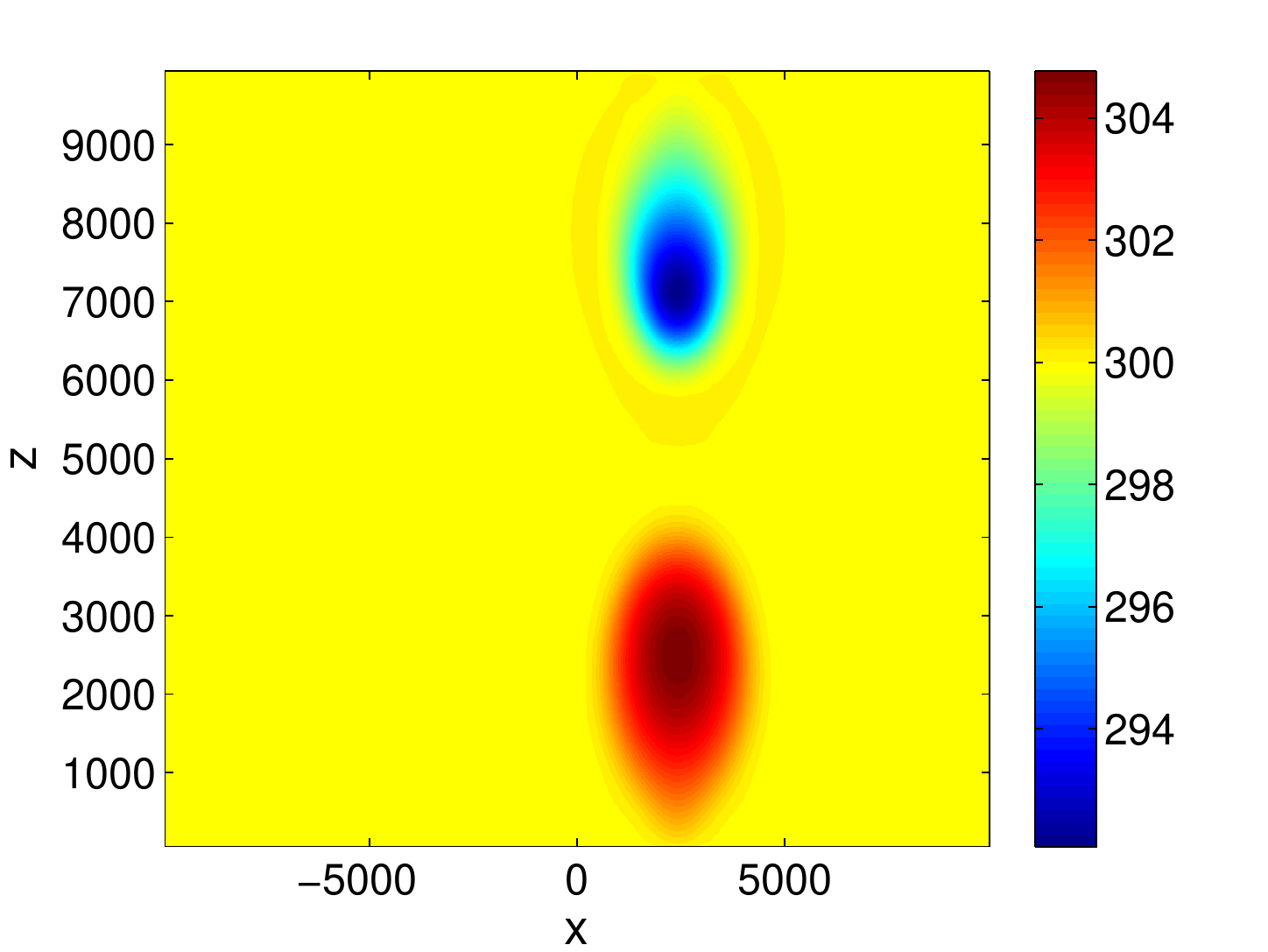}
\includegraphics[scale=0.45]{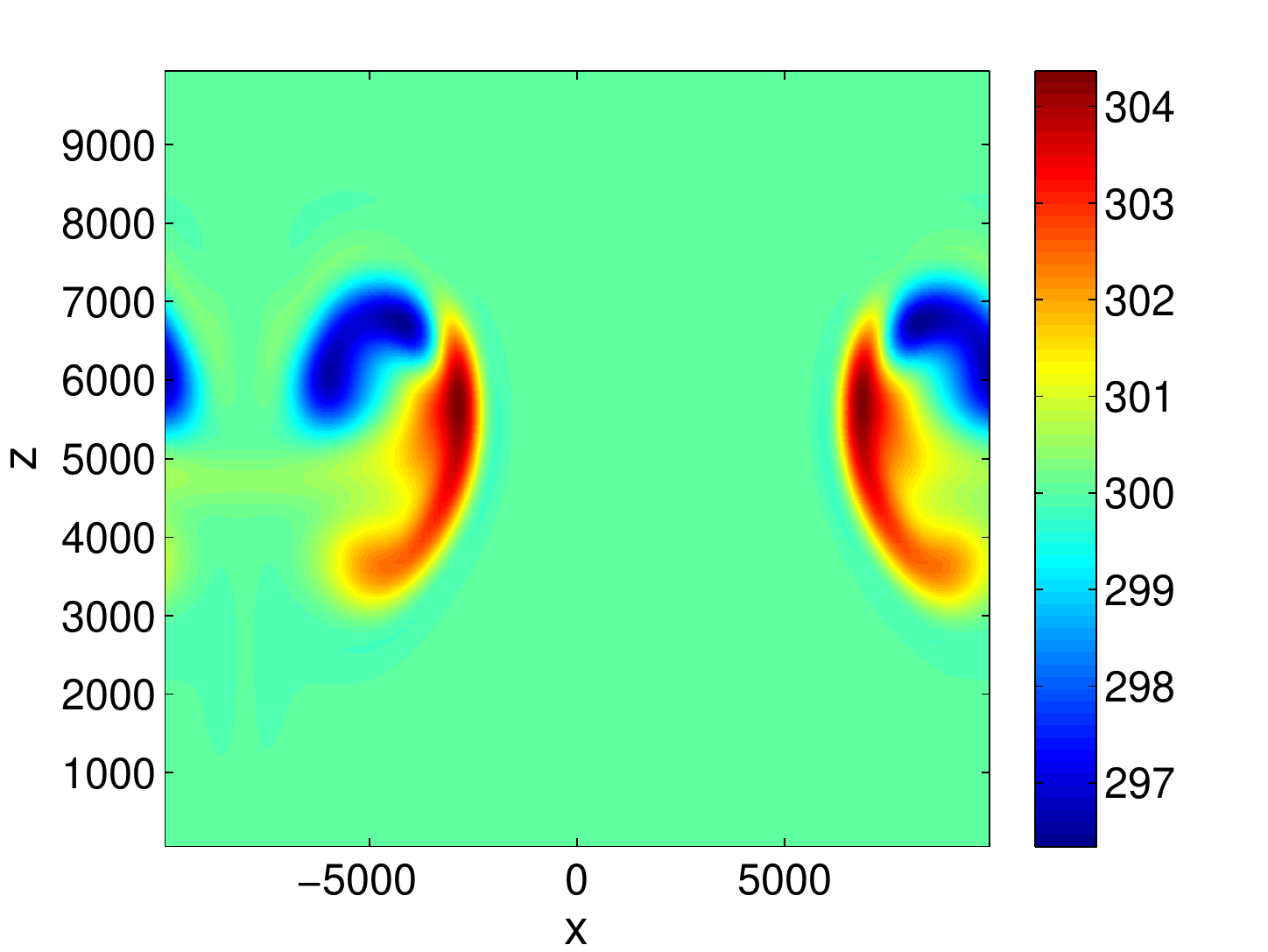}
\end{minipage}
\hspace{0.0cm}
\begin{minipage}[b]{0.5\linewidth}
\centering
\includegraphics[scale=0.45]{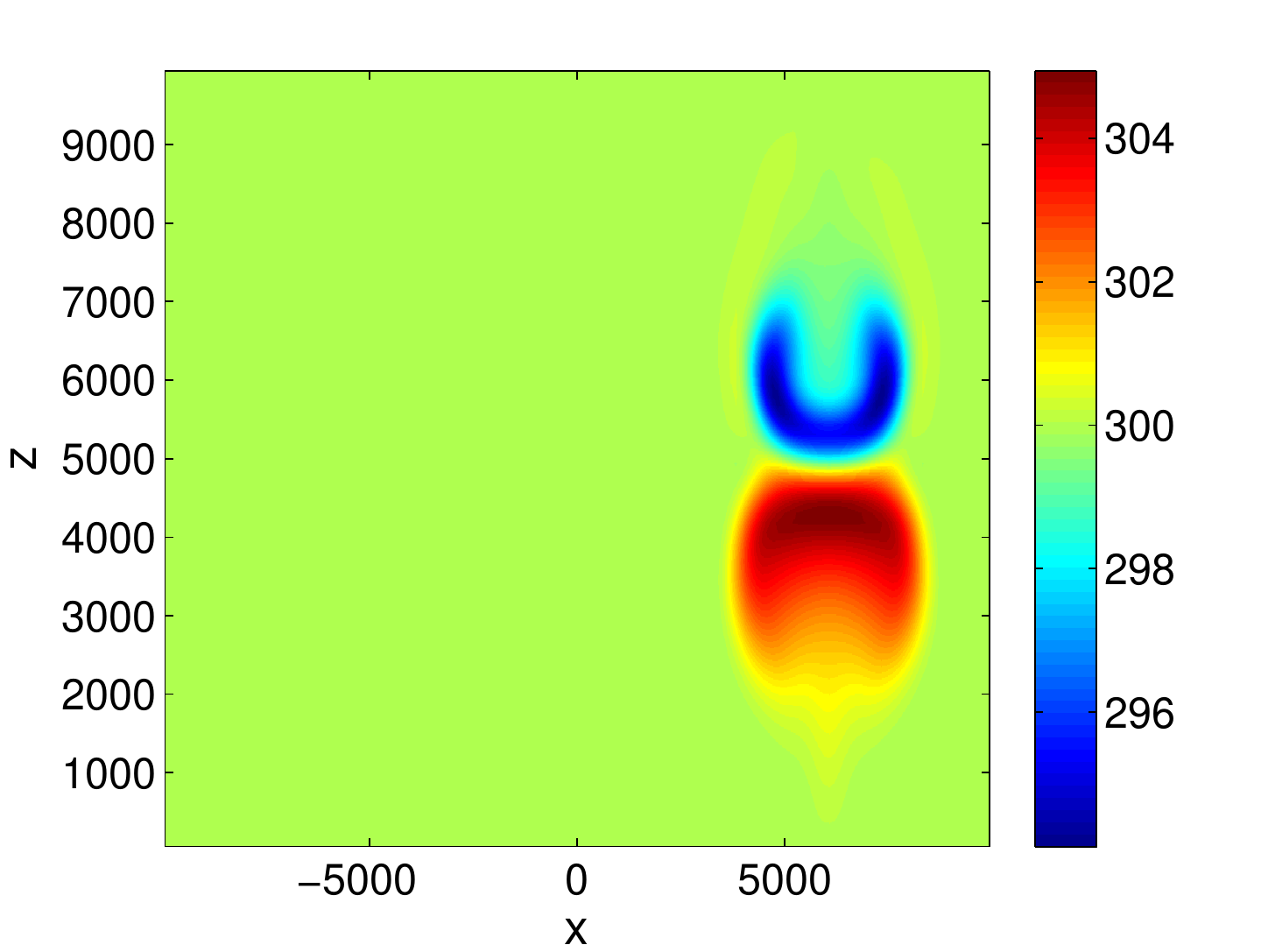}
\includegraphics[scale=0.45]{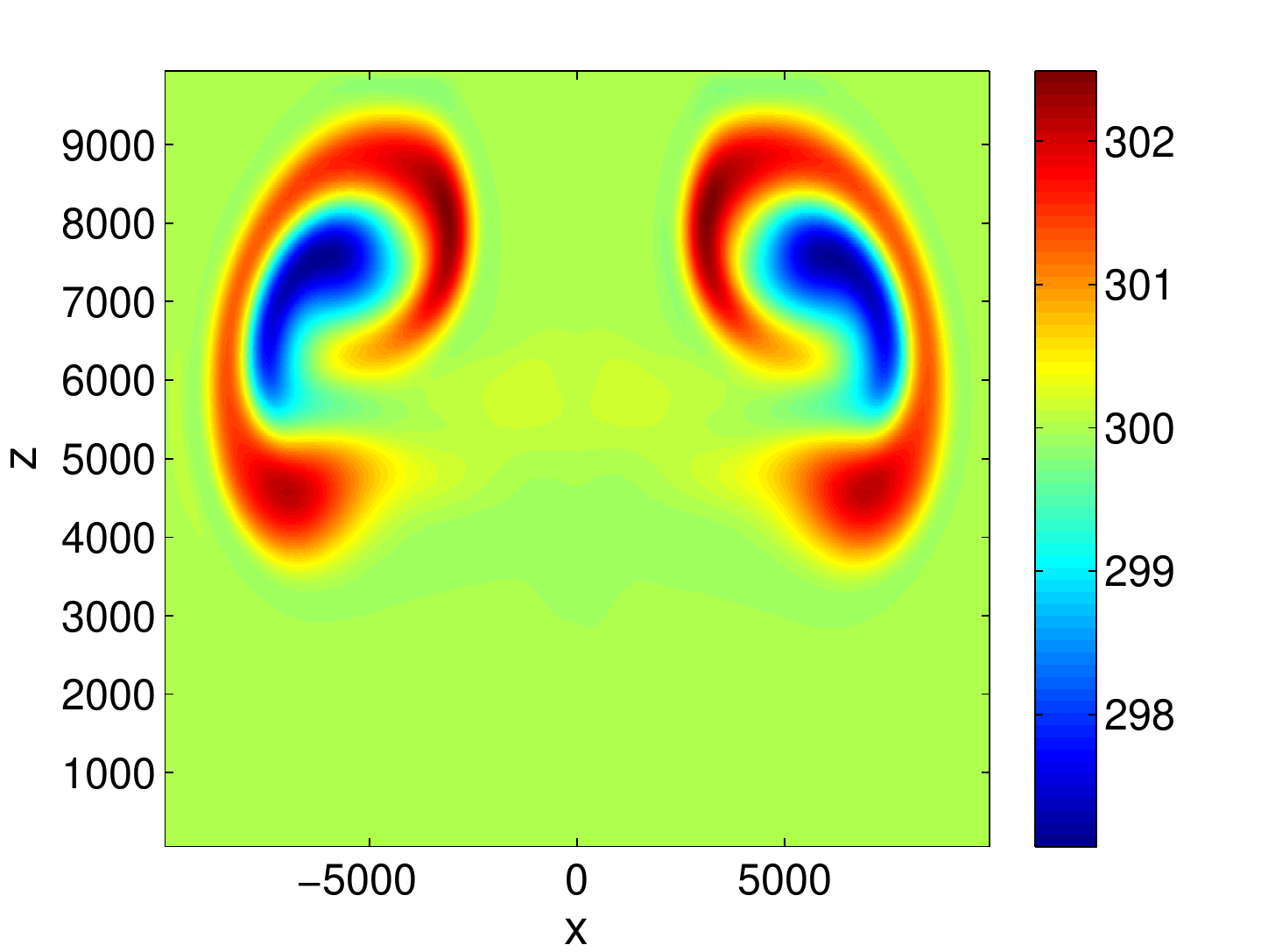}
\end{minipage}
\label{fig:figure7}
\end{figure}

\begin{figure}[H]
\caption{Potential temperature for the hot and cold bubbles test case. Colormap of the residual ($\theta_1-\theta_2$) at $t=120$ and 1000[s]. $\dx=\dz=125[m]$, $160\times 80$ elements.}

\begin{minipage}[b]{0.5\linewidth}
\centering
\includegraphics[scale=0.45]{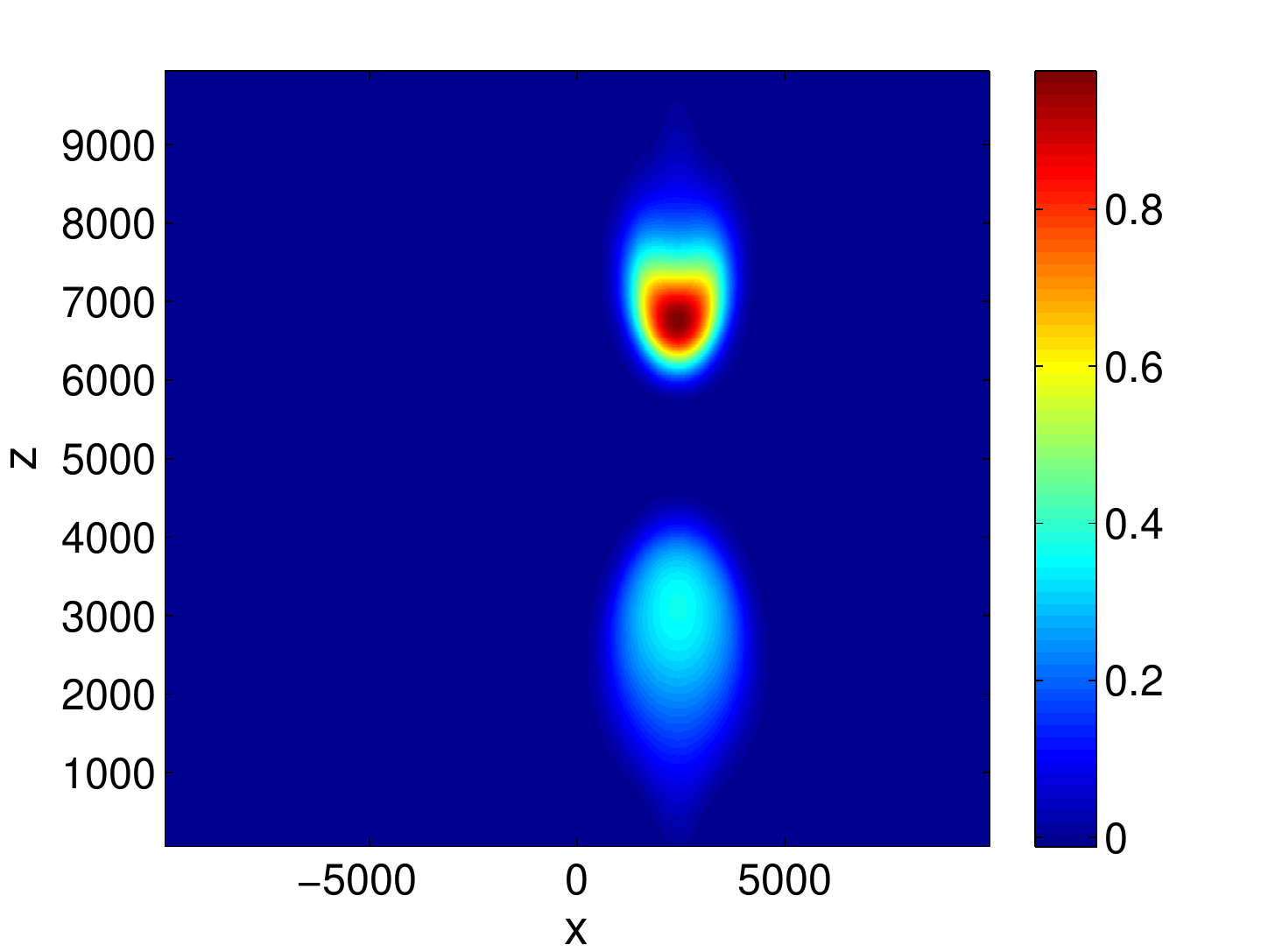}
\end{minipage}
\hspace{0.0cm}
\begin{minipage}[b]{0.5\linewidth}
\centering
\includegraphics[scale=0.45]{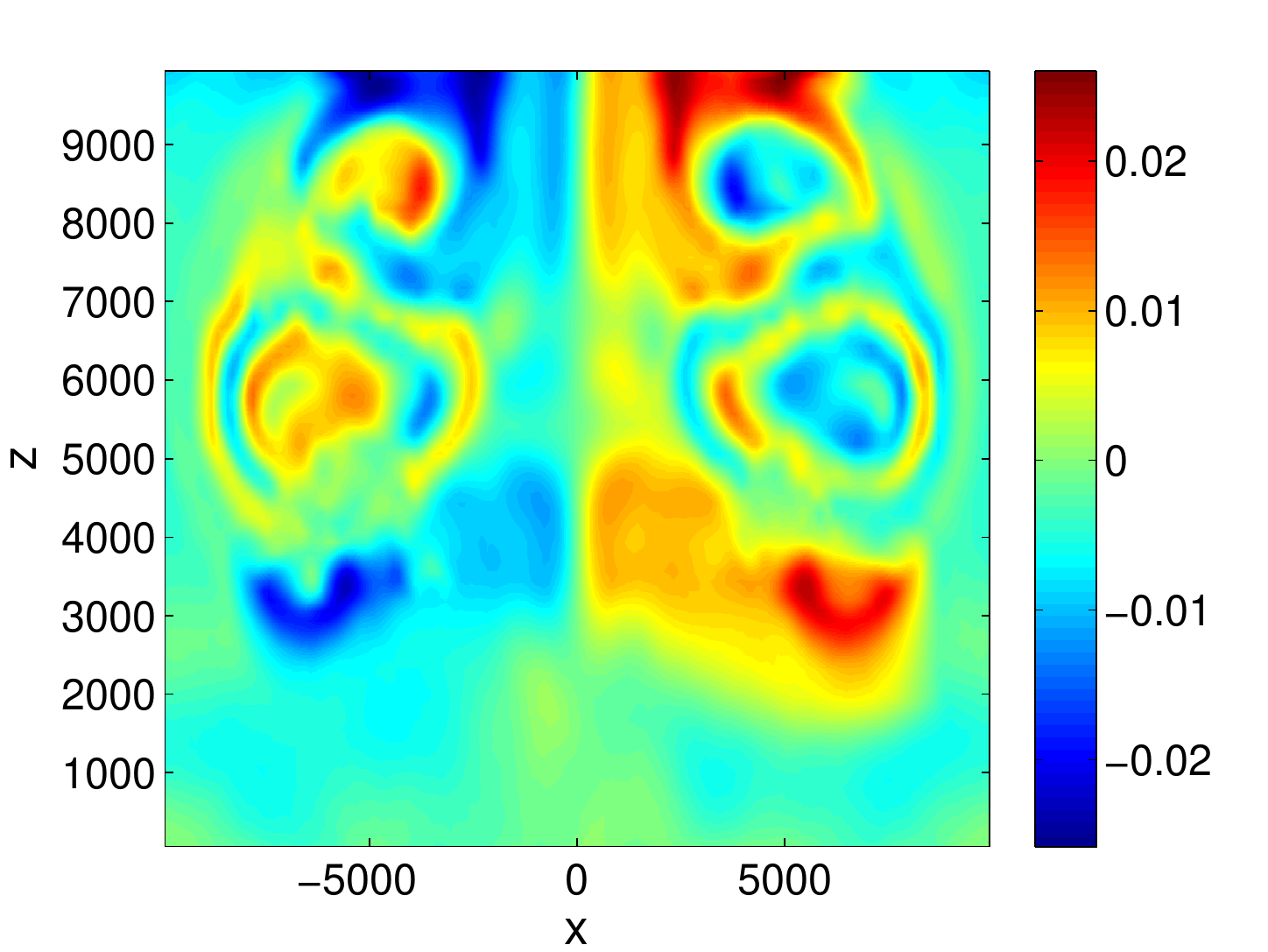}
\end{minipage}
\label{fig:figure8}
\end{figure}

\begin{figure}[H]
\caption{Velocity field vector plots at $t=120$, 300, 600 and 1000[s] for the second layer. $\dx=\dz=125[m]$, $160\times 80$ elements.}
\begin{minipage}[b]{0.5\linewidth}
\centering
\includegraphics[scale=0.45]{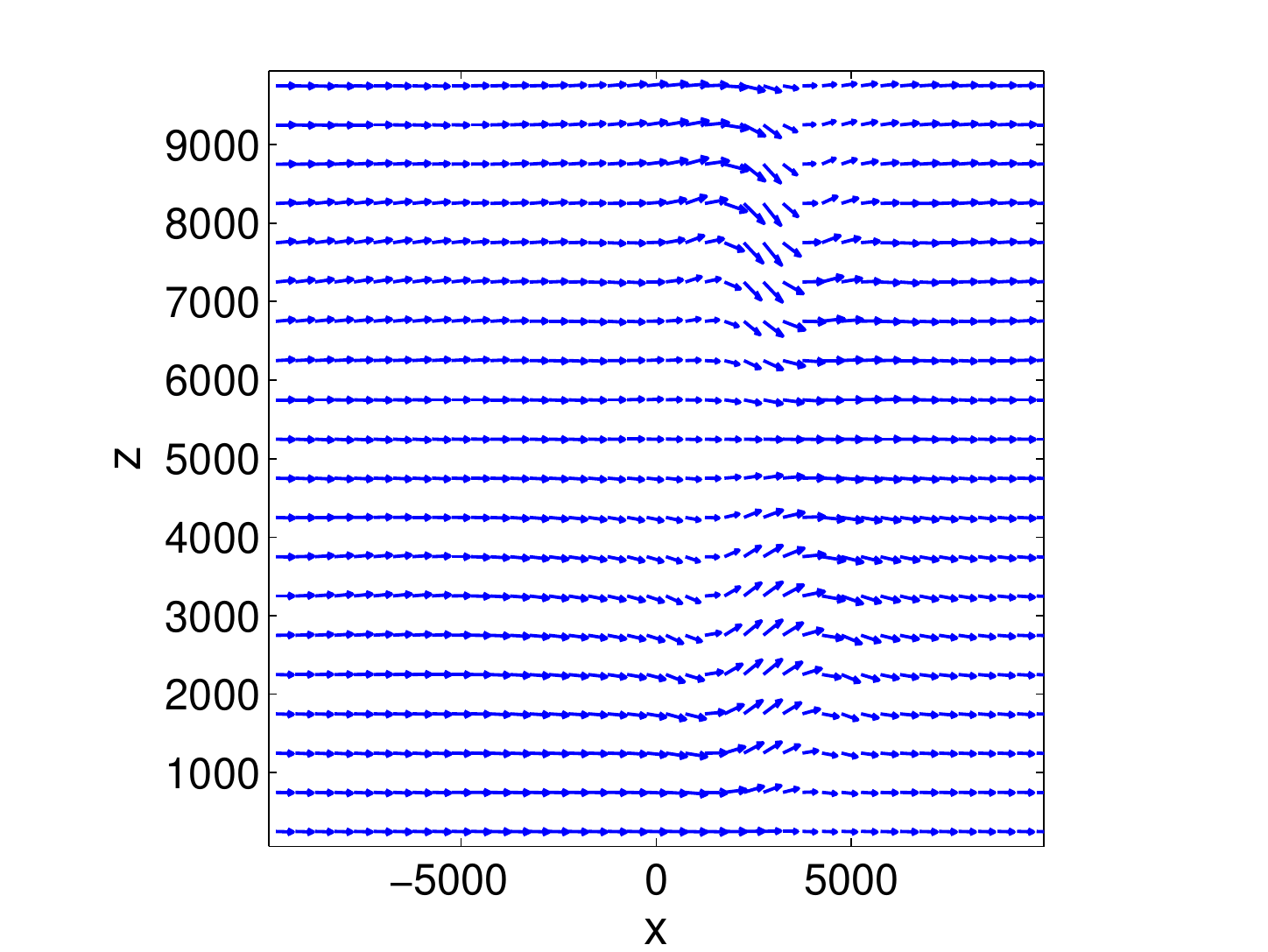}
\includegraphics[scale=0.45]{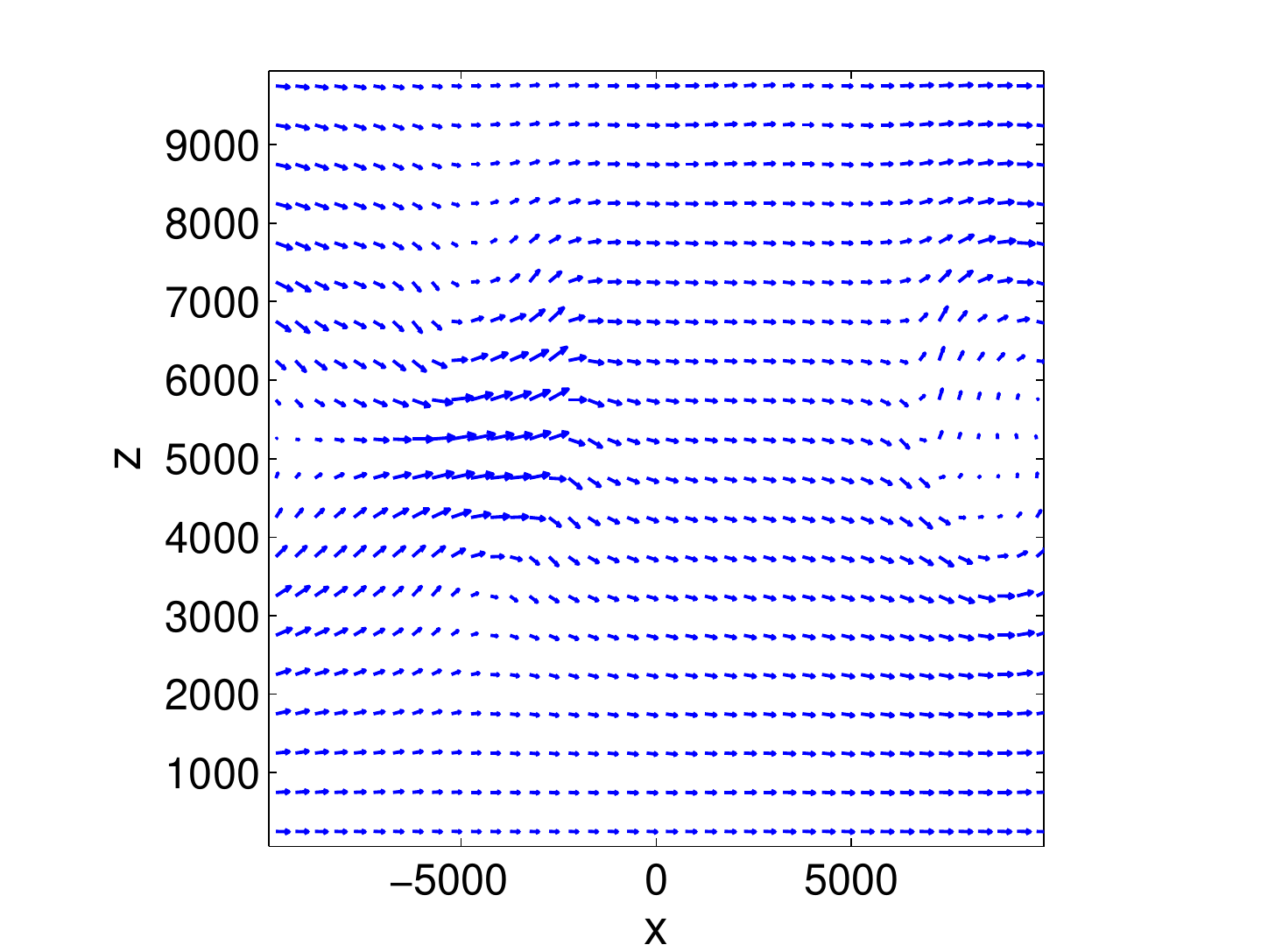}
\end{minipage}
\begin{minipage}[b]{0.5\linewidth}
\centering
\includegraphics[scale=0.45]{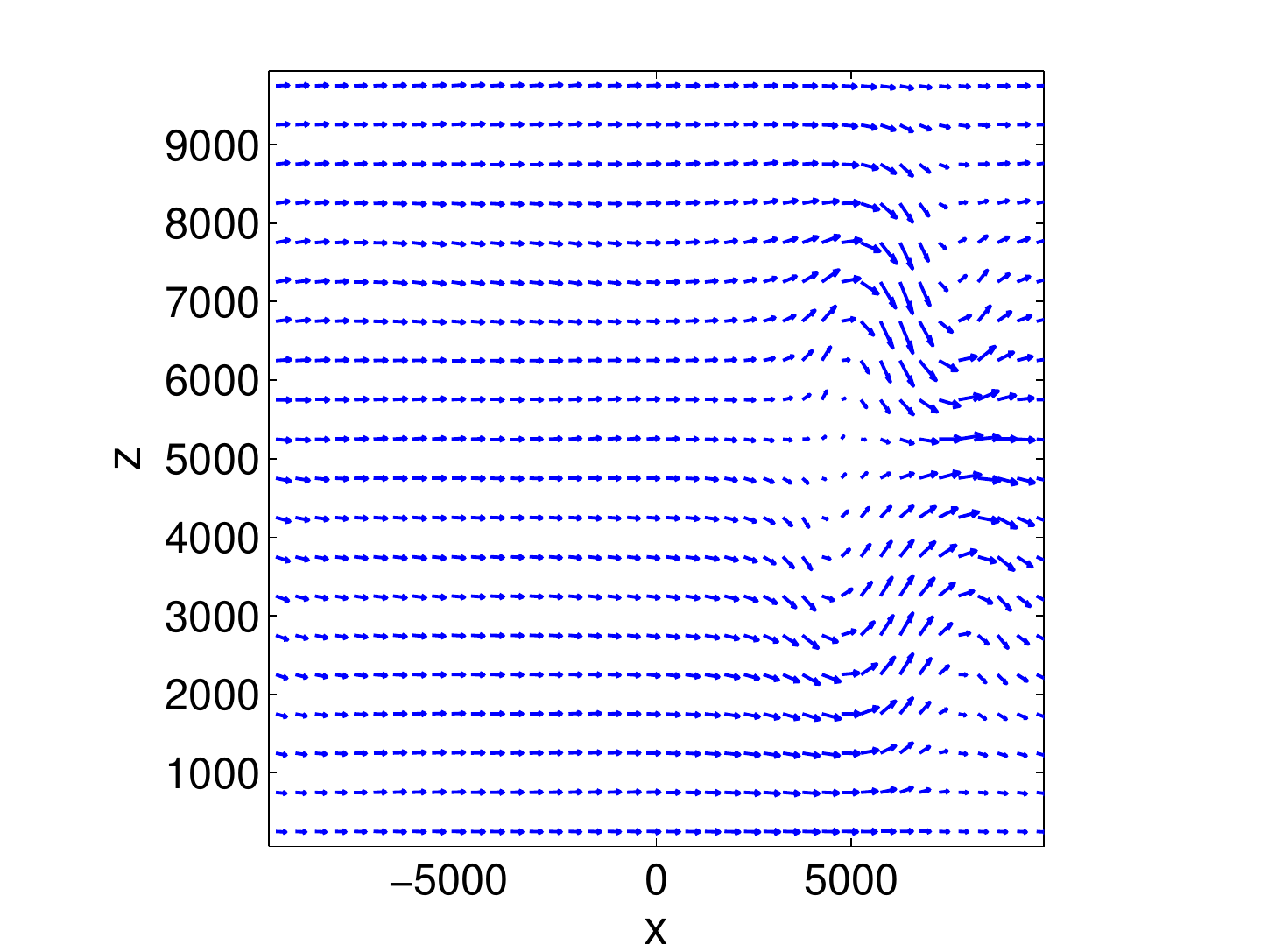}
\includegraphics[scale=0.45]{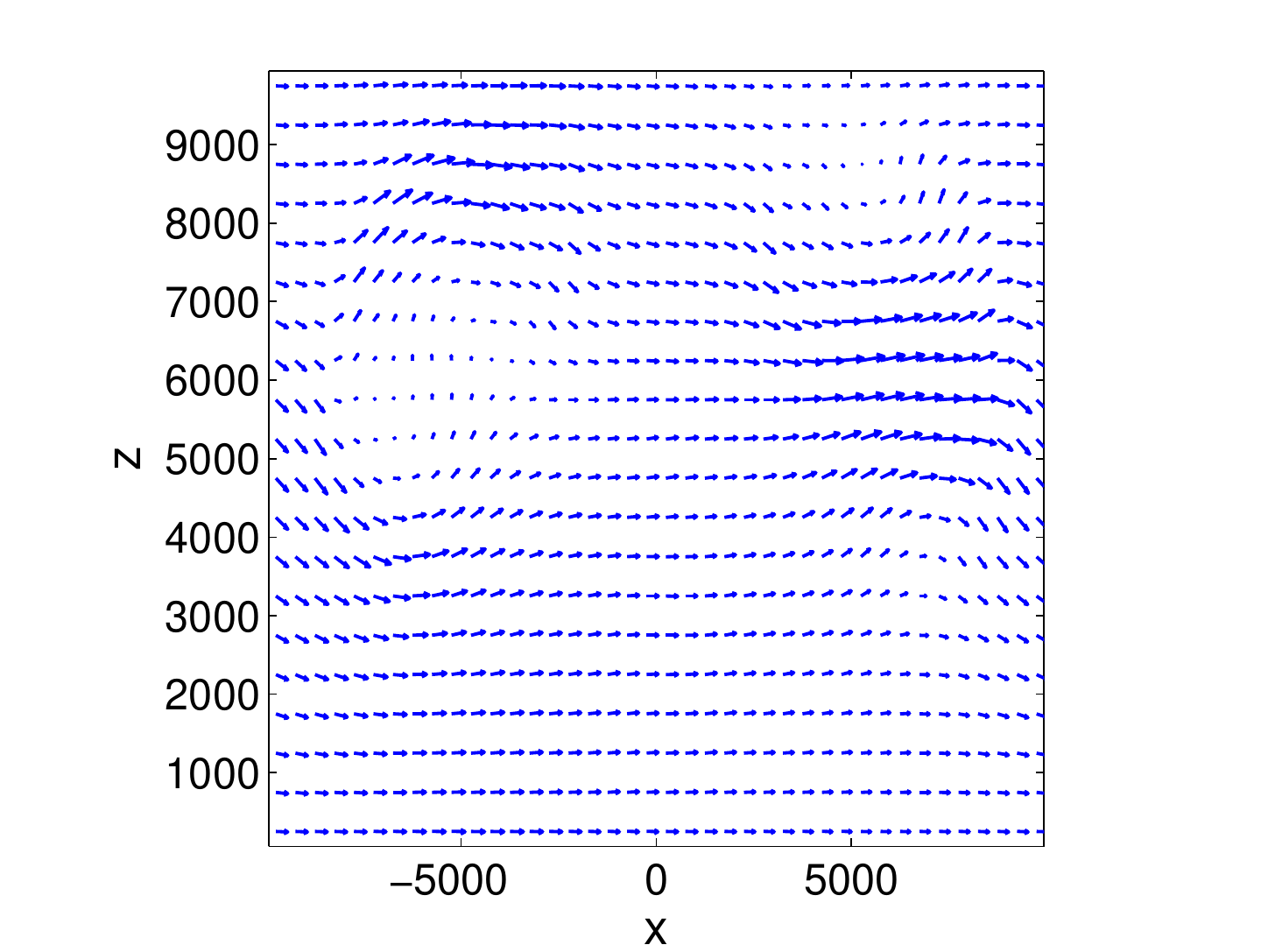}
\end{minipage}
\label{fig:figure10}
\end{figure}

\begin{figure}[H]
\caption{Velocity field colormaps at $t=300$[s] for the second layer. Left: horizontal velocity. Right: vertical velocity. $\dx=\dz=125[m]$, $160\times 80$ elements.}
\begin{minipage}[b]{0.5\linewidth}
\centering
\includegraphics[scale=0.45]{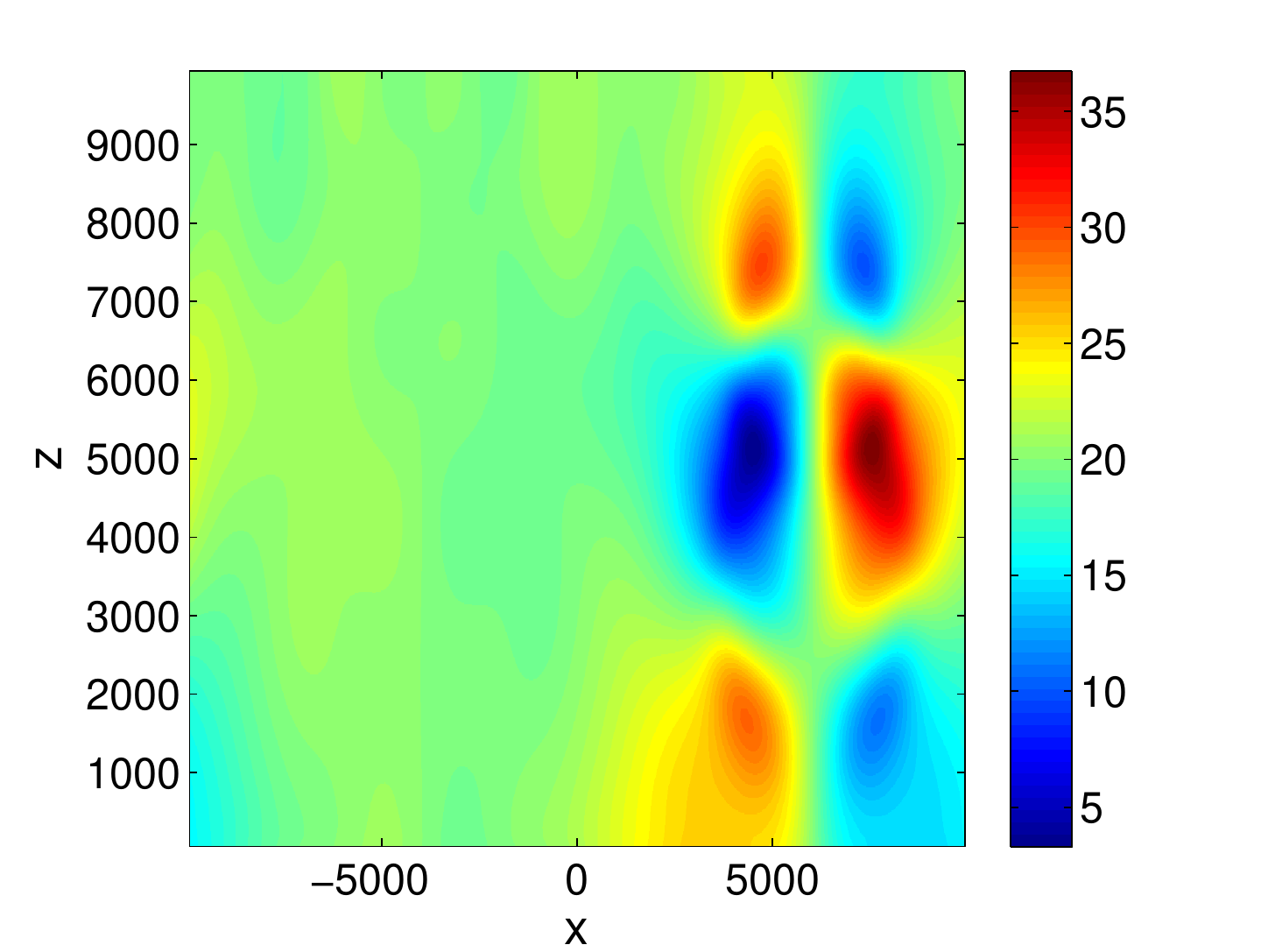}
\end{minipage}
\hspace{0.0cm}
\begin{minipage}[b]{0.5\linewidth}
\centering
\includegraphics[scale=0.45]{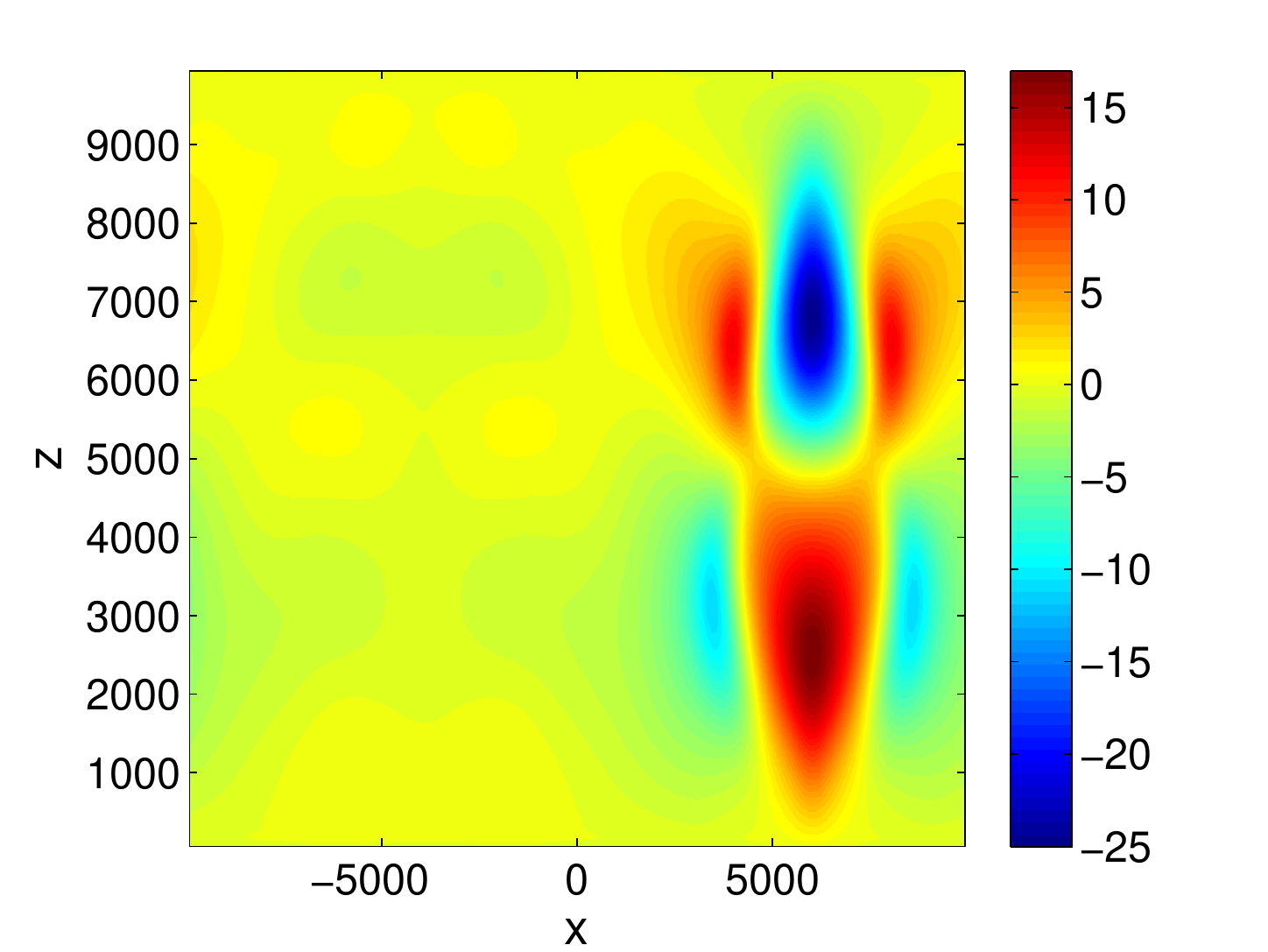}
\end{minipage}
\label{fig:figure9}
\end{figure}

\begin{figure}[H]
\caption{Potential temperature extrema evolution for the hot and cold bubbles test case.}
\centering
\includegraphics[width=0.5\textwidth,height=0.35\textheight]{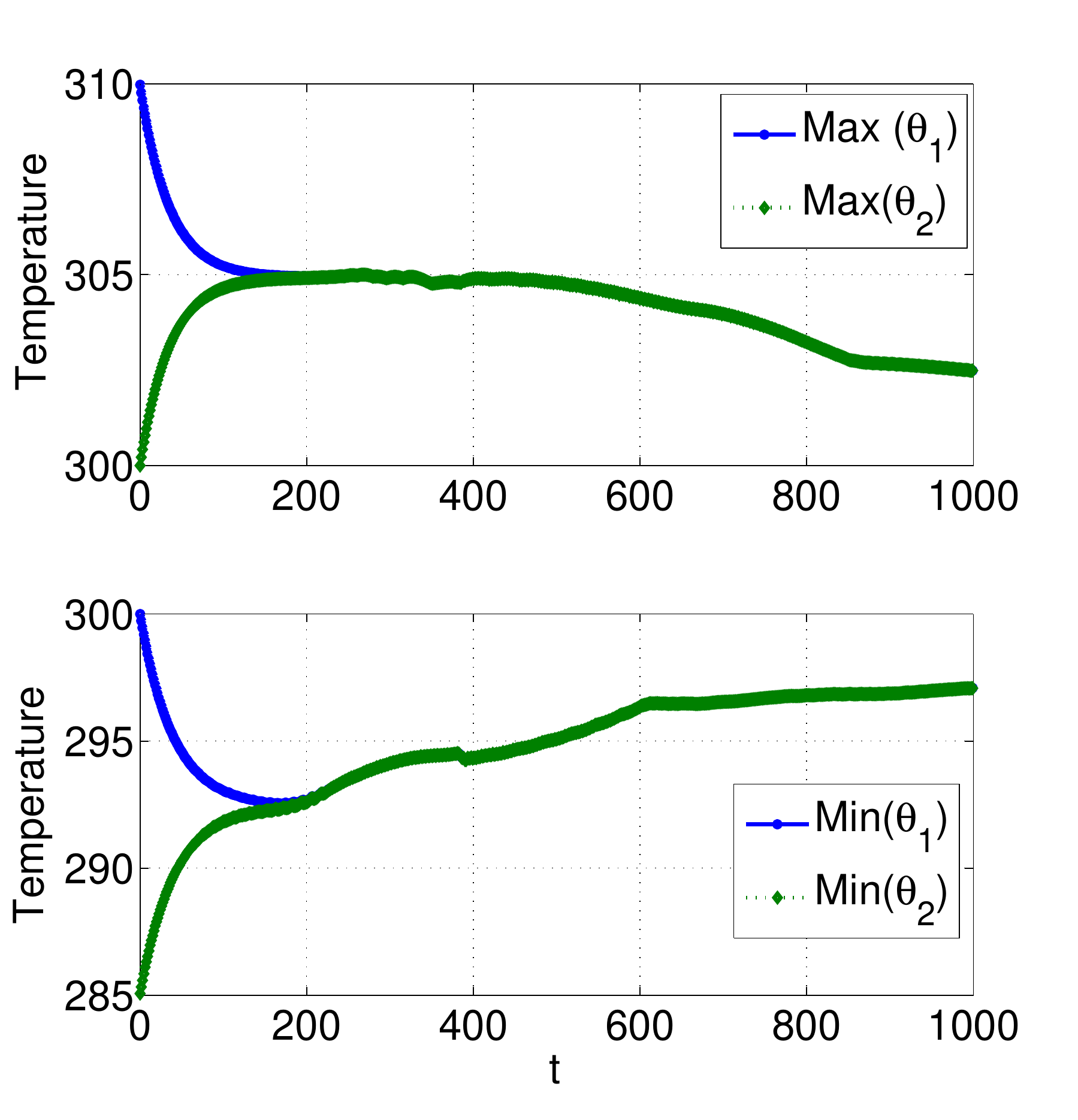}
\label{fig:figure11}
\end{figure}

\begin{figure}[H]
\caption{Total energy of the system for the hot and cold bubbles test case.}
\centering
\includegraphics[width=0.8\textwidth,height=0.17\textheight]{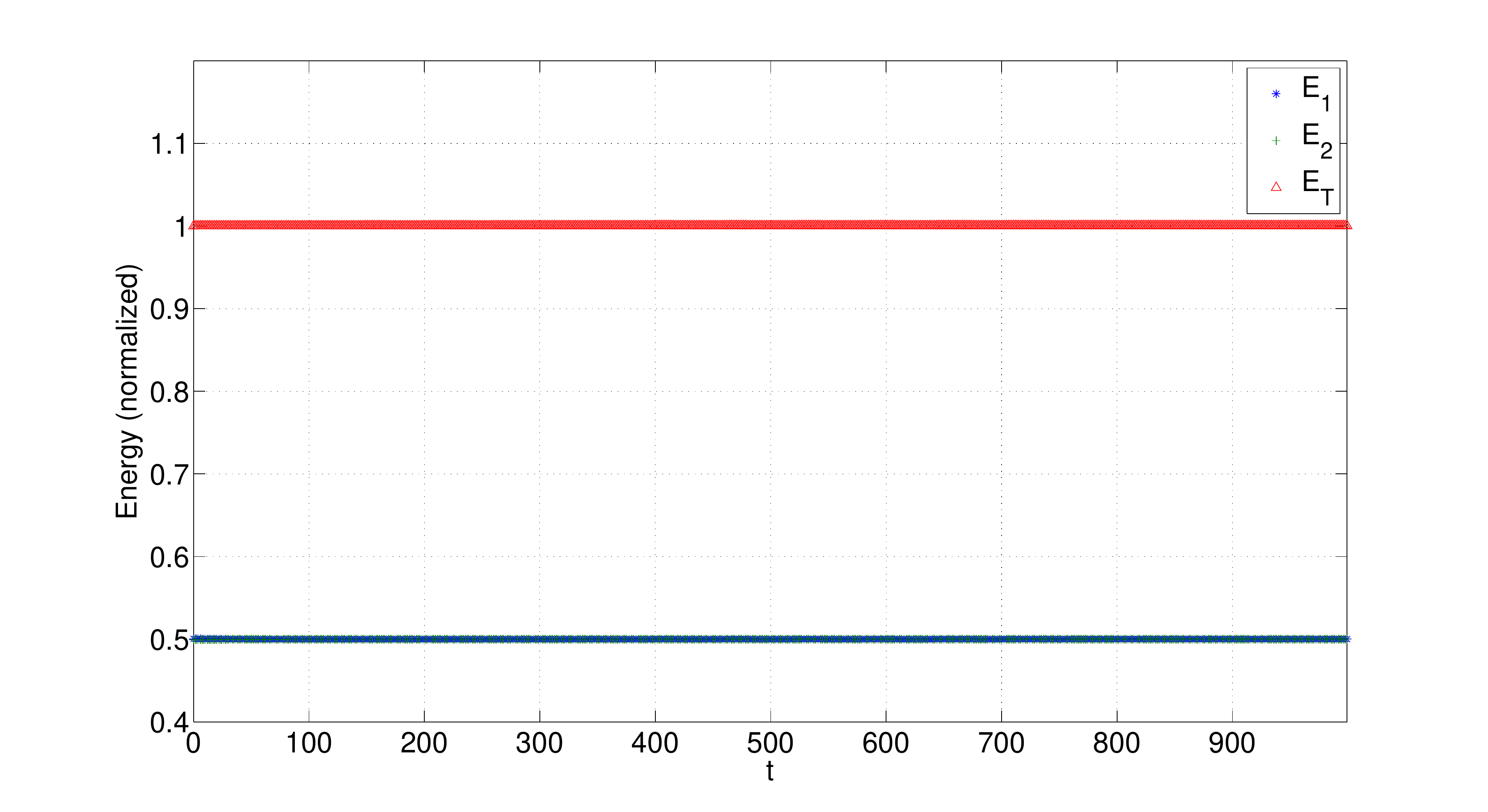}
\includegraphics[width=0.8\textwidth,height=0.17\textheight]{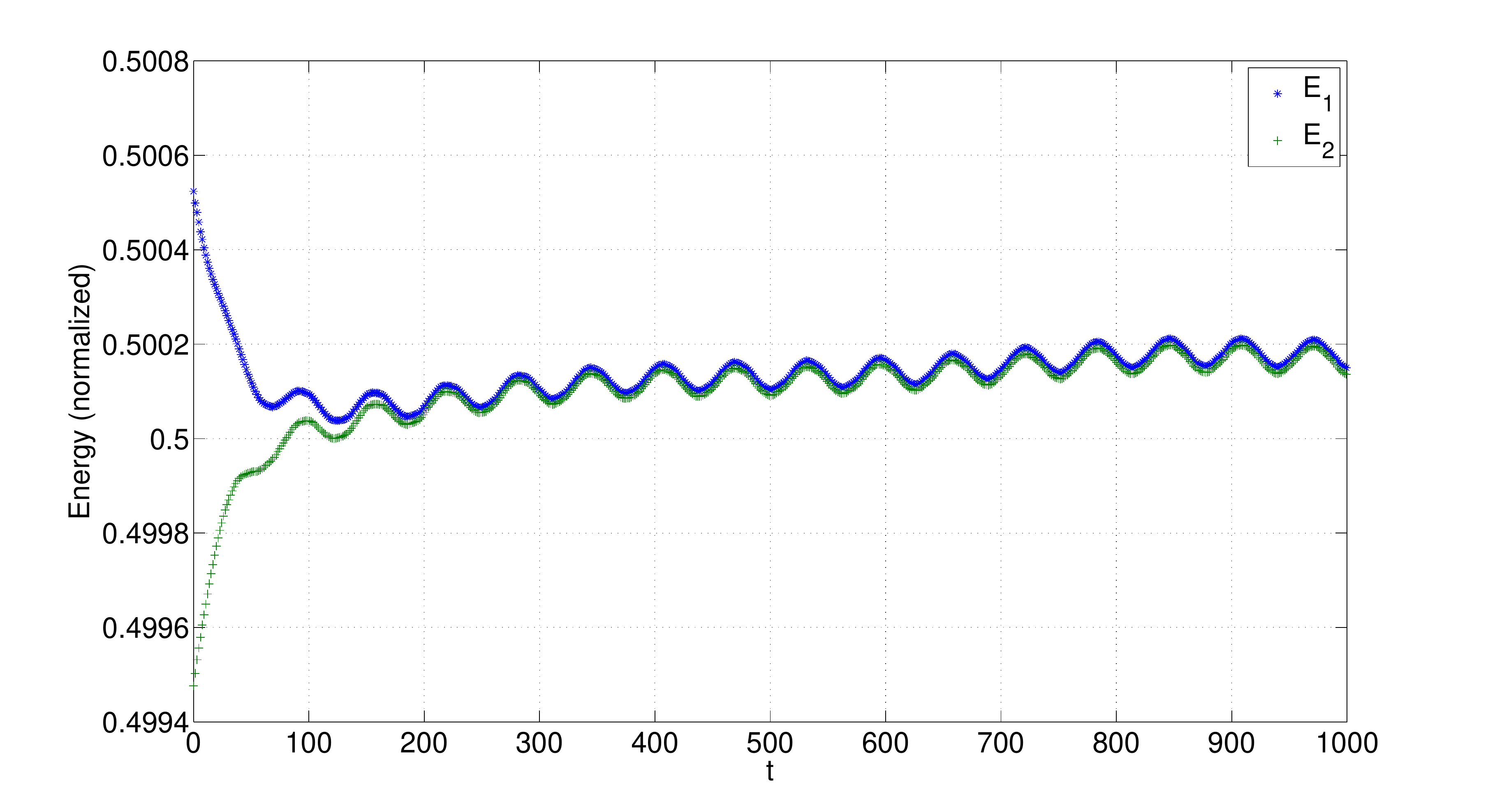}
\label{fig:figure12}
\end{figure}

\subsection{Adjustment of shear between two layers}

This third test case aims to study the level of interaction and adjustment between the
two layers in the model. We again consider the same domain parameters,
except for a final simulation time of $300[s]$. However, variables are initialized in a different way.
First of all, horizontal velocity $v$ is set to a value of $10[ms^{-1}]$ in the first layer,
while in the second layer it has a value of $-10[ms^{-1}]$.
The horizontal velocity $u$ is initialized in the first layer with a logarithmic profile
\begin{equation}
u_1=50\left(\log\left(\frac{z}{L^Z}+1\right)\right)^{1/2}
\end{equation}
and it is set to zero in the second layer. Vertical velocity in both layers is set to zero.
Thermodynamic reference states are different: no potential temperature perturbation
is considered in any case, the first layer having the same neutral atmosphere
as in the previous test cases, while the second is considered to be a stable atmosphere
with
\begin{equation}
\bar\theta_2=\theta_0\exp\left(\frac{\mathcal{N}^2}{g}z\right),
\end{equation}
where $\theta_0=\bar\theta_1=300[K]$ and $\mathcal{N}$ is the Brunt-V\"ais\"al\"a frequency
with a value of $\mathcal{N}=0.01[s^{-1}]$.
In both cases $\rho$ is computed as in the previous tests.
Boundary conditions are set again
to be periodic in the lateral $x$-direction and solid wall boundary conditions in $z$.
\begin{figure}[H]
\caption{Temporal evolution of the horizontal velocity in the $x$-direction $u$, vertical section at $x=0$, in layer 1 (left) and layer 2(right). }
\centering
\includegraphics[width=\textwidth]{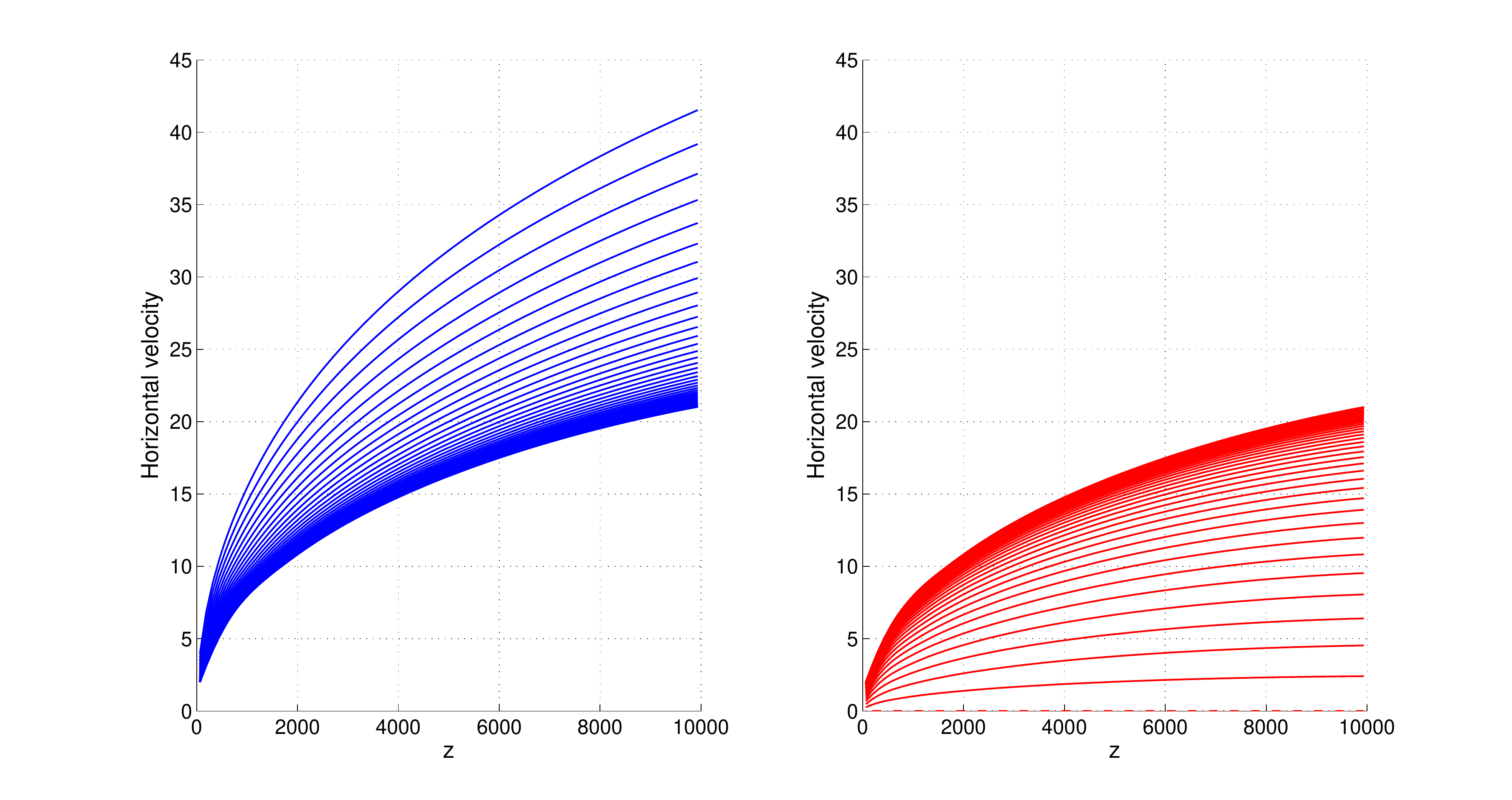}
\label{fig:figure13}
\end{figure}

\begin{figure}[H]
\caption{Temporal evolution of the potential temperature $\theta$, vertical section at $x=0$, in layer 1 (left) and layer 2(right).}
\centering
\includegraphics[width=\textwidth]{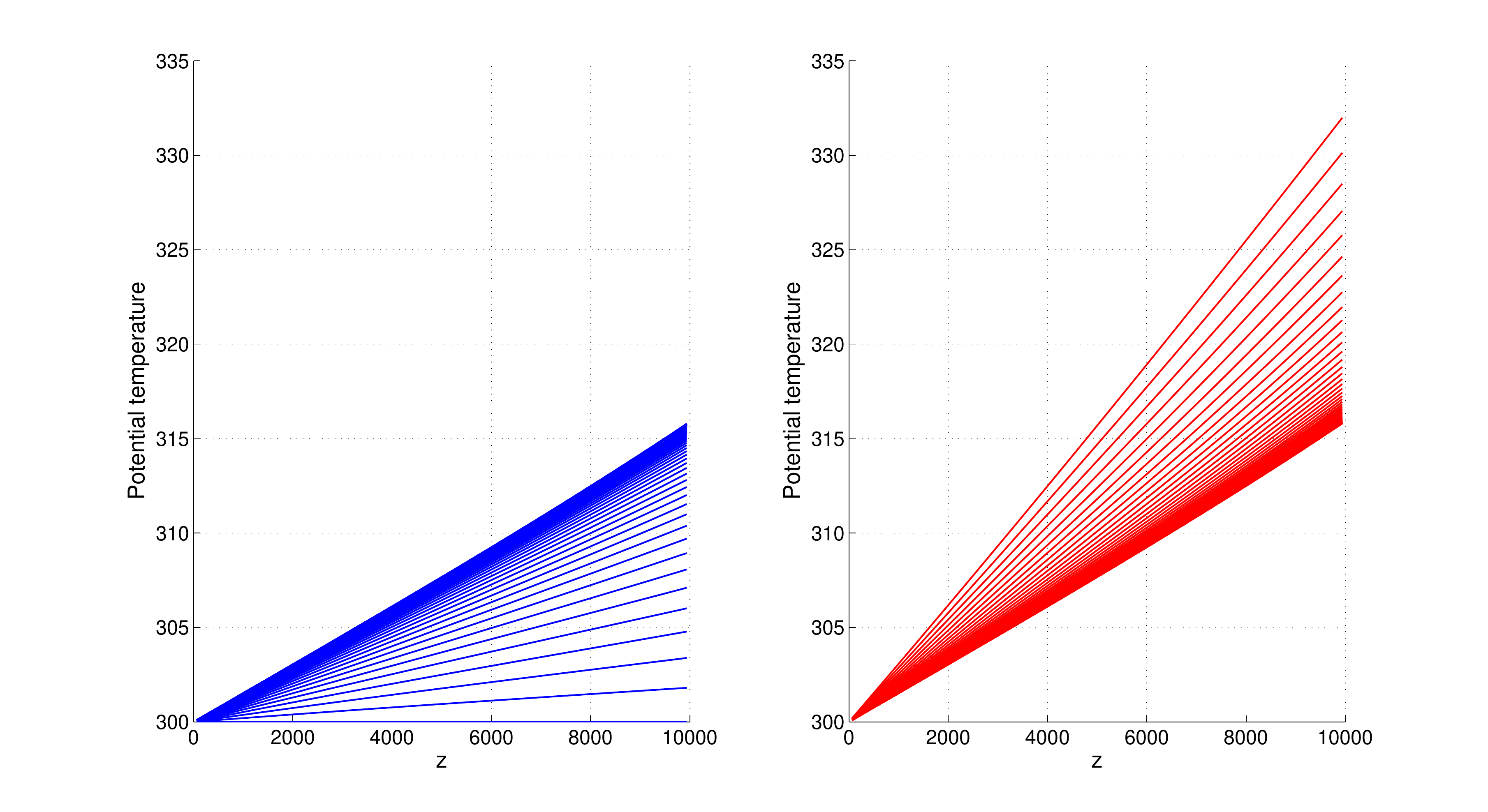}
\label{fig:figure14}
\end{figure}

\begin{figure}[H]
\caption{Temporal evolution of the horizontal velocity in the $y$-direction $v$, vertical section at $x=0$, in layer 1 (left) and layer 2(right).}
\centering
\includegraphics[width=\textwidth]{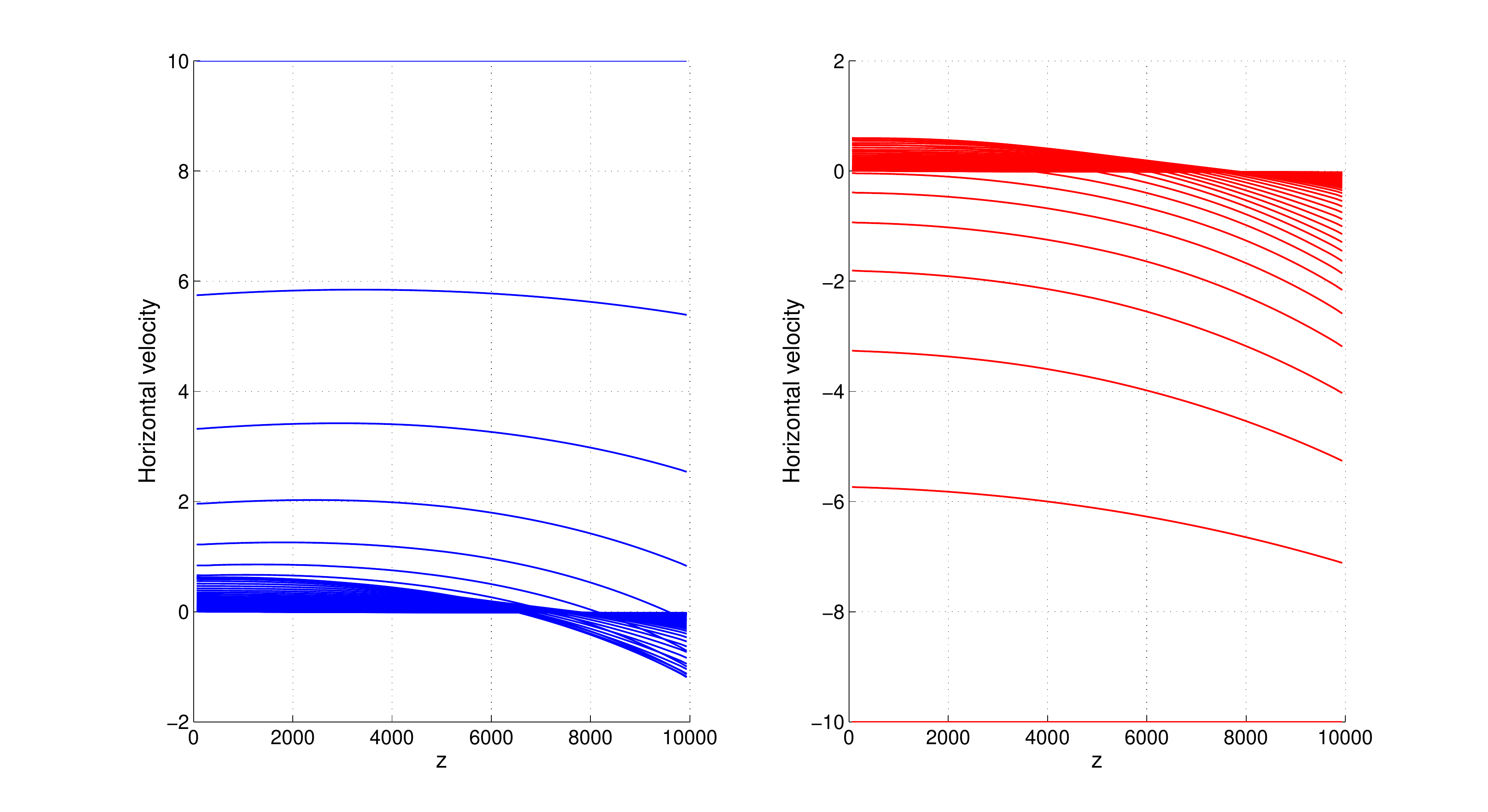}
\label{fig:figure15}
\end{figure}

In figure \ref{fig:figure13} we present the velocity profile for the velocity component in the
x-direction. We see that the profiles are adjusted to some logarithmic profile from
above in layer 1, and from below in layer 2. This is as expected.
The same adjustment is seen for the potential temperature in figure \ref{fig:figure14}, the
adjusted atmosphere becomes stable with a stability which is a ''mean value" between the
two layers.
In figure \ref{fig:figure15} we present the time evolution of the velocity component in the
y-direction, in a cross section at $x=0$. So here we expect to see a decay of an initial
shear, produced by the initial shear in the x-direction. And this is what figure
\ref{fig:figure15} shows, a decay from a positive value in layer 1, and from a negative value
in layer 2.

\begin{figure}[H]
\caption{Total energy of the system for the wind shear test case.}
\centering
\includegraphics[width=0.8\textwidth,height=0.17\textheight]{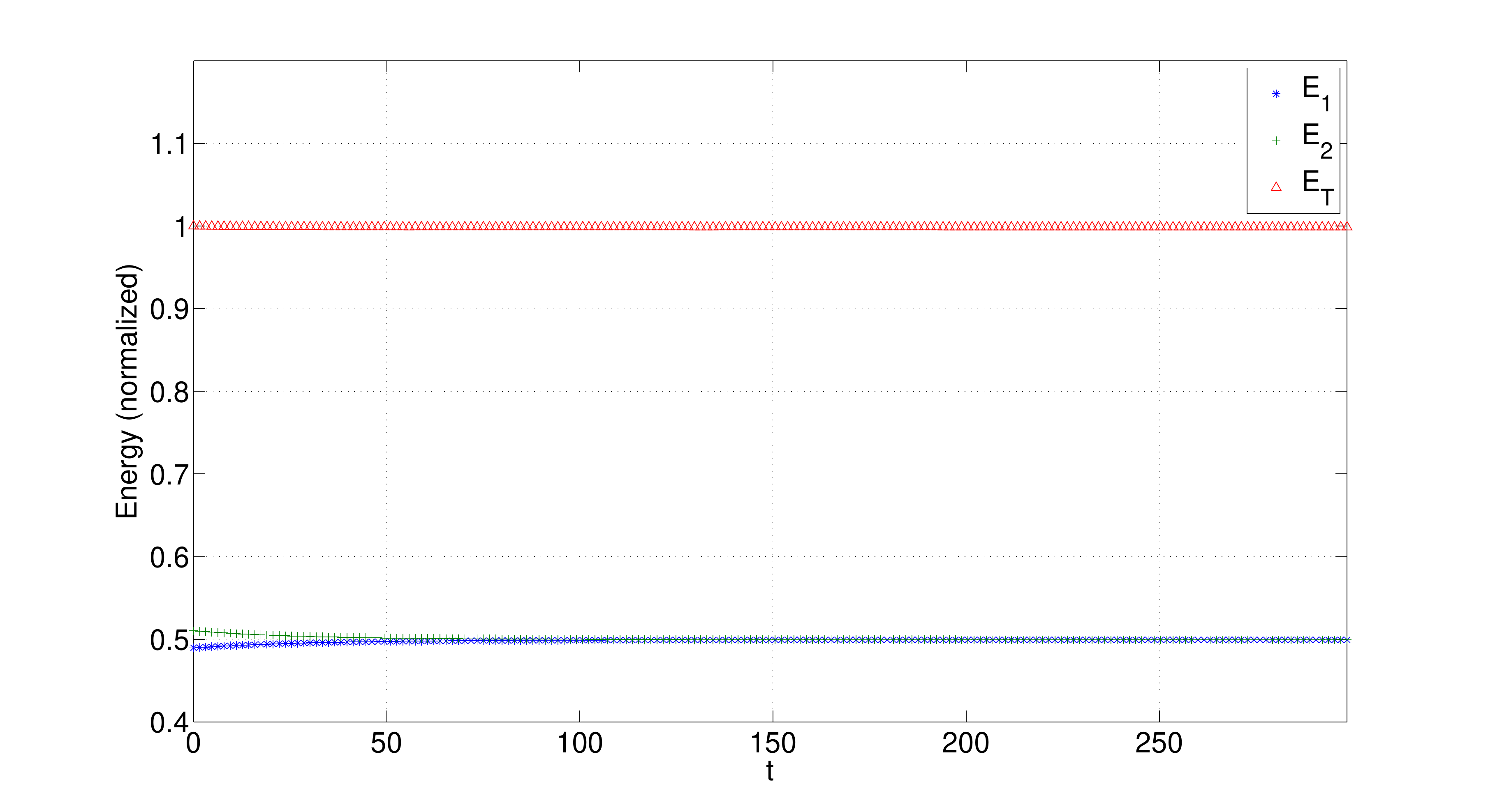}
\includegraphics[width=0.8\textwidth,height=0.17\textheight]{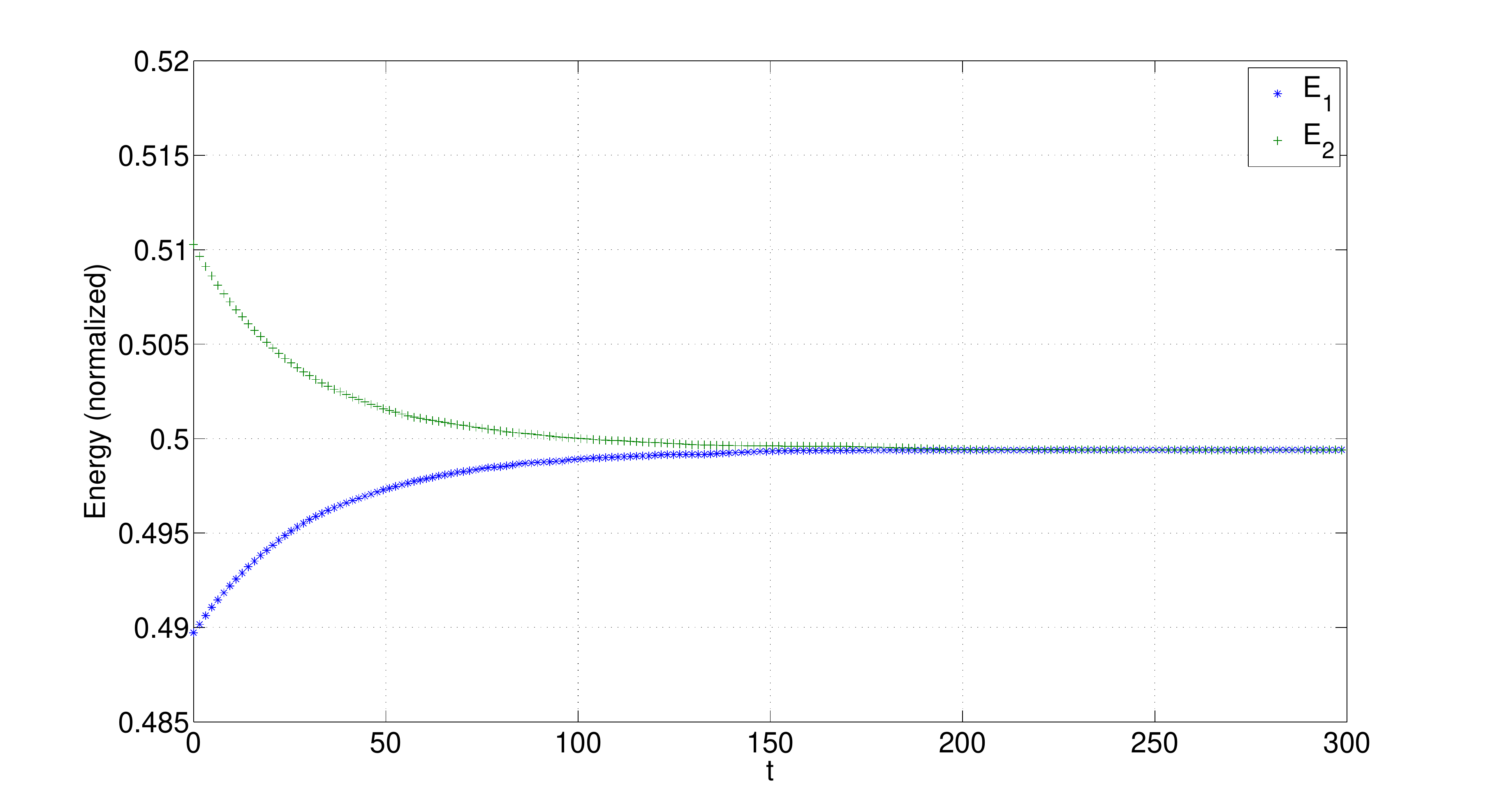}
\label{fig:figure16}
\end{figure}

Note the interesting ''undershoot'' in layer 1, and the ''overshoot'' in layer 2:
The profiles maintains their shape, before the final decay towards zero.
Recall that this test is without explicit friction or parameterizations of the shear
layer, so it is basically thermodynamics ($\theta\to \P\to u$) that controls
the process.
That the adjustment works well is also illustrated in figure \ref{fig:figure16},
which shows the conservation of total energy and the adjustment of the energy in
the two layers.

This test case can be extended to have a rising bubble in one of the layers as well
as other perturbations to simulate front behavior; one can also introduce parameterizations of the shear layer, for example by considering
a simple turbulence model (a Smagorinsky model is one choice).

\subsection{Inertia-gravity waves through a periodic channel}

\begin{figure}[H]
\caption{Horizontal velocity $u$ extrema evolution. Inertia-gravity waves case; combined test run.}
\centering
\includegraphics[width=\textwidth]{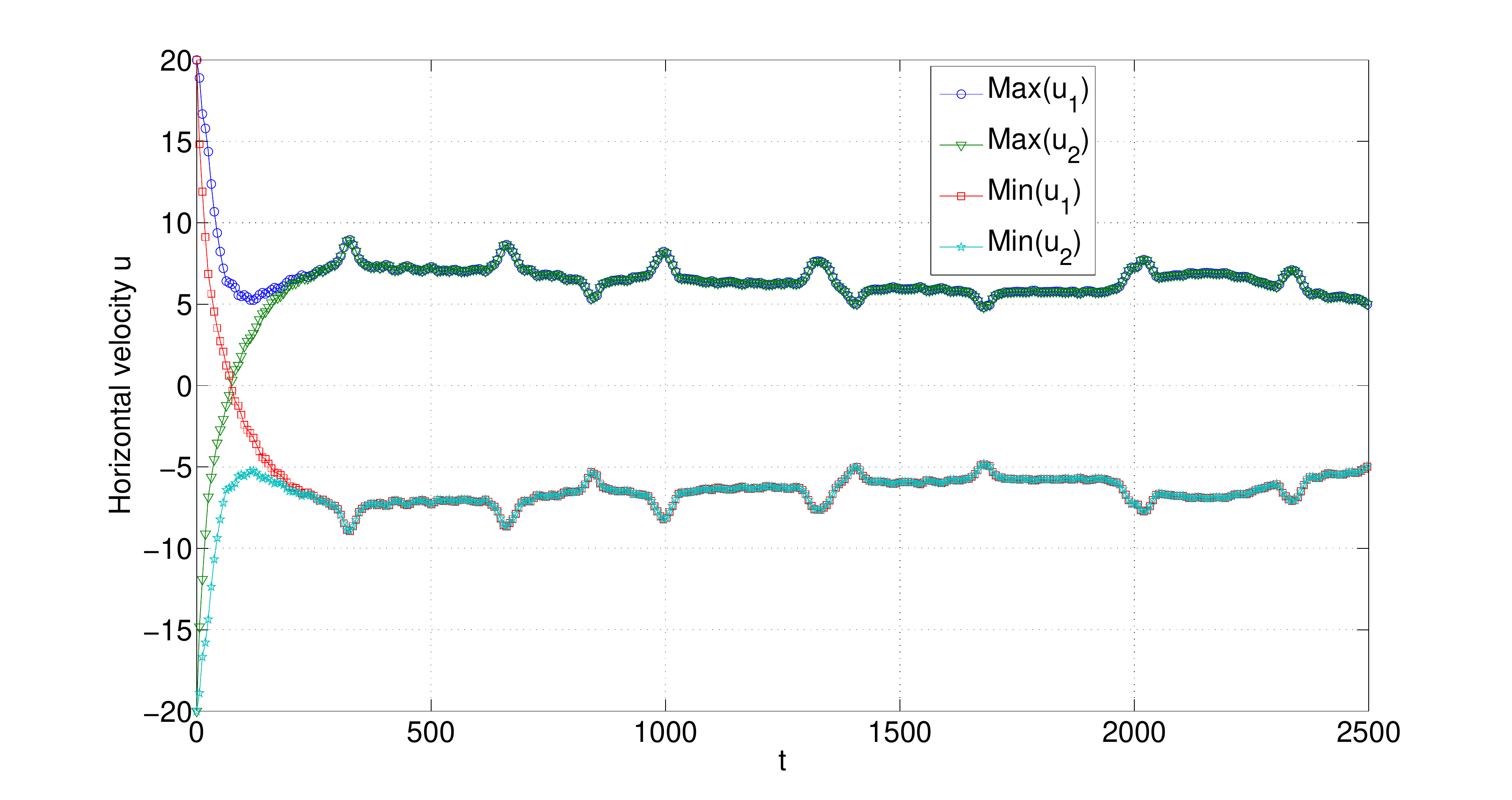}
\label{fig:figure20}
\end{figure}

\begin{figure}[H]
\caption{Inertia-gravity waves. Potential temperature colormaps at $t=0$ and $t=2500$ for separate tests. $\dx=\dz=500[m]$, $600\times 20$ elements.}
\centering
\includegraphics[width=\textwidth, height=0.9\textheight]{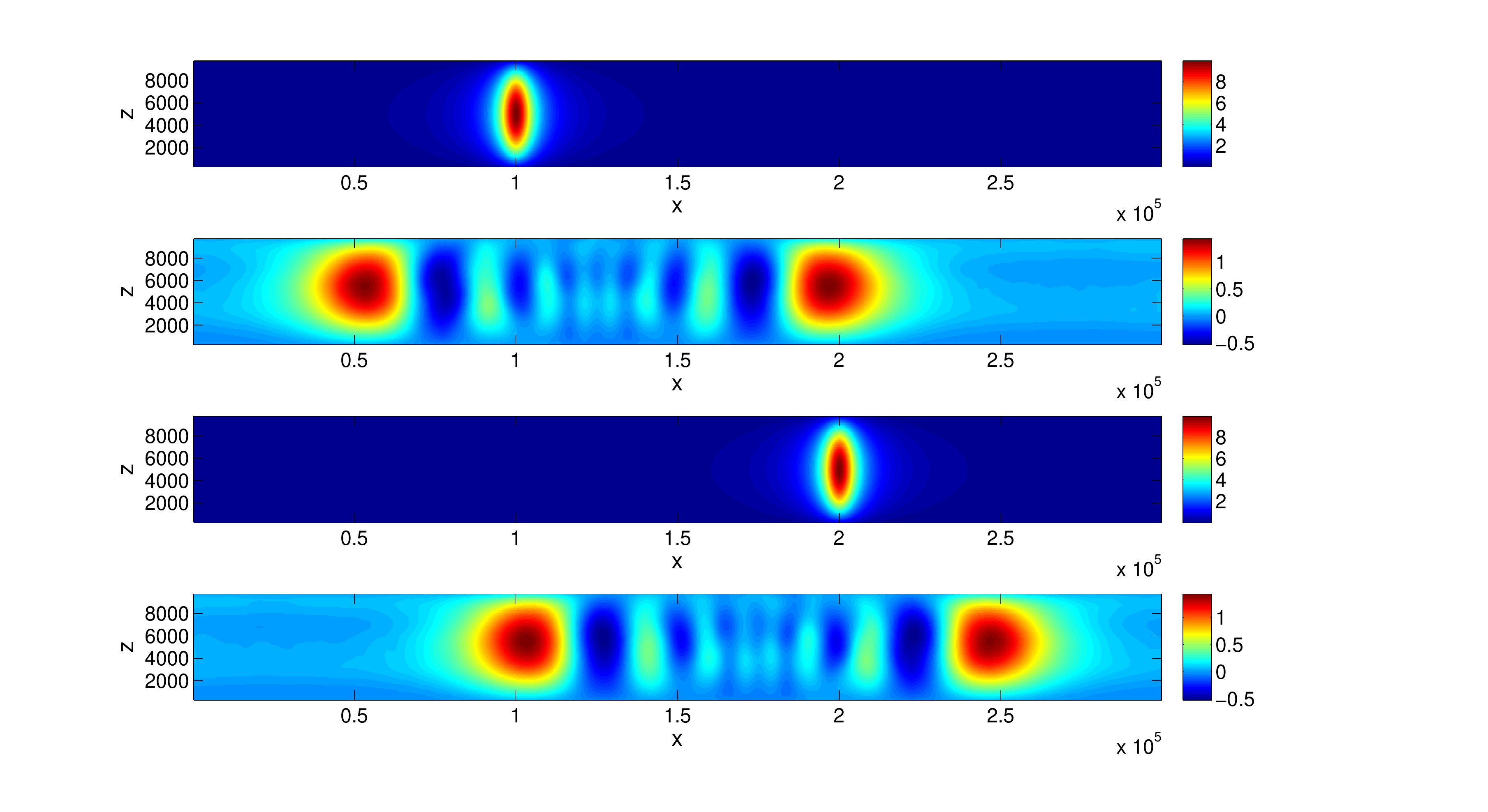}
\label{fig:figure17}
\end{figure}

\begin{figure}[H]
\caption{Inertia-gravity waves. Potential temperature colormaps for the first layer at $t=0$, 1000, 2000 and 2500$[s]$ in the combined test. $\dx=\dz=500[m]$, $600\times 20$ elements.}
\centering
\includegraphics[width=\textwidth, height=0.9\textheight]{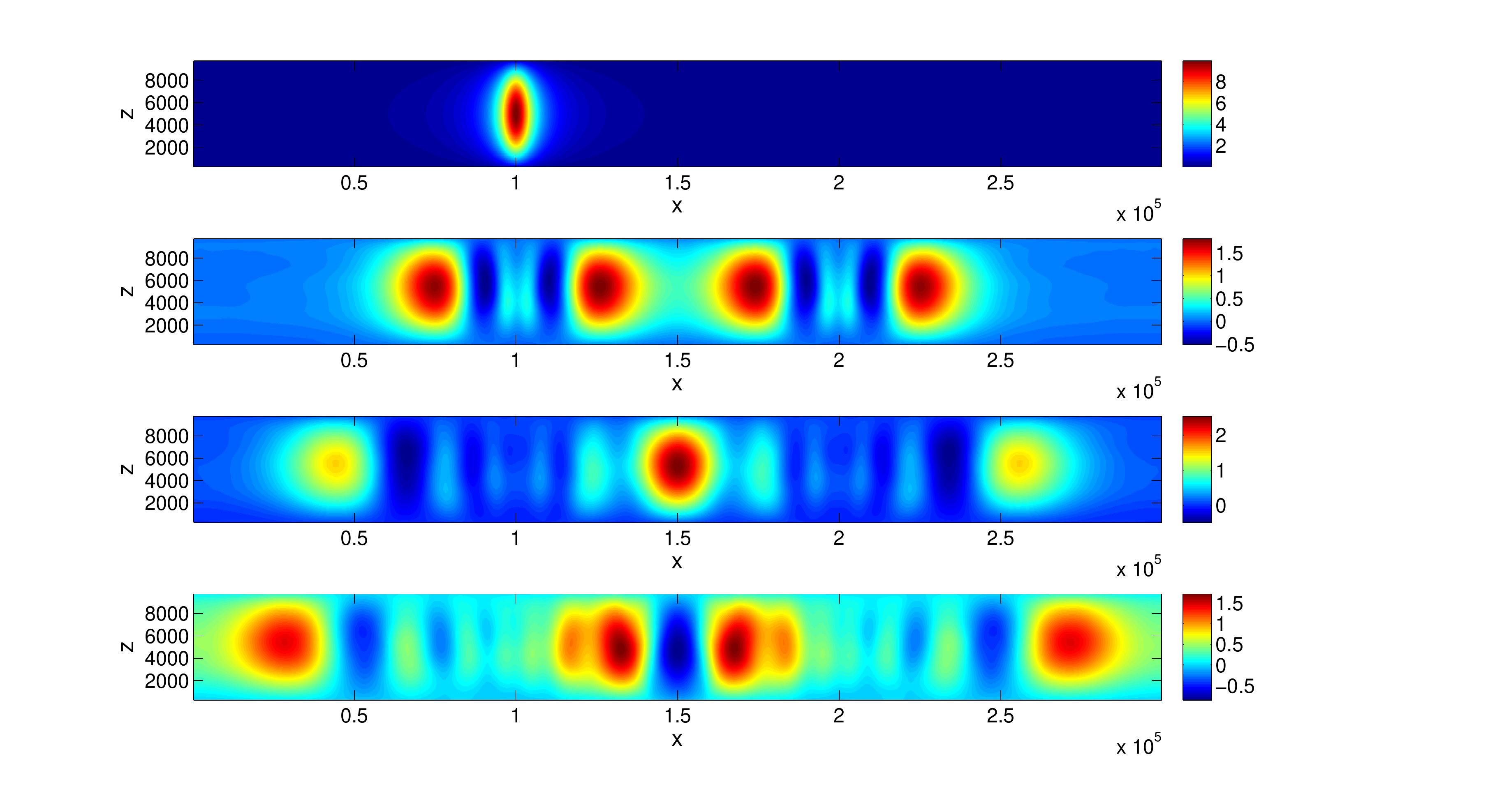}
\label{fig:figure18}
\end{figure}

\begin{figure}[H]
\caption{Inertia-gravity waves. $X-Z$ vector field plot for the first layer at $t=0$, 1000, 2000 and 2500$[s]$ in the combined test. $\dx=\dz=500[m]$, $600\times 20$ elements.}
\centering
\includegraphics[width=\textwidth, height=0.8\textheight]{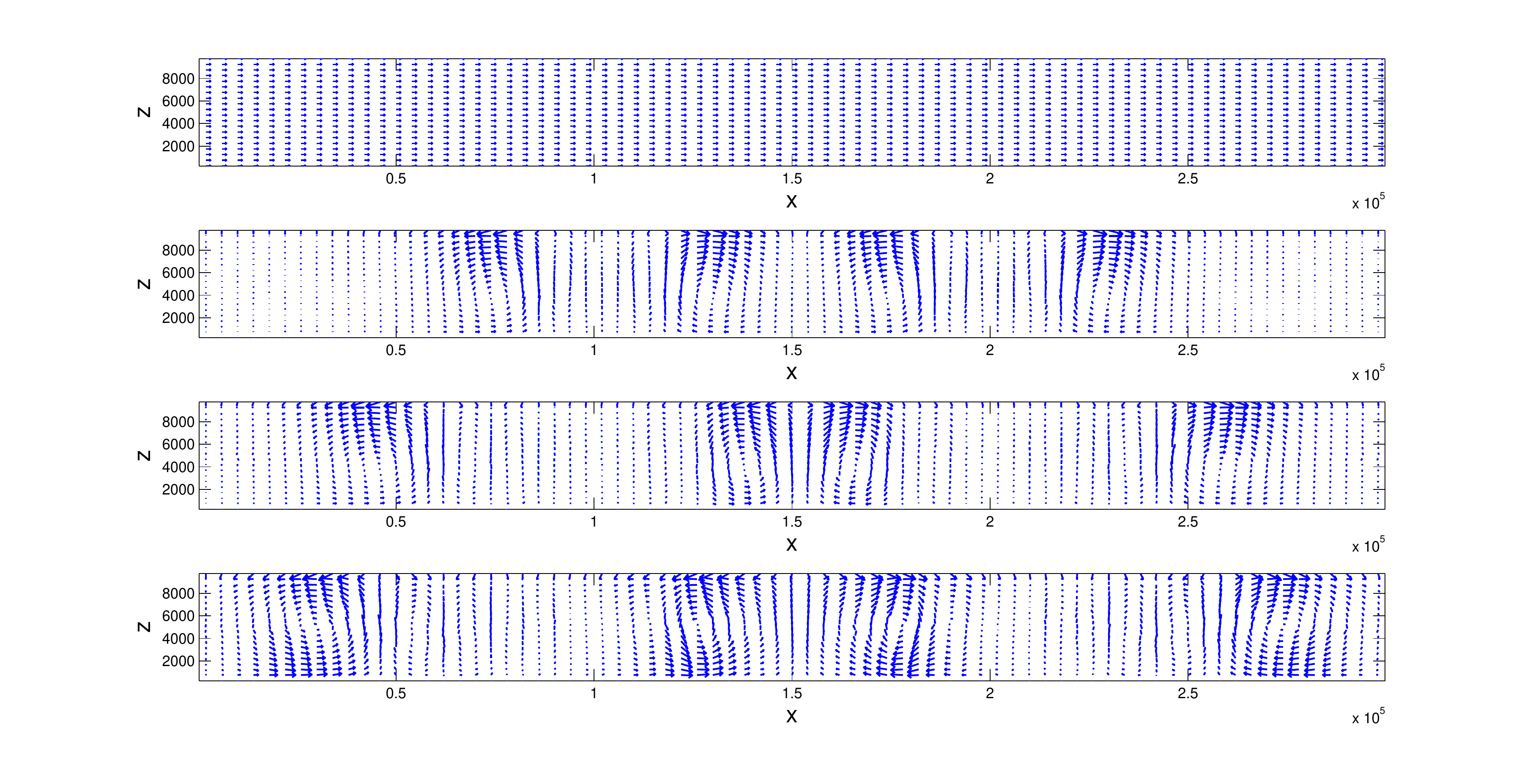}
\label{fig:figure19}
\end{figure}

\begin{figure}[H]
\caption{Total energy of the system for the inertia-gravity waves case; combined test run.}
\centering
\includegraphics[width=0.8\textwidth,height=0.17\textheight]{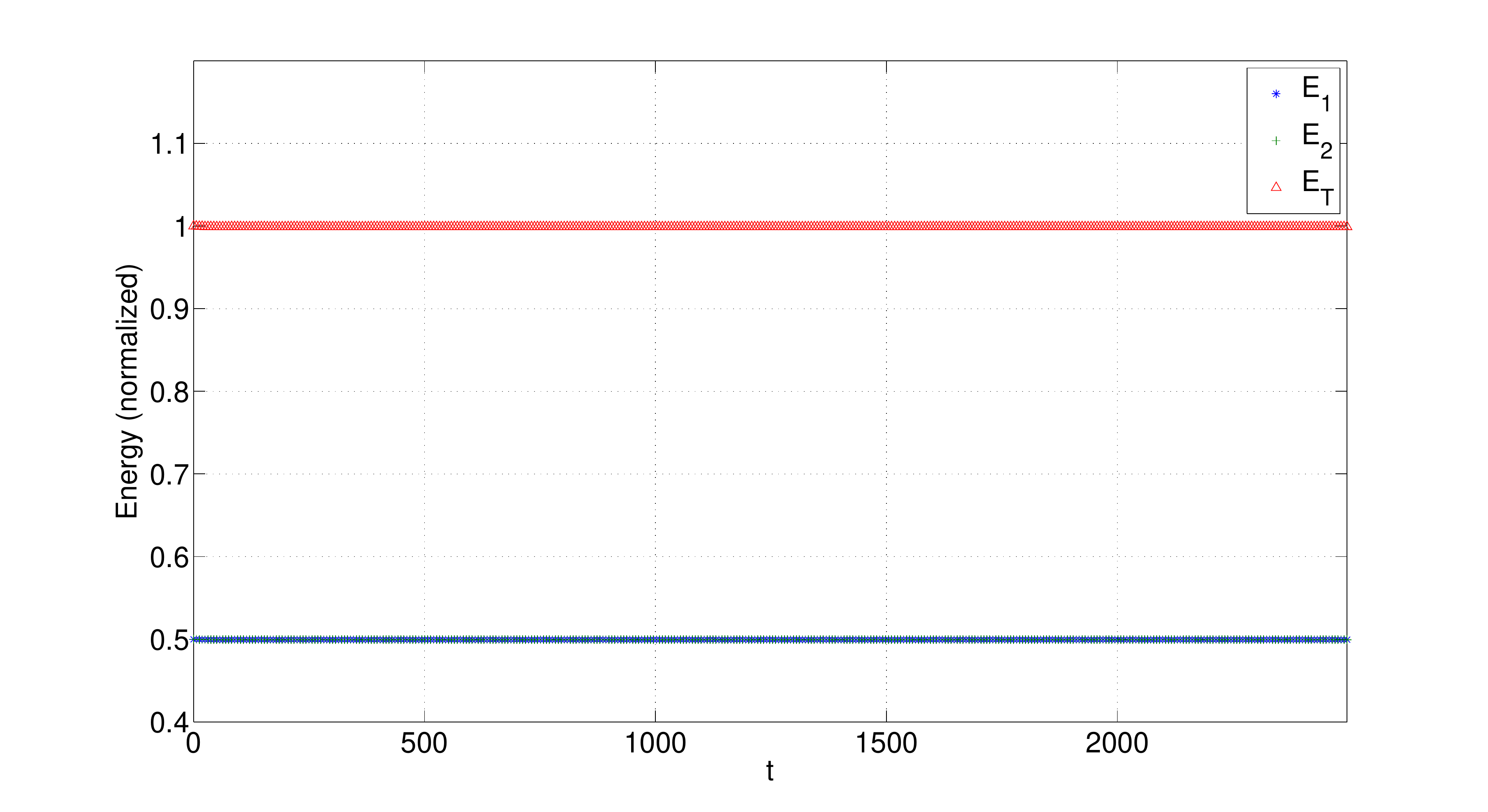}
\includegraphics[width=0.8\textwidth,height=0.17\textheight]{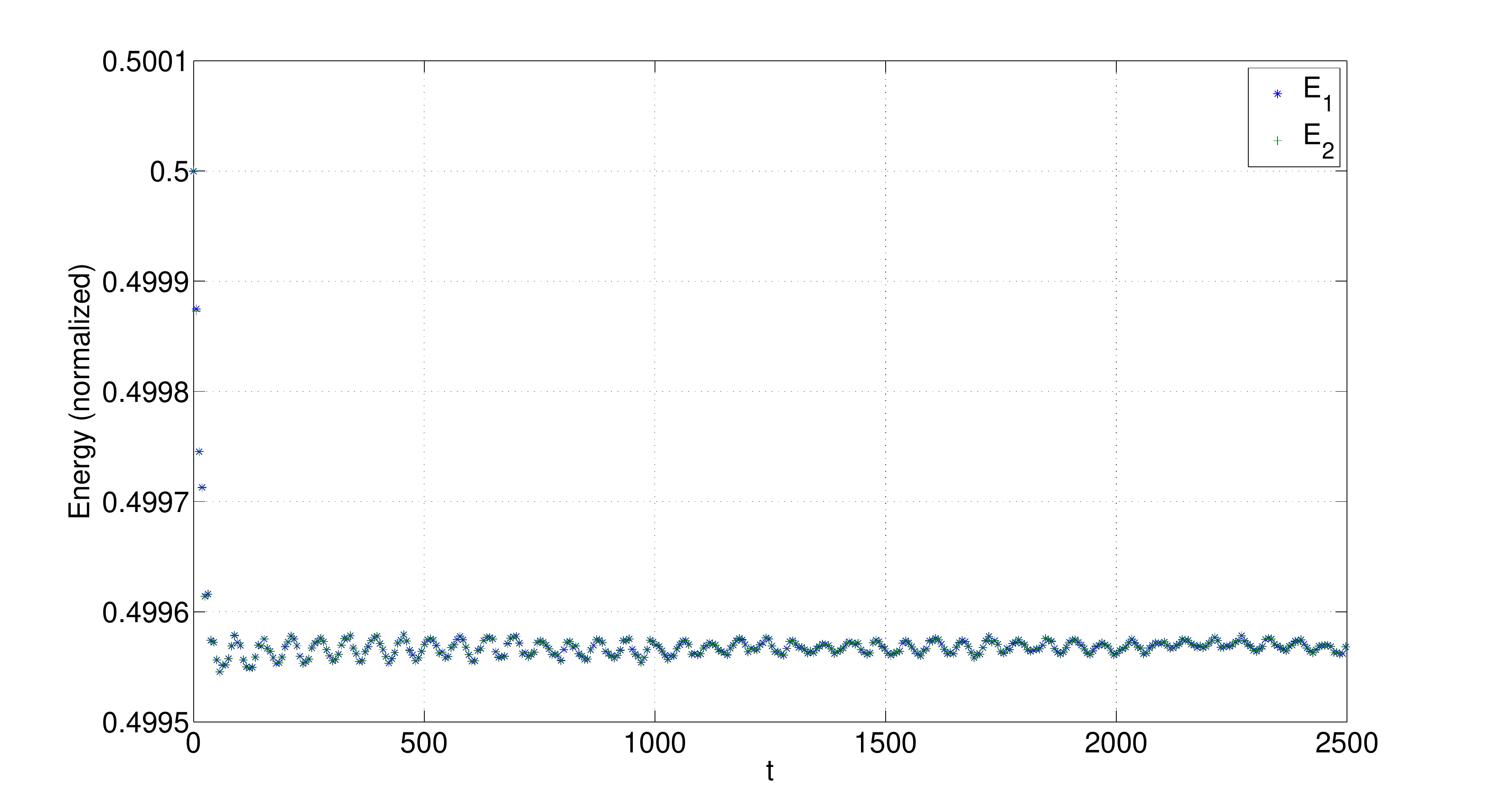}
\label{fig:figuresha}
\end{figure}

The last case that we study is a variation of the inertia-gravity waves test case
proposed by \cite{skamarock}. The purpose of this test is twofold:
first, we qualitatively compare the behavior of our model for the generation of wave motion
in a nonhydrostatic scale, but we also study the nonlinear effects of the coupling.
We set $\Omega=[0,\,300000]\times[-10000,\,10000]\times[0,\,10000]\times[0,\,3000]$,
with $\dx=\dz=500[m]$. In a first model run, the horizontal velocity $u$ is set to a value
of $20 [ms^{-1}]$ for the first layer, while any other velocity variable is set to zero.
Both layers are initialized with the same stable atmosphere profile as in the previous example,
and a temperature perturbation is added to the first layer, in the form,
\begin{equation}
\theta'_1=\theta_p\frac{\sin \left(\frac{\pi z}{L^Z}\right)}{1+\left(\frac{x-x_p}{a}\right)^2},
\end{equation}
with $\theta_p=10[K]$, $x_p=100,000[m]$, and $a_p=5,000[m]$.
In a second model run, we leave the first layer at rest and add a velocity $u=-20[ms^{-1}]$
to the second layer. We remove the temperature perturbation from the first layer
and assign it to the second, only changing the parameter $x_p=200,000$, thus expecting a similar
behavior as in the first model run but in the opposite direction.
A third and final run considers a combination of the previous model runs, i.e.,
one perturbation in every layer, the first layer with a positive horizontal displacement
and the second layer with a negative horizontal displacement.
\vspace{2mm}
Figure \ref{fig:figure17} shows the results of the first and second model runs, initial condition
and result at $t=2500$, first run in the two upper panels and second run in the two lower
panels. These results are in very good agreement with the results in \cite{skamarock}
and other test results in the literature.
In figure \ref{fig:figure18} we present the results of the third model run.
We see immediately that these results are not superpositions of the results in figure
\ref{fig:figure17}, so we see the effects of a nonlinear interaction between the two inertia-gravity
waves.
In figure \ref{fig:figure19} we present the velocity fields for the nonlinear interaction
case. We see time-dependent fan-like structures typical of nonlinear wave interaction
and almost perfect symmetry of the fields. Figure \ref{fig:figure20} shows the evolution of extremal values for the horizontal velocity in every layer. It can be seen that both layers tend to an equilibrium, while an interaction between maximal and minimal values occurs.
It is not easy to explain all the details in the results, for example the
development from $t=2000$ to $t=2500$, since the interaction is purely nonlinear.
However, it is possible to see effects of focussing, which is a purely nonlinear wave
phenomenon.
In the solutions the symmetry is well preserved indicating
that the is performing in a good manner; conservation of energy is shown in figure \ref{fig:figuresha} .

\section{Summary and concluding remarks}

We have derived a set of dimensionally-reduced fully nonlinear primitive equations to model atmospheric dynamics.
Our goal was to show that the 2.5D model is  able to reproduce well-known atmospheric phenomena properly.
A DG approach was used for the dimensional reduction, whereas a
WENO-TVD method was used for discretization. We think that finite volume methods are the most promising discretization method
for atmospheric models since it ensures discrete conservation. If one uses,
as we have done, equations in conservative form with conserved variables,
the discretized equations should model the conserved quantities in the atmosphere
with high accuracy.
The numerical experiments with extensively used test cases for atmospheric model cores
have shown that the dicretization is accurate, stable, energy-conserving
and do not produce spurious oscillations. No artificial diffusion or any other
form of smoothing is used.
We think that this is a strong indication that our set of equations are well-posed,
and that a higher order discretization is indeed valuable.
Our approach can be extended in several directions: we could use higher order elements
in the dimensional reduction and/or even higher order finite volume methods.
The purpose of this would be to investigate the accuracy versus the computational
costs. This can easily be done for well-known test cases where we have results
to compare with, but we would also like to use more complicated tests cases from
real weather situations. The latter test cases could be used to investigate
how one can model certain weather phenomena where one, from a physical viewpoint,
have a certain insight in the behavior.
There is also an option of using other methods than the DG method
for the dimensional reduction, but our opinion is that this method is good natural choice,
because it poses no difficulty  for handling nonlinear equations.
Also on the computational side, the actual computational costs has not been explicitly
addressed in this paper. However, it is worth to mention that every block of the code can be completely vectorized,
while the dimensionally-reduced equations have a natural form of parallelism by parallelizing
over the layers. This is an important topic since the existing dynamical cores of atmospheric models
which we would like to compare with in term of computational cost, are parallelized.

\bibliographystyle{model3-num-names}
\bibliography{bibtesis}

\end{document}